\documentclass[3p]{elsarticle}

\usepackage{tikz}
\newcommand*{\DashedArrow}[1][]{\mathbin{\tikz [baseline=-0.25ex,-latex, dashed,#1] \draw [#1] (0pt,0.5ex) -- (1.3em,0.5ex);}}%
\usepackage{comment}

\usepackage{amsfonts}
\usepackage{graphics}
\usepackage{color}
\usepackage{amssymb}
\usepackage{verbatim}
\usepackage{amsmath}

\usepackage{amssymb}

\usepackage{amscd}

\usepackage{euscript}

\newtheorem{Theorem}{Theorem}

\newtheorem{Corollary}[Theorem]{Corollary}

\newtheorem{Example}[Theorem]{Example}

\newtheorem{Definition}[Theorem]{Definition}
\newtheorem{remark}[Theorem]{Remark}

\newtheorem{Lemma}[Theorem]{Lemma}

\newtheorem{Proposition}[Theorem]{Proposition}

\newtheorem{Fundamental Theorem}{Fundamental Theorem}

\newenvironment{Proof}[1][Proof]{\textbf{#1.} }{\ \rule{0.5em}{0.5em}}

\usepackage[all]{xy}

\xyoption{v2}
\xyoption{all}
\xyoption{2cell}

\def \E {\mathbb{E}}

\def \Hom {{\rm hom}}

\def \D {\EuScript{D}_*}

\def \fo {\textrm{ for each }}
\def \an {\textrm{ and }}

\def \wh {\textrm{ where }}
\def \ad {\mathrm{ad}}

\def \T {\EuScript{C}{\rm omp}}

\def \fo {\textrm{ for each }}
\def \an {\textrm{ and }}

\def \wh {\textrm{ where }}
\def \A {\mathcal{A}}
\def \B {\mathcal{B}}
\def \T {\EuScript{T}}
\def \id {\mathrm{id}}

\def \d {\partial}
\def \de {\delta}
\def \t {\triangleright}
\def \tn {\otimes}
\def \ra {\xrightarrow}
\def \w {\omega}
\def \Z {\mathbb{Z}}
\def \im {\mathrm{im}}

\def \G {\mathcal{G}}

\def \f {\phi}
\def \p {\psi}
\def \a {\alpha}
\def \b {\beta}

\def \de {\delta}

\def \HOM {\mathrm{HOM}}

\def \K {\mathcal{K}}

\def \P {\EuScript{P}_*}

\def \Tet {\mathcal{TET}}

\def \SG {\mathcal{SG}}

\def \Gr { {\rm Gr} }
\def \Pr { {\rm Pr} }
\def \FG {{\mathcal{F}^{\scriptscriptstyle{\rm group}}}}

\def \e {\emptyset}

\newcommand{\sqspace}[9] {&#3\ar[rr]^{#8}&& #4 \\ & &#9\ar@{.>}
[r]^{\quad  d_{0}}\ar@{.>}[u]|{d_{1}}\ar@{.>}[l]_{d_{3}\quad }\ar@{.>}[d]|{d_{2}} & \\ &#1\ar[uu]^{#7} \ar[rr]_{#5} &&#2\ar[uu]_{#6}  }

\newcommand{\trispace}[7] {&&&&&#3 \\ &&&&\\ &&& &{#7}\ar@{.>}
[r]^{\quad  d_{0}}\ar@{.>}[l]^{d_{1}} \ar@{.>}[d]|{d_{2}} & \\ &#1\ar[uuurrrr]^{#6} \ar[rrrr]_{#4}& &&&#2\ar[uuu]_{#5}  }

\newcommand{\diskspace}[5] {&&&&&&\\&#1  \ar@/^3pc/[rrrr]^{{#4}}\ar@/_3pc/[rrrr]_{{#3}}  & &#5 \ar@{.>}[u]|{d_{1}} \ar@{.>}[d]|{d_{2}}& & #2 \\&&&&&&}

\newcommand{\pathspace}{&\A\ar[rrr]^{\stackrel{\Pr_0^\A}{   \DashedArrow[<-,densely dotted    ]   } \,\,\P(\A)\,\, \stackrel{\Pr_1^\A}  {   \DashedArrow[->,densely dotted    ] }}&&&{\A} } 

\newcommand{\pu}[2] {#1 \ra{#2} #1\d(#2)}

\newcommand{\pt}[3] {#1 \ra{(#2,#3)} #1#2\de(#3)}

\newcommand{\squ}[5] { &#1\d(#3)\ar[rr]^{(\d(#3)^{-1} \t #2)\,\, #4\quad \quad } & &#1\, \d(#2)\, \d(#3)\,\d(#4)\\
&#1   \ar@/^2pc/ @{{}{ }{}}[rr]|{#5} \ar[u]^{#3}\ar[rr]_{#2} & & #1\d(#2) \ar[u]_{#3\,#4\, \de(#5)}
}

\newcommand{\sqt}[8] {&#1 #4\de(#7) \ar[rrrrrrrrrr]^{\Big(   ( \d(#4^{-1} )\t #2) #5, \quad  \big( ( \d(#4^{-1} )\t #2) #5\big)^{-1} \t' l^{-1}\,\,\,\, #5^{-1} \t'\big (\{#4^{-1}, #2^{-1} \}^{-1}\,\, #4^{-1} \t' #3 \big )\,\, #6\,\,  l\,\, #8 \Big)\quad \quad \quad \quad} &&&&&&&&& & #1#2\de(#3) #4#5\de(#6)\de(#7) \de(#8)\\
&#1 \ar[u]^{(#4,#7)}\ar[rrrrrrrrrr]_{(#2,#3)} & &&&&&&&&& #1#2\de(#3) \ar[u]_{(#4\,#5 \de(#6),#7#8)}
}

\newcommand{\trit}[8] {& & & #1#2\de(#3) #5\de(#6) \de(#8)\\ & & &\\
&#1  \ar[rruu]^{\big(    #2 #5, #5^{-1} \t' #3 \,\, #6 #8 \big) \quad }\ar[rr]_{(#2,#3)} & & #1#2\de(#3) \ar[uu]_{(#5 \de(#6),#8)}
}


\newcommand{\tritl}[8] {& & &  \psi(e)\,s(\d(e))\,\de(t(e))\,s'(\d(e))\,\de(t'(e))\\ & & &\\
&#1  \ar[rruu]^{\Big(    s(\d(e))\,s'(\d(e))\,\big(\de(\w^{(s,s')}(\d(e))\big)^{-1}, \w^{(s,s')}(\d(e))\,s'(\d(e))\t' t(e) \,\,t'(e) \Big)\quad \quad \quad \quad\quad \quad \quad }\ar[rr]_{\big(#2,#3\big)} & & #1\,#2\,\de(#3) \ar[uu]_{\big(s'(\d(e)),t'(e)\big)}
}


\newcommand{\trixu}[5] {& & & #1 \d(#2) \d(#4)\\ & & &#5 \quad \quad \quad \quad \quad \\
&#1  \ar[rruu]^{#2#4}\ar[rr]_{#2} & & #1 \d(#2) \ar[uu]_{#4 \de(#5)} 
}


\newcommand{\trixuvb}[5] {& & & #1 \d(#2) \d(#4)\\ & & &#5 \quad \quad \quad \quad \quad \\
&#1  \ar[rruu]^{#2#4}\ar[rr]_{#2} & & #1 \d(#2) \ar[uu]_{#4 } 
}


\newcommand{\trixun}[5] {& & & #1 \\ & & &#5 \quad \quad \quad \quad \quad \\
&#1  \ar[rruu]^{1}\ar[rr]_{#2} & & #1 \d(#2) \ar[uu]_{\overline{s}(g)} 
}


\newcommand{\trixumod}[5] {& & & #1 \d(#2) \d(#4)\\ \\
&#1 \ar @/^2pc/ @{{}{ }{}} [rrrr]|{\boxed{#5}} \ar[rruu]^{#2#4}\ar[rrrr]_{#2} & &&& #1 \d(#2) \ar[uull]_{ #4 \de(#5)} 
}

\makeatother
\newcommand{\trixumodb}[5] {& & & #1 \d(#2) \\ \\
&#1 \ar @/^2pc/ @{{}{ }{}} [rrrr]|{\boxed{#5}} \ar[rruu]^{#2#4}\ar[rrrr]_{#2} & &&& #1 \d(#2) \ar[uull]_{ 1} 
}
\makeatother
\newcommand{\trixumodf}[5] {& & & \phi''(g) \d(s''(g)) \d(s(g))\\ \\
&#1 \ar @/^2pc/ @{{}{ }{}} [rrrr]|{\boxed{#5}} \ar[rruu]^{#2#4\quad }\ar[rrrr]_{#2} & &&& #1 \d(#2) \ar[uull]_{ s'(g)} 
}

\makeatother
\newcommand{\trixumodg}[5] {& & &  \phi''(g) \d(s''(g)) \d(s(g))  \\ \\
&#1 \ar @/^2pc/ @{{}{ }{}} [rrrr]|{\boxed{#5}} \ar[rruu]^{#2#4\quad \quad }\ar[rrrr]_{#2} & &&& #1 \d(#2) \ar[uull]_{s(g) k(g)} 
}

\makeatother
\newcommand{\trixumods}[5] {& & & \phi(g)\, \d(s(g))\, \d(s''(g))\\ \\
&#1 \ar @/^2pc/ @{{}{ }{}} [rrrr]|{\boxed{#5}} \ar[rruu]^{#2#4\quad \quad \quad }\ar[rrrr]_{#2} & &&& #1 \d(#2) \ar[uull]_{ s''(g)} 
}

\makeatother
\newcommand{\trixumoda}[5] {& & &   \phi(g)\, \d(s(g))\, \d(s''(g))\\ \\
&#1 \ar @/^2pc/ @{{}{ }{}} [rrrr]|{\boxed{#5}} \ar[rruu]^{#2#4\quad \quad \quad}\ar[rrrr]_{#2} & &&& #1 \d(#2) \ar[uull]_{\quad \quad \de(k(g))\,s''(g)\, \de\big(s''(g)^{-1} \t' k(g)^{-1}\big)} 
}

\makeatother
\newcommand{\trixumodn}[5] {& & & #1 \d(#2) \\ \\
&#1 \ar @/^2pc/ @{{}{ }{}} [rrrr]_{\boxed{#5}} \ar[rruu]^{#2}\ar[rrrr]_{#2} & &&& #1 \d(#2) \ar[uull]_{1} 
}

\makeatother
\newcommand{\trixumodp}[5] {& & & #1 \d(#2) \\ \\
&#1 \ar @/^2pc/ @{{}{ }{}} [rrrr]|{\boxed{#5}} \ar[rruu]^{#2}\ar[rrrr]_{#2} & &&& #1 \d(#2) \ar[uull]_{#4 } 
}

\makeatother
\newcommand{\trixumodd}[5] {& & & \phi(g)\, \d(s(g))\, \d(s'(g)\,\d(s'(g)) \\ \\
&#1 \ar @/^1pc/ @{{}{ }{}} [rrrr]|{\boxed{#5}} \ar[rruu]^{(s \tn s' \tn s'')(g)}\ar[rrrr]_{#2} & &&& #1 \d(#2) \ar[uull]_{(s' \tn s'')(g)} 
}


\newcommand{\trixumodddd}[5] {& & & #1 \d(#2) \d(#4)\\ \\
&#1 \ar @/^2pc/ @{{}{ }{}} [rrrr]|{\boxed{#5}} \ar[rruu]^{#2#4}\ar[rrrr]_{#2} & &&& #1 \d(#2) \ar[uull]_{#4 } 
}

\newcommand{\trixumoddd}[5] {& & & \phi(g)\, \d(s(g))\, \d(s'(g)\,\d(s'(g)) \\ \\
&#1 \ar @/^1pc/ @{{}{ }{}} [rrrr]^{\boxed{\w^{(s,s' \tn s'')}(g)}} \ar[rruu]^{(s \tn s' \tn s'')(g)}\ar[rrrr]_{#2} & &&& #1 \d(#2) \ar[uull]_{(s' \tn s'')(g)} 
}

\newcommand{\trixuu}[5] {& & & #1 \d(#2) \d(s'(g))\\ & & &#5 \quad \quad \quad \quad \quad \\
&#1  \ar[rruu]^{(s \tn s')(g)}\ar[rr]_{#2} & & #1 \d(#2) \ar[uu]_{s'(g)} 
}


\newcommand{\trixul}[5] {& & & \phi(g) \, \d(s(g))\,\d(s'(g)) \\ & & &#5 \quad \quad \quad \quad \quad \\
&#1  \ar[rruu]^{#2#4}\ar[rr]_{#2} & & #1 \d(#2) \ar[uu]_{s'(g)} 
}


\newcommand{\dt} {
&a  \ar@/^2pc/[rr]^{\big(    e \de(k')^{-1},  k 'k   \big)}\ar@/_2pc/[rr]_{(e,k)} & & ae\de(k) 
}


\newcommand{\du} {
&g   \ar@/^2pc/[rr]^{x\de(k)^{-1}}\ar@/_2pc/[rr]_{x} & k  & g \d(x) }

\newcommand{\te}[9] { & & & #1 \d(#2) \d(#4) \d(#6)  \\ \\ & && #1 \d(#2) \d(#4)\ar[uu]|{#6 \de(#7) \de(#9)} \\ & #1\ar[uuurr]^{#2 #4 #6} \ar[urr]_{#2 #4 } \ar[rrrr]_{#2  }\ar @/^1pc/ @{{}{ }{}} [rrrr]|{\boxed{#5}}  \ar@/_2.5pc/ @{{}{ }{}} [rruuu]^{\boxed{#7 #9}}  & && & #1 \d(#2) \ar[llu]^{#4 \de(#5)} \ar[lluuu]_{#4 \de(#5)#6 \de(#7)} \ar@/^1pc/ @{{}{ }{}} [lluuu]|{ \boxed{#9}}   }

\newcommand{\tew}[9] { & & & #1 \d(#2) \d(#4)   \\ \\ & && #1 \d(#2) \d(#4)\ar[uu]_{1} \\ & #1\ar[uuurr]^{#2 #4 #6} \ar[urr]_{#2 #4 } \ar[rrrr]_{#2  }\ar @/^1pc/ @{{}{ }{}} [rrrr]|{\boxed{#5}}  \ar@/_0pc/ @{{}{ }{}} [rruuu]_>>>>>>>>>>>>{\boxed{#9}}  & && & #1 \d(#2) \ar[llu]^{#4} \ar[lluuu]_{#4 #6 } \ar@/^0pc/ @{{}{ }{}} [lluuu]^>>>>>>>>>>>>>{ \boxed{#9}}   }

\newcommand{\teq}[9] { & & & \phi''(g) \d(s''(g)) \d(s(g)) \\ \\ & && \phi''(g) \d(s''(g)) \d(s(g))\ar[uu]_{1} \\ & #1\ar[uuurr]^>>>>>>>>>>>>>>>>>>>>>>>{#2 #4 #6 \quad \quad } \ar[urr]|{(s'' \tn s)(g) } \ar[rrrr]_{#2  }\ar @/^1pc/ @{{}{ }{}} [rrrr]|{\boxed{#5}}  \ar@/_0pc/ @{{}{ }{}} [rruuu]_>>>>>>>>>>>>>>>{\boxed{\big((s'' \tn k)(g)\big)^{-1}}}  & && & #1 \d(#2) \ar[llu]|{s(g)} \ar[lluuu]_{s(g)k(g)} \ar@/^0pc/ @{{}{ }{}} [lluuu] ^>>>>>>>>>>>>>>>{ \boxed{#9}}   }


\newcommand{\tee}[9] { & & & \phi(g) \d(s(g)) \d(s''(g)) \\ \\ & && \phi(g) \d(s(g))\ar[uu]|{s''(g)} \\ & #1\ar[uuurr]^{#2 #4 #6} \ar[urr]_{#2 #4 } \ar[rrrr]_{#2  }\ar @/^1pc/ @{{}{ }{}} [rrrr]|{\boxed{#5}}  \ar@/_3.5pc/ @{{}{ }{}} [rruuu]^{\boxed{#7 }}  & && & #1 \d(#2) \ar[llu]^{1} \ar[lluuu]_{s''(g)} \ar@/^1pc/ @{{}{ }{}} [lluuu]|{ \boxed{#9}}   }

\newcommand{\ter}[9] { & & & \phi(b)\,\d(s(b))\,\d(s''(b)) \\ \\ & && #1 \d(#2) \d(#4)\ar[uu]|{s''(b)} \\ & #1\ar[uuurr]^{#2 #4 #6} \ar[urr]_{#2 #4 } \ar[rrrr]_{#2  }\ar @/^1pc/ @{{}{ }{}} [rrrr]|{\boxed{#5}}  \ar@/_1pc/ @{{}{ }{}} [rruuu]|<<<<<<<<<<<<<<<<<<<{\boxed{1}}  & && & #1 \d(#2) \ar[llu]^{1} \ar[lluuu]_{s''(b)} \ar@/^1pc/ @{{}{ }{}} [lluuu]|{ \boxed{#9}}   }

\newcommand{\tei}[9] { & & & \phi(g)\,\d(s(g))\, \d(s''(g)) \\ \\ & && \phi(g)\,\d(s(g))\, \d(s''(g)) \ar[uu]_{1} \\ & #1\ar[uuurr]^<<<<<<<<<<<<<<<<{(s \tn s'')(g)\,\ #6\quad } \ar[urr]_{(s \tn s'')(g) } \ar[rrrr]_{#2  }\ar @/^1pc/ @{{}{ }{}} [rrrr]|{\boxed{#5}}  \ar@{{}{ }{}} [rruuu]_>>>>>>>>>>>>>>>>>>>{\boxed{#7}}  & && & #1 \d(#2) \ar[llu]^{s''(g)} \ar[lluuu]_{s''(g)} \ar@/^1pc/ @{{}{ }{}} [lluuu]|{ \boxed{#9}}   }

\newcommand{\teo}[9] { & & & #1 \,\d(#2)\, \d(#4) \\ \\ & && #1\, \d(#2)\, \d(#4)\ar[uu]_{1} \\ & #1\ar[uuurr]^<<<<<<<<<<<<<<<<<<{#2 #4 #6\quad } \ar[urr]_{#2 #4 } \ar[rrrr]_{#2  }\ar @/^1pc/ @{{}{ }{}} [rrrr]|{\boxed{#5}}  \ar@{{}{ }{}} [rruuu]_>>>>>>>>>>>>>>>>>>>>{\boxed{#7}\quad\quad\quad}  & && & #1 \d(#2) \ar[llu]^{#4 } \ar[lluuu]_{#4} \ar@/^1pc/ @{{}{ }{}} [lluuu]|{ \boxed{#9}}   }

\newcommand{\tep}[9] { & & & #1\, \d(#6)  \\ \\ & && #1 \ar[uu]|{#6 } \\ & #1\ar[uuurr]^{#2 } \ar[urr]_{ 1} \ar[rrrr]_{#2  }\ar @/^1pc/ @{{}{ }{}} [rrrr]|{\boxed{#5}}  \ar@/_2.5pc/ @{{}{ }{}} [rruuu]^{\boxed{1}}  & && & #1 \, \d(#2) \ar[llu]^{\overline{s}(g)} \ar[lluuu]_{1} \ar@/^1pc/ @{{}{ }{}} [lluuu]|{ \boxed{#9}}   }

\newcommand{\ten}[9] { & & & #1\, \d(#6)  \\ \\ & && #1 \ar[uu]|{#6 } \\ & #1\ar[uuurr]^{#2 } \ar[urr]_{1 } \ar[rrrr]_{#2  }\ar @/^1pc/ @{{}{ }{}} [rrrr]|{\boxed{#5}}  \ar@/_2.5pc/ @{{}{ }{}} [rruuu]^{\boxed{1}}  & && & #1\, \d(#2) \ar[llu]^{#4 } \ar[lluuu]_{1 } \ar@/^1pc/ @{{}{ }{}} [lluuu]|{ \boxed{#9}}   }

\newcommand{\teu}[9] { & & & #1 \d(#2) \d(#4) \d(#6)  \\ \\ & && #1 \d(#2) \d(#4)\ar[uu]|{#6 } \\ & #1\ar[uuurr]^{#2 #4 #6} \ar[urr]_{#2 #4 } \ar[rrrr]_{#2  }\ar @/^1pc/ @{{}{ }{}} [rrrr]^{\boxed{#5}}  \ar@/_2pc/ @{{}{ }{}} [rruuu]^{\boxed{#7} }  & && & #1 \d(#2) \ar[llu]^{#4 } \ar[lluuu]_{#4 #6 } \ar@/^1pc/ @{{}{ }{}} [lluuu]|{ \boxed{#9}}   }

\newcommand{\ted}[9] { && & & \phi(g) \,\d(s(g)) \,\d(s'(g)) \,\d(s''(g)) )  \\ \\ && && \phi(g) \,\d(s(g)) \,\d(s'(g)) \ar[uu]|{s''(g) } \\ & #1\ar[uuurrr]^{(s\tn s' \tn s'')(g)} \ar[urrr]_{(s \tn s')(g) } \ar[rrrrrr]_{#2  }\ar @/^1pc/ @{{}{ }{}} [rrrrrr]|{\boxed{#5}}  \ar@/_2pc/ @{{}{ }{}} [rrruuu]|<<<<<<<<<<<<<<<<<<<<<<{\boxed{\w^{(s \tn s',s'')}(g)}}  & &&&& & #1 \d(#2) \ar[lllu]^{s'(g)} \ar[llluuu]_{(s'\tn s'')(g)} \ar@/^2pc/ @{{}{ }{}} [llluuu]|<<<<<<<<<<<<<<<<<<<<<<<<<<{ \boxed{#9}}   }
\makeatletter

\begin{document}

\author[A]{Bj\"{o}rn Gohla}
\ead{b.gohla@gmx.de}

\author[B]{Jo\~{a}o Faria Martins}
\ead{jn.martins@fct.unl.pt}\

\address[A]{   Centro de Matem\'{a}tica da Universidade do Porto,
Departamento de Matem\'{a}tica da FCUP;\\ Rua do Campo Alegre, 687 ; 4169-007 Porto, Portugal  }

\address[B]{Departamento de Matem\'{a}tica and Centro de Matem\'{a}tica e Aplica\c{c}\~{o}es\\ Faculdade de Ci\^{e}ncias e Tecnologia (Universidade Nova de Lisboa),
Quinta da Torre, 2829-516 Caparica, Portugal}


 \title{Pointed homotopy and pointed lax  homotopy of 2-crossed module maps}

\begin{abstract}
We address the (pointed) homotopy theory of 2-crossed modules (of groups), which are known to  faithfully represent Gray 3-groupoids, with a single object, and also connected homotopy 3-types. The homotopy relation between 2-crossed module maps will be defined in a similar way to Crans' 1-{transfor}s between strict Gray functors, however being pointed, thus  this corresponds to  Baues' homotopy relation between quadratic module maps. Despite the fact that this homotopy relation between 2-crossed module morphisms is not, in general, an equivalence relation, we prove that if  $\A$ and $\A'$ are 2-crossed modules, with the underlying group $F$ of $\A$ being free (in short $\A$ is free up to order one), then  homotopy  between 2-crossed module maps $\A \to \A'$ yields, in  this case, an equivalence relation. Furthermore, if a chosen basis $B$ is specified for $F$, then  we can define a 2-groupoid  $\HOM_B(\A,\A')$ of 2-crossed module maps $\A \to \A'$, homotopies connecting them, and 2-fold homotopies between homotopies, where the latter  correspond to (pointed) Crans' 2-{transfor}s between 1-transfors.

We define a partial resolution $Q^1(\A)$, for a 2-crossed module $\A$, whose underlying group is free, with a {canonical} chosen basis, together with a projection map ${\rm proj}\colon Q^1(\A) \to \A$, defining isomorphisms at the level of 2-crossed module homotopy groups. This resolution (which is part of a comonad) leads to a {weaker notion of homotopy (lax homotopy)} between 2-crossed module maps, which we fully develop and describe. In particular, given 2-crossed modules $\A$ and $\A'$, there exists a 2-groupoid $\EuScript{HOM}_{\rm LAX}(\A,\A')$ of (strict) 2-crossed module maps $\A \to \A'$, and their lax homotopies and lax 2-fold homotopies, leading  to the question of whether the category of 2-crossed modules and strict maps can be enriched over the monoidal category {\emph Gray}.

The associated notion of a (strict) 2-crossed module map $f\colon \A \to \A'$ to be a lax homotopy equivalence has the two-of-three property, and it is closed under retracts. This discussion leads to the issue of whether there exists a  model category structure in the category of 2-crossed modules (and strict maps) where weak equivalences correspond to lax homotopy equivalences, and any free up to order one 2-crossed module is cofibrant. 
 \end{abstract}
\maketitle

\section{Introduction and simplicial group background / context}

Let $\G=\big(G_n,d_i^n,s_i^n; i\in\{0,1,\dots,n\}, n=0,1,2,\dots\big) $ be a simplicial group; \cite{Q,GJ,M,C}. As usual, see for example \cite{MuPo,GM}, we say that $\G$ is free if each group $G_n$  of $n$-simplices is a free group, with a chosen basis, and these basis are stable under the degeneracy maps $s_i^n\colon G_n \to G_{n+1}$. Recall that the Moore complex \cite{Po,Po2,Co} $N(\G)$ of a simplicial group $\G$ is given by the (normal) complex of groups $(\dots \to  A_n\ra{\d_n} A_{n-1} \to \dots \to A_0=G_0)$, 
where: $$A_n=\bigcap_{i=0}^{n-1} \ker(d_i^n),$$
and $\d_n\colon A_n \to A_{n-1}$ is the restriction of the boundary map $d^n_n \colon G_n \to G_{n-1}$. We say that the Moore complex of $\G$ has length $n$ if the unique (possibly) non trivial components of $N(\G)$ are $A_{n-1} \to A_{n-2}\to \dots \to A_0$. (Here ``length'' correspond to the number of groups, rather than the number of arrows, which is the usual convention).  Not surprisingly, this Moore complex has a lot of extra structure, defining what a hyper crossed complex is \cite{CC}, which retains enough information to recover the original simplicial group, up to isomorphism. This contains two well known results, stating that the categories of simplicial groups with Moore complexes of length  two and three (respectively) are equivalent to the categories of crossed modules and of 2-crossed modules  of groups (respectively), see \cite{Po,Po2,Co}, the latter being exactly hyper crossed complexes of length two and three (respectively); we will go back to this issue below. Hyper crossed complexes therefore generalise both crossed modules and 2-crossed modules.

 Looking at the last two stages of the Moore complex $N(\G)$  of a simplicial group $\G$, namely $\d=\d_1\colon N_1(\G) \to N_0(\G)$, one has an induced action of $N_0(\G)$ on $N_1(\G)$ by automorphisms, and  also the action of $N_0(\G)$ on itself by conjugation, and the boundary map $\d\colon N_1(\G)  \to N_0(\G)$ preserves these actions; in other words one has a pre-crossed module (\cite{BHS,B1,B2,MuPo}), called the pre-crossed module associated to the simplicial group $\G$.

The homotopy groups of a simplicial group $\G$ are, by definition, given by the homology groups of its Moore complex $N(\G)$ (which is a normal complex of, not necessarily abelian, groups). These correspond to the homotopy groups of the simplicial group $\G$ seen  as a simplicial set (despite the fact that $\pi_0(S)$, for $S$ a simplicial set, is not in general a group but a set). Note that simplicial groups are Kan complexes, and therefore {their} homotopy groups are well defined \cite{M,C}.

It is a fundamental result of Quillen \cite{Q,GJ} that the category of simplicial groups is a model category, where weak equivalences are the simplicial group maps $f\colon \G \to \G'$, inducing isomorphisms at the level of homotopy groups, and fibrations are the simplicial groups maps $f\colon \G \to \G'$ whose induced map on Moore complexes $(f_i,i=0,1,2\dots)\colon N(\G) \to N(\G')$ is surjective for all $i>0$, and in particular any object is fibrant. Notice that, for $f=(f_i)$ to be a fibration,  we do not impose that the induced map $f_0\colon A_0=G_0 \to A_0'=G_0'$ {be} surjective; however if $f$ is a {weak} equivalence and a fibration then certainly $f_0$ is surjective. Cofibrations are defined as being the maps that have the left lifting property with respect to all acyclic fibrations. In particular any  free simplicial group is cofibrant \cite{Q,GM}. The pre-crossed module associated to a free simplicial group is of the form $F_1 \to F_0$ where $F_0$ is a free group and $F_1 \to F_0$ is a free pre-crossed module, \cite{MuPo}. Such a pre-crossed module is what is called in \cite{MuPo} a totally free pre-crossed module.

Let $\SG$ denote the category of simplicial groups and $\SG_n$ denote the full subcategory of simplicial groups with Moore complex of length   $n$. The former is a reflexive subcategory of $\SG$ and we denote the reflexion functor (the $n$-type or $n^{\rm th}$-Postnikov section) by  $P_{n}\colon \SG \to \SG_n$, a left adjoint to the inclusion functor $\SG_n\to \SG$. At the level of Moore complexes $(A_m,\d_m)$ is sent, via $P_n$, to (see \cite{GM}): 
$$A_{n-1}/\d(A_{n}) \to A_{n-2} \to A_{n-3} \to \dots \to A_0 .$$  
This adjunction induces a closed model category structure on $\SG_n$ where fibrations (weak equivalences) are the maps whose underlying simplicial group map is a fibration (weak equivalence), and cofibrations are the maps which have the left lifting property with respect to all acyclic fibrations; all of this is explained for example in \cite{CG}. Therefore the $n$-type functor $P_n$ preserves cofibrations, \cite{CG1}. This will give model category structures in the categories of crossed modules and of 2-crossed modules (of groups); \cite{CG,CG1}. The case of crossed modules of groupoids in treated in \cite{No,MS,BGo}, not appealing directly to simplicial group(oid) techniques.

A crossed module $(\d\colon E \to G, \t)$ is given by a group morphism $\d\colon E \to G$, together with a left action $\t$ of $G$ on $E$ by automorphisms, so that $\d$ preserves the actions, where $G$ is given the adjoint action (thus $(\d\colon E \to G,\t)$ is a pre-crossed module), such that, furthermore,   for each $e,f \in E$ the Peiffer pairing $\langle e,f\rangle\doteq\big(efe^{-1                                                                                         }\big)\,\, \big(\d(e) \t f^{-1}\big) \in E$ vanishes.  The category of crossed modules is equivalent to the category of simplicial groups with Moore complex of length two, where morphisms of crossed modules $(\d\colon E \to G, \t) \to (\d\colon E' \to G', \t)$   are given by chain maps $(\psi\colon E \to E', \phi\colon G \to G')$, preserving the actions. Crossed modules form  a model category \cite{CG,MS}, where weak equivalences are the maps inducing isomorphisms on homotopy groups (the homology groups of the underlying complexes) and fibrations $(E \to G) \to (E' \to G')$ are the  crossed module maps $(\psi\colon E \to E', \phi\colon G \to G')$, such that  $\psi \colon E \to E'$ is surjective {(note that we do not require $\phi\colon G \to G'$ to be surjective)}. The homotopy category of {crossed modules} is equivalent to the homotopy category of pointed connected 2-types \cite{L,Po2,B1,B2} (where an $n$-type is a space $X$, {homotopic to a CW-complex,} such that $\pi_i(X)=\{0\}$, if $i >n$).
Since the reflexion functor $P_2\colon \SG \to \SG_2$ preserves cofibrations, a crossed module {$(\partial\colon E \to G)$} is cofibrant when $G$ is a free group, since in this case one can prove that there exists a free simplicial group (with Moore complex of length three) whose second Postnikov section is {$(\partial\colon  E \to G)$.} 

{The fact that, if $(\d\colon E\to G,\t)$ is a crossed module, with $G$ a free group, then it is cofibrant in the model category of crossed modules of groupoids is directly proven for example in \cite{No}. We note however that, considering the inclusion of the category of crossed modules of groups into the model category \cite{No,BG,MS} of crossed modules of groupoids, it does not hold that a fibration of crossed modules of groups is necessarily  a fibration of crossed modules of groupoids, even though the same is true for weak equivalences. In fact a map $(\psi,\phi)\colon (E \to G) \to (E' \to G')$, of crossed modules of groups,  is a fibration of crossed modules of groupoids if, and only if, both $\psi\colon E \to E'$ and $\phi\colon G \to G'$ are surjective.} {This apparent contradiction between the two  definitions of fibrations of crossed modules of groups should not come as a surprise, since crossed modules of groups (and their maps) model pointed homotopy classes of maps connecting pointed 2-types, whereas crossed modules of groupoids model non-pointed 2-types and  homotopy classes of maps between them.} {We note that a pointed fibration of pointed  simplicial sets is also not necessarily a fibration of simplicial sets; \cite[V-Lemma 6.6]{GJ}. Considering the nerve functor \cite{BHS}, from the category of crossed modules to the category of pointed simplicial sets, the nerve of a fibration $(\psi,\phi)\colon (E \to G) \to (E' \to G')$, of crossed modules of groups, is a fibration of simplicial sets if, and only if, $\phi(G/\d(E))=G'/\d(E')$, by \cite[V - Corollary 6.9]{GJ}, which is the same as saying that $\phi(G)=G'$.}

The 2-crossed modules were defined by Conduch\'{e} \cite{Co}, who proved that the category of 2-crossed modules is equivalent to the category of simplicial groups whose Moore complex has length three. A 2-crossed module is given by a complex of groups $(L \ra{\de} E \ra{\d} G)$, with a given action of $G$ on $E$ and $L$ by automorphisms, such that $(\d\colon E \to G)$ is a pre-crossed module, and also we have a map $\{,\}\colon E \times E \to L$ (the Peiffer lifting), lifting the Peiffer commutator map $<,>\colon E \times E \to E$, where, as before, $\langle e,f \rangle =(efe^{-1})\,(\d(e) \t f^{-1})$. This  lifting has to satisfy very natural properties, satisfied by the Peiffer pairing itself.

The category of 2-crossed modules is a Quillen model category, \cite{CG,CG1}, where a map $({\mu}, \psi, \phi) \colon (L \to E \to G) \to (L' \to E' \to G')$ is a fibration if, and only if, ${\mu}\colon L \to L'$ and $\psi\colon E \to E'$ each are surjective, and a weak equivalence if it induces isomorphisms of homotopy groups. Since the reflection functor $P_3\colon \SG \to \SG_3$ preserves cofibrations, one can see that any 2-crossed module {$\A=( L \to E \to G)$, with $(E\to G)$ being a totally free pre-crossed module (in short $\A$ is free up to order two),} is cofibrant in this model category. Similar results are proved in \cite{La}, where an analogous model category structure for Gray categories is constructed, and it is proven that a Gray category is cofibrant if, and only if, its underlying sesquicategory \cite{St} is free on a computad. We will go back to this issue below. {The homotopy category of 2-crossed modules  is equivalent to the homotopy category of 3-types: pointed CW-complexes $X$ such that $\pi_i(X)=\{0\}$ if $i>3$, where the Whitehead products $\pi_2(X) \times \pi_2(X) \to \pi_3(X)$ {are encoded} in the Peiffer lifting.}

{Crossed squares, defined in  \cite{GL}, were proven in \cite{L} to be models for homotopy 3-types, furthermore being equivalent to $cat^2$-groups. Ellis constructed in \cite{E} the fundamental crossed square of a pointed CW-complex, containing all relevant 3-type information.  By using Brown-Loday's theorem \cite{BL1,BL2}, this fundamental crossed square can be proven to be totally free, thus giving a combinatorial description of the first three homotopy groups of a CW-complex.}

{The quadratic modules, defined by Baues in \cite{B2} {(being  models for homotopy 3-types)} are a special case of 2-crossed modules, satisfying however additional nilpotency conditions, ensuring for example that the Peiffer liftings  $\{a,\langle b,c\rangle\}$ and $\{\langle a,b\rangle,c\}$ are trivial. This does simplify the calculations  in some cases, since the main difficulty in working with 2-crossed modules has to do with performing complicated computations with Peiffer liftings. It is proven in \cite{AU} that the category of quadratic modules is a reflexive subcategory of the category of 2-crossed modules. {Even though quadratic modules  retain all of the necessary flexibility from 2-crossed modules, in order to model homotopy 3-types, they have the drawback that the functor constructed in \cite[page 216]{B2} from the category of CW-complexes (and cellular maps) to the category of quadratic modules is somehow un-intrinsic, being inductively defined while attaching cells.} {On the other hand, the functor constructed in \cite{FM3}, which associates to a pointed CW-complex $X$  its fundamental 2-crossed module (constructed via Ellis' fundamental crossed square \cite{E} of a pointed CW-complex), which also models the homotopy 3-type of $X$, was intrinsically defined by using usual relative and triadic homotopy groups. This fundamental 2-crossed module of a reduced CW-complex is totally free. Therefore it has a combinatorial description, permitting (as done in \cite{B2,E}) a combinatorial description of the first three homotopy groups of a reduced CW-complex, and the remaining 3-type information, such as the Whitehead products $\pi_2 \times \pi_2 \to \pi_3$. In addition, maps between 3-types can be completely described, up to homotopy, by maps between fundamental 2-crossed modules (up to homotopy); see \cite{FM3}. The fundamental quadratic module of a CW-complex is a quotient of this fundamental 2-crossed module.}

 {An advantage that 2-crossed modules have, in comparison with the (very similar) crossed squares, is that there is only one group appearing in dimension two, which facilitates the definition of homotopy between maps. On the model category theoretical side, since the category of  2-crossed modules is equivalent to the category of simplicial groups with Moore complex of length three, this means that it comes with a wealth of homotopical algebra structures, transported  from the category of simplicial groups \cite{Q}.}

{In addition to simplicial groups with Moore complex of length three, there are several categories equivalent to the category of 2-crossed modules, for example the category of braided regular crossed modules \cite{BG} and the category of Gray 3-groupoids \cite{KP}, mentioned below. We emphasise that the category of 2-crossed modules is not equivalent to the category of crossed squares, even though there is a mapping cone functor from the latter to the former. The comparison between some of these algebraic models for 3-types appears in \cite{AU}. }

Recall, see \cite{Cr,Gu,GPS}, that a Gray 3-category $\mathcal{C}$ (or Gray enriched category) is a category enriched over the monoidal category of 2-categories, with the Gray tensor product. These can be given a more explicit definition, see \cite{Cr,KP}. In particular, one has sets $C_0$, $C_1$, $C_2$ and $C_3$ of objects, morphisms, 2-morphisms and 3-morphisms, and several operations between then. In particular, objects, 1-morphisms and 2-morphisms of $\mathcal{C}$ form a sesquicategory \cite{St}, a structure similar to a 2-category, but where the interchange law \cite{KP,HKK} does not hold. Nevertheless, the interchange law holds up to a chosen tri-morphism: the ``interchanger''.   Gray 3-categories correspond to the strictest version of tri-categories, in the sense that any tricategory can be strictified to a tri-equivalent Gray-category, \cite{Gu,GPS}. For this reason Gray 3-categories are also called semistrict tri-categories.

 Gray 3-groupoids can be defined analogously: there exists an inclusion of the category of 2-groupoids into the category of $\w$-groupoids \cite{BHS} (equivalent to crossed complexes of groupoids), which has a left adjoint $T$ (the cotruncation functor) similar to the $2$-type functor. The tensor product of $\w$-groupoids is treated in \cite{BHS,BH}, and from now on called {the Brown-Higgins tensor product. Composing} with the cotruncation functor yields a tensor product in the category of 2-groupoids (which is simply the restriction of the Gray tensor product of 2-categories to 2-groupoids), part of a monoidal closed structure. A Gray 3-groupoid is a groupoid enriched over this monoidal category of 2-groupoids. These are therefore  Gray 3-categories where any $i$-cell ($i \ge 1$) is strictly invertible.
It is folklore, and explicitly proven for example in \cite{KP,BG}, that any 2-crossed module defines a Gray 3-groupoid with a single object (a ``Gray 3-group''), and conversely. In particular the interchanger is derived from the Peiffer lifting in the given 2-crossed module.

The notion of a homotopy between $\w$-groupoid maps is treated in \cite{BH,BHS,Br4}. Considering $n$-fold homotopies between $\w$-groupoids maps defines an internal-hom  ``$\HOM$'' in the category of $\w$-groupoids, which, together with the Brown-Higgins tensor product of $\omega$-groupoids, induces a monoidal closed structure. By applying the cotruncation functor this yields a monoidal closed structure in the category of 2-groupoids.  If $A$ and $B$ are 2-groupoids then $\HOM(A,B)$ is simply the 2-groupoid of strict  functors $A \to B$, pseudo-natural transformations between them and their modifications. Considering 2-groups (2-groupoids with a single object) $A$ and $B$, and pointed {pseudo} natural transformations and modifications, yields simply a groupoid with objects the maps $A \to B$, and morphisms  their pointed pseudo-natural transformations. The latter groupoid is what we want to generalise for Gray 3-groups.

In \cite{Cr}, the notion of a pseudo-natural transformation (a 1-{transfor}) between  strict functors of Gray 3-categories was addressed. Unlike the notion of pseudo-natural transformations between strict functors of 2-categories, this does not define (in general) an equivalence relation between Gray functors, the problem being that we cannot \cite{Cr}, in general, concatenate, {or invert,}  1-transfors, due to the lack of the interchange condition in Gray 3-categories. To have an equivalence relation between strict Gray functors one needs the full force of tricategories, leading to a less restrictive notion of 1-transfors, see \cite{Gu,GPS,Go}. In particular given Gray 3-categories $A$ and $B$ one has a Gray 3-category $\mathcal{HOM}(A,B)$ of strict Gray functors and their lax 1-, 2- and 3-{transfor}s.

Given that strict 1-transfors between Gray 3-category maps $f,g\colon A \to B$ can be modeled by maps from the tensor product (of Gray categories \cite{Cr}) $A \tn I$ into $B$, where $A \tn I$ is a cylinder object for $A$, the fact that two maps being related by a 1-transfor does not yield, in general,  an equivalence relation is not at all surprising, given that in a model category such a construction would only be an equivalence relation, in general, if $A$ were cofibrant and $B$ {were} fibrant. As we have mentioned before $A$ is cofibrant if, and only if,  the underlying sesquicategory of $A$ is free on a computad, \cite{La}.

In this article, we will analyse the notion of 1-transfors connecting Gray 3-groupoid maps, as well as 2-transfors connecting these,  in  context of 2-crossed modules.  Similarly to \cite{B2,FM3}, we will restrict our discussion to the pointed case, and, since we are working with complexes of groups, we will name (strict) 1-{transfor}s and 2-transfors as being homotopies and 2-fold homotopies. These are very  similar to the usual notions of homotopies between chain complex maps, and 2-fold homotopies connecting them,  however adapted to the non-abelian case.

 To this end, given a 2-crossed module $\A'$, we will define a {(functorial)} path-space $\P(\A')$ for it, together with two fibrations ${\rm Pr}_0^{\A'},{\rm Pr}_1^{\A'}\colon \P(\A') \to \A'.$ This will be a good path space \cite{DS}, in the model category of 2-crossed modules. We use $\P(\A')$ to define the homotopy relation between 2-crossed module maps $\A \to \A'$. This coincides with the notion defined ad-hoc in \cite{FM3}.  {Dually we note that cylinder objects were constructed in the category of quadratic modules in \cite{B2} and in the category of crossed chain algebras in \cite{T}. The latter generalise 2-crossed modules.}

Let $\A=(L \to E \to G)$ and $\A'=(L'\to E' \to G')$ be 2-crossed modules (or more precisely the underlying complexes of them, since we also have actions $\t$ of $G$ on $E$ and  $L$ and a lifting $\{,\}\colon E \times E \to L$ of the Peiffer pairing, and the same for $\A'$.) A homotopy between the 2-crossed module maps $f,f'\colon \A \to \A'$, explicitly $f=({\mu}\colon L \to L', \psi \colon E \to E',\phi\colon G \to G')$ and $f'=({\mu}'\colon L \to L',\psi'\colon E \to E',\phi'\colon G \to G')$, is given by two set maps $s\colon G \to E'$ and $t\colon E \to L'$, satisfying appropriate properties (defining what we called a  quadratic $f$-derivation), resembling the notion of homotopy between quadratic module maps, treated in  \cite{B2}. For example, for all $g,h \in G$ we must have $s(gh)=\phi(h)^{-1}\t s(g) \,\,s(h)$, thus $s\colon G \to E'$ is to be a derivation.

{As expected by the discussion above, the notion of homotopy between 2-crossed module maps $\A \to \A'$ is not an equivalence relation (we give in \ref{fail} an explicit example to show that this is the case).  This issue can be fixed, for example, if  we consider the case when the underlying pre-crossed module {$(\partial \colon E \to G)$} of $\A$ is totally free, thus $\A$ is cofibrant in the model category of 2-crossed modules, which renders all objects fibrant. This point of view was the one considered in \cite[Lemma 4.7]{B2}, while proving that the homotopy relation between quadratic chain complex maps is an equivalence relation, in the totally free  case.}

What is surprising (and will be the main result of this paper) is that if $\A=(L \to E \to F)$ is a 2-crossed module, {such that $F$ is a  free group (in short $\A$ is free up to order one)}, then homotopy between 2-crossed module maps $\A \to \A'$ defines an equivalence relation. Moreover, if a free basis $B$ of $F$ is specified, then we can define a groupoid, whose objects are the 2-crossed module maps $\A \to \A'$ and the morphisms are their homotopies, represented as (for $(s,t)$ a quadratic $f$-derivation):
$$f \ra{(f,s,t)} f' .$$
This process can be continued, to define a 2-groupoid  $\HOM_B(\A,\A')$, with objects being the 2-crossed module maps $\A \to \A'$, the morphisms being  the homotopies between then, and the 2-morphisms being 2-fold homotopies between homotopies. We will present very detailed calculations.

Consider two  2-crossed modules: $$  \A=\left (L \ra{\de} E \ra{\d} F,  \t,\{,\}\right)  \an \A'=\left (L' \ra{\de} E'\ra{\d} G',  \t,\{,\}\right),$$ such that $F$ is a free group, with a chosen basis $B$. If the quadratic $f$-derivation  $(s\colon F \to E',t\colon E \to L')$ connects $f\colon \A \to \A'$ and $f'\colon \A \to \A'$ and $(s'\colon F \to E' ,t'\colon E \to L')$ connects $f'\colon \A \to \A'$ and  $f''\colon \A \to \A'$, diagrammatically:
$$f \ra{(f,s,t)} f' \textrm{ and } f' \ra{(f',s',t')} f'',$$
then we explicitly construct  a quadratic $f$-derivation $(s \tn s',f \tn f')$, such that: 
$$f \ra{(f,s \tn s',t \tn s'')} f''.$$
This {concatenation} ``$\tn$'' of homotopies is associative and it has inverses. 
 The derivation $s \tn s' \colon F \to E'$   is the unique derivation $F \to E'$ which, on the chosen basis $B$ of $F$, has the form $b \longmapsto s(b)s'(b)$. In the case when {$(\partial\colon  E' \to G')$} is a crossed module, and the Peiffer lifting of $\A'$ is trivial, we have that $(s\tn s')(g)=s(g)s'(g)$ for each $g \in F$. Otherwise one has a map $\w^{(s,s')}\colon  F\to L'$, measuring the difference between $s\tn s'$ and the pointwise product of $s$ and $s'$; namely we have:
$$(s \tn s')(g)= s(g) \,\, s'(g)\,\, \de\big(\w^{(s,s')}(g)\big)^{-1}, \textrm{ for each } g \in F.$$
This map $\w^{(s,s')}\colon F  \to L' $ has a prime importance in the construction of the {concatenation} of homotopies. 

Let $G$ be a group. The free group on the underlying set of $G$ is denoted by $\FG(G)$. The set inclusion $G \to \FG(G)$ is written as $g \in G \longmapsto [g] \in \FG(G)$. We have the obvious basis $[G]=\{[g], g \in G\}$ of $\FG(G)$. Consider the obvious  projection group map $p\colon \FG(G) \to G$, thus $p([g])=g$.
There exists a very natural partial resolution functor $Q^1$, from the category of 2-crossed modules to itself,  which, to a 2-crossed module 
$\A=\left (L \ra{\de} E \ra{\d} G ,  \t,\{,\}\right),$
associates the  2-crossed module:
$$Q^1(\A)=\left(L \ra{\de'} E {{}_\d \times p} \FG(G) \ra{\d'} \FG(G),\t,\{,\}\right),$$
where:
$$ E {{}_\d \times p} \FG(G)=\{(e,u)\in E \times \FG(G): \d(e)=p(u)\}.$$
Therefore $Q^1(\A)$ is 
free up  to order one. Moreover there is a projection ${\rm proj}= (r,q,p) \colon Q^1(\A) \to \A$ which   yields isomorphisms at the level of homotopy groups. It has the form:
\begin{equation*}
{\rm proj} \quad  =  \xymatrix{ &L \ar[r]^-{\de'}\ar[d]_r& E {{}_\d \times p} \FG(G) \ar[r]^-{\d'}\ar[d]_q& \FG(G)\ar[d]^p\\
&L \ar[r]_\de & E \ar[r]_\d  & G}
\end{equation*}
where $r=\id$ and $q(e,u)=e$.
Therefore  ${\rm proj}\colon Q^1(\A) \to \A$ is surjective and a weak equivalence (thus an acyclic fibration in the model category of 2-crossed modules).

The map  ${\rm proj}\colon Q^1(\A) \to \A$, resembling a cofibrant replacement, is proven to be part of  a comonad in \cite{Go} {(we note that the results in this paper are independent from the ones in \cite{Go}).}
 Its co-Kleisli category \cite{ML} leads to  weaker notions of maps $\A \to \A'$ between 2-crossed modules \cite{Ga,Go}, and of homotopies between strict 2-crossed module maps $\A \to \A'$, as well as of their 2-fold homotopies, yielding a 2-groupoid $\EuScript{HOM}_{\rm LAX}(\A, \A')$, of strict maps $\A \to \A'$, lax homotopies and lax 2-fold homotopies, which we fully describe.
Let us be specific. Let $\A$ and $\A'$ be 2-crossed modules. Let $\hom(\A,\A')$ denote the set of 2-crossed module maps $\A \to \A'$. Given that ${\rm proj}\colon Q^1(\A) \to \A $ is surjective, the map $f \in \hom(\A,\A') \longmapsto f\circ {\rm proj} \in \hom(Q^1(\A),\A')$ is an injection. A 2-crossed module map $Q^1(\A) \to \A'$ is said to be strict if it factors (uniquely) through ${\rm proj} \colon Q^1(\A) \to \A$. Then we define $\EuScript{HOM}_{\rm LAX}(\A,\A')$ as being the full sub-2-groupoid of $\HOM_{[G]}(Q^1(\A),\A')$, with objects  being the strict maps $Q^1(\A) \to \A'$ {(and we call the 1- and 2-morphisms of $\EuScript{HOM}_{\rm LAX}(\A,\A')$ lax 1- and 2-fold homotopies)}.  After presenting  $Q^1(\A)$ combinatorially (by generators and relations),  we will give a fully combinatorial description of lax homotopies  between strict  2-crossed module maps, and their lax 2-fold homotopies, therefore explicitly constructing the 2-groupoid  $\EuScript{HOM}_{\rm LAX}(\A, \A')$.

 Lax homotopies between strict 2-crossed module maps behave well with respect to composition by strict 2-crossed module maps. Therefore it is natural to conjecture that the category of 2-crossed modules, strict  2-crossed module maps, lax homotopies and  their 2-fold homotopies is a Gray 3-category. 

We say that a 2-crossed module map $f\colon \A \to \A'$ is a lax homotopy equivalence if there exists a 2-crossed module map $g\colon \A' \to \A$ such that $f\circ g$ and $g\circ f$ each are lax homotopic to $\id_{\A'}$ and $\id_{\A}$, respectively. Since we can concatenate lax homotopies between 2-crossed module maps, the class of lax homotopy equivalences has the two-of-three property; \cite{DS}. Given that we can compose lax homotopies with strict 2-crossed module maps, any retract of a lax homotopy equivalence is a lax homotopy equivalence.  All of this discussion leads to the issue of whether there exists a model category structure in the category of 2-crossed modules  {(different from the one already referred to, which has  as cofibrant objects the retracts of free up to order two 2-crossed modules)} where weak equivalences correspond to  lax homotopy equivalences, and where free up to order one 2-crossed modules are cofibrant, and so that the path-space constructed in this paper is still a good path-space.

\section{Preliminaries on pre-crossed modules and 2-crossed modules}
\subsection{Pre-crossed modules and crossed modules }\label{pcm}
\begin{Definition}[Pre-crossed module]\label{cm}
A {\it pre-crossed module} $(\d\colon E \to G,\t)$ is given by a group morphism $\d\colon E \to G$, together with a left action $\t$ of $G$  on $E$ by automorphisms,  such that the following relation, called ``first Peiffer relation'', holds:
 $$\d(g \t e)=geg^{-1}, \textrm{ for each } g \in G \textrm{ and each } e \in E.$$
A {\it crossed module} $(\d\colon E \to G,\t)$ is a pre-crossed module satisfying, further,  the ``second Peiffer relation'': $$\d(x) \t y =xyx^{-1}, \textrm{ for each }  x,y \in E.$$
\end{Definition}
Note that  in a crossed module $(\d\colon E \to G, \t)$ the subgroup $\ker (\d)\subset E$ is central in $E$.

Let $(\d\colon E \to G, \t)$ be a pre-crossed module.
Given $x,y  \in E$,  their {\it Peiffer commutator} is given by
$$\langle  x,y\rangle   =\big(xyx^{-1}\big) \big(\d(x) \t y^{-1}\big).$$
Thus a pre-crossed module  is a crossed module if, and only if, all of its Peiffer commutators are the identity of $E$. In any pre-crossed module it holds that for each $x,y \in E$ (and where $1_G$ is the identity of $G$):
$$\d(\langle x,y\rangle)=1_G. $$

A morphism $f=(\p,\f)$ between the pre-crossed modules $(\d\colon E \to G, \t)$ and $(\d'\colon E' \to G', \t')$ is given by  a pair of group morphisms $\p\colon E \to E'$ and $\f\colon G \to G'$ making the diagram:
$$\begin{CD} E @>  \d>>   G\\
             @V\p VV  @VV\f V \\
             E' @>  \d'>>   G'
  \end{CD}
 $$
commutative, and such that  $$\p(g \t e)=\f(g) \t' \p(e), \fo  e \in E \an g \in G.$$
Morphisms of crossed modules are defined analogously, and therefore the category of crossed modules  is a full subcategory of the category of pre-crossed modules.

\begin{Example}[The underlying group  functor $\Gr_0$]\label{RGR1} There is an underlying group functor $\Gr_0$ sending a pre-crossed module $(E \to G,\t)$ to the group $G$. This has a right adjoint $R$ sending a group $G$ to the crossed module $(\id\colon G \to G,\ad)$, where we consider the adjoint action $\ad$ of $G$ on $G$. The unit of this adjunction, yields a pre-crossed module map $\eta_{\G}=(\d,\id)\colon (\d\colon E \to G,\t) \to (\id\colon G \to G,\ad)\doteq \eta(\partial\colon  E \to G,\t)$, for each pre-crossed module $\G=(\d\colon E \to G,\t)$. 
\end{Example}

\begin{Definition}[The principal group functor $\Gr_1$]\label{princ}
There is  a principal group functor $\Gr_1$ from the category of pre-crossed modules to the category of groups, sending  a pre-crossed module $(E \to G,\t)$  to  $E$.
\end{Definition}
The category of pre-crossed modules is fibered \cite{ML,BHS} over the category of groups, by considering the underlying group functor $\Gr_0$. Cartesian maps are easy enough to describe, \cite{AAO}:
\begin{Lemma}\label{cart}
 Let ${\cal G}=(\d\colon E \to G,\t)$ and ${\cal P}=(\d\colon M \to P,\t)$ be pre-crossed modules. Consider a pre-crossed module map $(\psi\colon E \to M,\phi\colon G \to P)\colon {\cal G} \to {\cal P}$. Suppose that  $E$ is isomorphic to the obvious  pullback, namely that there exists a group isomorphism: $$f\colon E\to G {{}_{\phi} \times_ \d} M\doteq \{(g,m)\in G \times M: \phi(g)=\d(m)\},$$ and that $\d\colon E \to G$ and $\psi\colon E \to M$ are the compositions of $f$ with the obvious projections $G {{}_{\phi} \times_ \d} M\to G$ and $ G {{}_{\phi} \times_ \d} M\to M$. Moreover, suppose that $f(g \t e)=g \t f(e)$, where $g \in G$ and $e \in E$, where we put $g \t (h,m)=(ghg^{-1},\phi(g) \t m)$, for each $g \in G$ and $(h,m) \in  G {{}_{\phi} \times_ \d} M$.
 Consider another pre-crossed module ${\cal H}=(\d\colon W \to H). $ Suppose we are given group maps  $\zeta\colon H \to G$ and $\phi'\colon H \to P $ with $\phi\circ \zeta=\phi'$ and a group map $\psi'\colon W \to E$, making $(\psi',\phi')\colon {\cal H} \to {\cal P}$ a pre-crossed-module morphism.
Then there exists a unique group morphism $\zeta^*\colon W \to E$, making the diagram:
\begin{equation}\label{gfg}
\xymatrix{ && E \ar[r]^{\psi} \ar[dd]|<<<<<\hole^{\d}& M\ar[dd]^{\d} \\&W \ar@{-->}[ru]^{\zeta^*} \ar[rru]_<<<<{\psi'}\ar[dd]_{\d} &&\\ && G\ar[r]^{\phi} & P \\ &H\ar[ru]^{\zeta} \ar[rru]_{\phi'}&& }
\end{equation}
commutative. Moreover $(\zeta^*,\zeta)\colon {\cal H} \to {\cal G}$ is a pre-crossed module map. 
\end{Lemma}
\begin{remark}\label{car2}
 Note that, given the form of pull-backs in the category of groups, there exists a unique set map $\zeta^*\colon W \to E$ that makes the diagram \eqref{gfg} commutative. Therefore a set map $\zeta^*\colon W \to E$ that makes  \eqref{gfg} commutative will immediately give a pre-crossed module map $(\zeta^*,\zeta)\colon {\cal H} \to {\cal G}$. 
\end{remark}
\begin{Proof}
 That the group map $\zeta^*\colon W \to E$ exists and is unique follows immediately since $E$ is a pull-back. In particular $\d\circ \zeta^*=\zeta\circ \d$.  Given $g \in H$ and $w \in W$ then
\begin{align*}
 \psi(\zeta^*(g \t w))&=\psi'(g \t w)=\phi'(g) \t \psi'(w)=\phi (\zeta(g)) \t \psi'(w).\\
 \d(\zeta^*(g \t w))&=\zeta(\d(g \t w))=\zeta(g \,\d(w)\, g^{-1})=\zeta(g)\,\, \zeta(\d(w))\,\, \zeta(g)^{-1}.
\end{align*}
Therefore if  $g \in H$ and $w \in W$ then $f\big(\zeta^*(g \t w)\big)=\zeta(g) \t f(\zeta^*(w)).  $
\end{Proof}
\subsection{Definition of 2-crossed modules and elementary properties. The secondary action  $\t'$ of a 2-crossed module}\label{d2cm}
We will follow the conventions of \cite{Co,FMPi} for the definition of a 2-crossed module. Important references on 2-crossed modules are  \cite{MuPo,KP,BG,Po}.
\begin{Definition}[2-crossed module]\label{2cmlg}
A 2-crossed module (of  groups) is given by {a chain complex} of  groups:
$$L \ra{\de} E\ra{\d} G$$
together with left actions $\t$, by automorphisms, of $G$ on $L$ and  $E$ (and on $G$ by conjugation), and a $G$-equivariant function {$\left\{,\right \} \colon E \times E \to L$} (called the Peiffer lifting). {Here $G$-equivariance means:}
\begin{equation*}
g \t \{e,f\}= \{g \t e, g \t f\}, \fo g \in G \an e,f \in E.
\end{equation*}
{These are to satisfy:}
\begin{enumerate}
 \item $L \ra{\de} E\ra{\d} G$ is {a chain complex} of $G$-modules (in other words $\d$ and $\de$ are $G$-equivariant {and $\d \circ \de =1$}.)
\item $\de(\left\{e,f\right \})=\left <e,f\right >,$ for each $e,f \in E$. {Recall that  $\left <e,f\right >=(efe^{-1})( \d(e) \t f^{-1}). $}
\item $[l,k]=\left\{\de(l),\de(k)\right \}, $ for each $l,k \in L$. Here $[l,k]=lkl^{-1}k^{-1}$.
\item $\left\{\de(l),e\right \}\left\{e, \de(l)\right \}=l(\d(e) \t l^{-1})$, for each $e\in E$ and $l \in L$.
\item $\left\{ef,g\right \}=\left\{e,fgf^{-1}\right \}\d(e)\t \left\{f,g\right \}$, for each $e,f,g \in E$.
\item $\left\{e,fg\right \}=\left \{e,f\right\}(\d(e) \t f) \t' \left \{e,g\right\}$, where {$e,f,g \in E$}.
\end{enumerate}
Here we have put:
\begin{equation}\label{tprime}
e \t'l=l\left\{\de(l)^{-1},e\right\}, \wh l\in L \an e \in E.
\end{equation}
\end{Definition}
The following is very well know; see \cite{Co,Po}. 
\begin{Lemma}[Secondary action $\t'$ of a 2-crossed module]\label{sa}
Let  $\A=(L \ra{\de} E \ra{\d} G, \t,\{,\}) $ be a 2-crossed module. The map $(e,l) \in G \times E \longmapsto e \t' l$ of equation \eqref{tprime} is a left action of $E$ on $L$ by automorphisms, called the ``secondary action of ${{\A}}$''. Together with the map $\de\colon L \to E$, the action $\t'$  defines a crossed module.
\end{Lemma}

 Note that in particular it follows that $\ker (\de) \subset L$ is central in $L$. We also have:
\begin{equation}
\left\{\de(l)^{-1},e\right\}^{-1}l^{-1}=(e \t' l)^{-1}=e \t'l^{-1}=l^{-1}\left\{\de(l),e\right\},
\end{equation}
\begin{align}\label{abc}
& e \t' l=\{\de(l),e\}^{-1}\,l,
 &\d(e) \t l=(e\t'{l})\,\{e,\de(l)^{-1}\}, &&\d(e) \t l =\{e, \de(l)\}^{-1} \,\, e \t' l.\end{align}
{Therefore $\d(a) \t l=a \t'l$ if $a \in E$ and $l \in \ker \de$; equation (\ref{v}). For each $a,b,c \in E$ we have:}
\begin{equation} { a \t' \{b,c\}=\d(a) \t \{b,c\}\big \{a,\left <b,c \right>^{-1}\big\}^{-1}= \{\d(a) \t b, \d(a) \t c\}\{a,(\d(b) \t c) bc^{-1}b^{-1}\}^{-1}.}
\end{equation}

A morphism $f=(\mu,\psi,\phi)$ between the 2-crossed modules $\A_1=(L_1 \to E_1 \to G_1, \t_1, \{,\}_1)$ and $\A_2=(L_2 \to E_2 \to G_2, \t_2, \{,\}_2)$ is given by group morphisms ${\mu}\colon L_1 \to L_2, \psi \colon E_1 \to E_2$ and $\phi\colon G_1 \to G_2$, defining a chain map between the underlying complexes, such that, for each $e,f \in E$, $g \in G$ and $k \in L$: 
$${\mu}(\{e,f\}_1)=\{\psi(e),\psi(f)\}_2, \quad {\mu}(g \t_1 k)=\phi(g) \t_2 {\mu}(k) \an  \psi(g \t_1 e)=\phi(g) \t_2 \psi(e).$$ 
The set of 2-crossed module morphisms $\A_1 \to \A_2$ is denoted by $\hom(\A_1,\A_2)$.

For a proof of the following  lemma we refer to  \cite{FM3,FMPi,Co}. 
\begin{Lemma}\label{rnn}In a 2-crossed module  $(L\ra{\de} E \ra{\d} G, \t, \{,\})$  we have, for each $e,f,g \in E$, $a\in G$ and $k \in L$:
\begin{align}
\{e,1_E\}&=\{1_E,e\}=1_L,  &    a \t (e \t' k)&=(a \t e) \t' (a \t k),   \label{v}                                                                                                                                      \\\label{imppp}
\{ef,g\}&= (e \t' \{f,g\} ) \{e,\d(f) \t g\},  &\{e,fg\}&=\left((efe^{-1}) \t' \{e,g\} \right)\{e,f\},\\
\{e,f\}^{-1}&=\d(e) \t \{e^{-1},efe^{-1}\}, 
&\{e,f\}^{-1}&=(efe^{-1}) \t'\{e,f^{-1}\}, \label{iNv}\\
\{e,f\}^{-1}&=(\d(e) \t f) \t' \{e,f^{-1}\}, \label{iNvv}
&\{e,f\}^{-1}&=e \t' \{e^{-1},\d(e) \t f\}.
\end{align}
\end{Lemma}

\begin{Example}\label{gh}
Given a pre-crossed module $E \to G$, consider the Peiffer subgroup $\langle  E,E\rangle \subset E$, generated by the Peiffer commutators $\langle  a,b \rangle$; see subsection \ref{pcm}. Then $\langle  E,E\rangle  \to E \to G$ is a 2-crossed module, where the Peiffer lifting is $\{a,b\}=\langle  a, b \rangle $. 
\end{Example}

\begin{Example}[The underlying pre-crossed module functor $T$]
  There is a truncation functor $T$, or underlying pre-crossed module functor, sending a 2-crossed module $\A=(L \ra{\de} E \ra{\d} G, \t,\{,\}) $ to its underlying pre-crossed module $(\d\colon E \to G,\t)$. This has a right adjoint sending a pre-crossed module $(\d\colon E \to G,\t)$ to $(\ker(\d) \to E \to G,\t)$, where we considered the inclusion map $\ker(\d) \to E$, and the Peiffer lifting is as in the previous example. A left adjoint to the truncation functor was constructed in \cite{FM3}.
\end{Example}

\begin{Example}[The underlying group $\Gr_0$  and principal group $\Gr_1$ functors]\label{RGR2}
 There is an underlying group functor $\Gr_0$ sending a 2-crossed module $(L \to E \to G)$ to the group $G$. This has a right adjoint  sending a group $G$ to the 2-crossed module $(\{1\} \to G \to G)$, considering the identity map $G\to G$, the adjoint action of $G$ on $G$; and the trivial Peiffer lifting. Compare with example \ref{RGR1}. On the other hand the principal group functor $\Gr_1$ sends a  a 2-crossed module $(L \to E \to G)$ to the group $E$.
\end{Example}

\begin{Definition}[Freeness up to order one] {We say that {a} 2-crossed module $\A=(L \ra{\de} E \ra{\d} F, \t,\{,\})$ is free up to order one if $F=\Gr_0(\A)$ is a free group. In this paper, free up to order one 2-crossed modules $\A$ will always come equipped with a specified chosen basis of the underlying group  $\Gr_0(\A)$.}
\end{Definition}

\begin{Definition}[Homotopy  groups of a 2-crossed module]
Given a 2-crossed module $\A=(L \ra{\de} E \ra{\d} G, \t,\{,\})$ then both $\im(\d) \subset G$ and $\im(\de)\subset E$ are normal subgroups. This permits us to define the homotopy groups $\pi_i(\A)$, where $i=1,2,3$ as the first three homology groups of the underlying complex of $\A$.
\end{Definition}
 \subsection{The path space  of a  2-crossed module}\label{tps}

Let ${\A}=(L \ra{\de} E \ra{\d} G, \t,\{,\})$ be a 2-crossed module. Let us define the path space $\P({\A})$ of it, together with two surjective  2-crossed module morphisms $\P({\A})\substack{\xrightarrow{{\rm Pr}_1^{{\A}}} \\ \ra[{\rm Pr}_0^{{\A}}]{} }  {{\A}} ,$ and an inclusion $i_{{\A}}\colon {{\A}} \to \P({{\A}})$, with ${\rm Pr}_1^{{\A}}\circ i_{{\A}}={\rm Pr}_0^{{\A}}\circ i_{{\A}}=\id_{{{\A}}}$.

\subsubsection{The derived action $*$ and the first and second lifted actions $\bullet$ of a 2-crossed module}
Most of this discussion appeared in \cite{FM3}.
 \begin{remark}[Convention on semidirect products]\label{sd} Given a left action $\t$ of the group $G$ on the group $E$ by automorphisms, the convention  for the semidirect product $G \ltimes_{\t} E$ is: $$(g,e)(g',e')=\big(gg',(g'^{-1} \t e) e'\big).$$ 
In particular given $g \in G$ and $e \in E$ we have:
$$(g,e)^{-1}=\big(g^{-1},g \t e^{-1}\big). $$
\end{remark}

 We resume the notation of \ref{d2cm}. 
{Consider the left action  of $E$ on $L$ (the secondary action $\t'$ of ${{\A}}$) given by {$e \t' k\doteq k\{\de(k)^{-1},e\}$,} where $e \in E$ and $k\in L$; lemma \ref{sa}.  For the following  see \cite{Co,FM3}}.

\begin{Lemma}[Derived action]\label{derived}
Let ${{\A}}=(L \ra{\de} E \ra{\d} G, \t,\{,\})$ be a 2-crossed module. There exists a left action $*$ of $E$ on  $ E \ltimes_{\t'} L$, by automorphisms (called the ``derived action of ${{\A}}$''), with the form:
\begin{align*}
b*(e,k)&= \left (\d(b) \t e, \left ( b \t'\{b^{-1}, \d(b) \t e^{-1}\}\right) \,\, b \t' k\right ) = \left ( \d(b) \t e, \{b,  e^{-1}\}^{-1}\,\, b \t' k \right ), \wh e,b \in E \an k \in L,
\end{align*}
{see equation (\ref{iNvv}).} Note that if $e \in E$ and $k \in L$:
\begin{equation}\label{refnow}
 b*(\de(k),k^{-1})=\big(\de( \d(b) \t k), \d(b) \t k^{-1}\big).
\end{equation}
\end{Lemma}
Consider the group $E \ltimes_*( E \ltimes_{\t'} L)$, whose group law is (remark \ref{sd}):
\begin{equation}\label{gl}
\begin{split}
(a,e,k)(a',e',k')&=\left(aa', ( \d(a'^{-1} )\t e) e', \Big((a'e')^{-1} \t'\big (\{a',\d(a')^{-1} \t e^{-1} \} k \big ) k' \right)\\
&=\left(aa', ( \d(a'^{-1} )\t e) e', \Big(e'^{-1} \t'\big (\{a'^{-1}, e^{-1} \}^{-1}\,\, a'^{-1} \t' k \big ) k' \right).
\end{split}
\end{equation} 

Particular cases of the multiplication are:
\begin{align}
(a,1,k)(a',1,k')&=\left(aa',1,(a'^{-1} \t k) k'\right), \wh a,a' \in E \an k,k'\in L,\\
(1,e,k)(1,e',k')&=\left(1,ee',(a'^{-1} \t k) k'\right), \wh e,e' \in E \an k,k'\in L. \label{a1}
\end{align}
Thus since $(\de\colon L \to E,\t')$ is a crossed module:
\begin{align}
\big(\de(l),1,k\big)\big(\de(l'),1,k'\big)&=\big(\de(l)\de(l'),1,l'^{-1}kl'k'\big) , \wh k,k',l,l' \in L;\label{k4}\\
\big(1,\de(l),k\big)\big(1,\de(l'),k'\big)&=\big(1,\de(l)\de(l'),l'^{-1}kl'k'\big), \wh k,k',l,l' \in L.\label{a2}
\end{align}
 {Put $a=(a,1,1)$, $e=(1,e,1)$ and $k=(1,1,k)$,  and the same for their images under $\t$ and $\t'$. We have:}
\begin{equation}\label{mnn}
 (a,e,k)=aek, \quad a ka^{-1}=a \t' k, \quad eke^{-1}=e \t' k .
\end{equation}
Putting $l=(1,1,l)$, for $l \in L$:
\begin{equation}\label{lkj}
 klk^{-1}=\de(k) \t' l= (\de(k),1,1)\,(1,1, l) \,(\de(k),1,1)^{-1}=(1,\de(k),1) \,(1,1, l) \, (1,\de(k),1)^{-1}.
\end{equation}

Consider the semidirect product $G \ltimes_\t E$. The following essential lemma appeared in \cite{FM3}.
\begin{Lemma}[First lifted action]\label{lifted}
Let ${{\A}}=(L \ra{\de} E \ra{\d} G, \t,\{,\})$ be a 2-crossed module of groups. 
There exists a left action by automorphisms $\bullet$ of $G \ltimes_\t E$ on $E \ltimes_*( E \ltimes_{\t'} L)$ (the first lifted action of ${{\A}}$), with:
\begin{align*}
(g,x)\bullet (a,e,k)&=\Big(g \t a, g \t \big((\d(a)^{-1} \t x)\, e\, x^{-1}  \big), 
g \t \Big ((x e^{-1}) \t' \Big \{a^{-1},  x^{-1}\Big\}^{-1} \Big)  \, g \t \Big \{ x, e^{-1}a^{-1} \Big\}   \big(g \d(x)\big) \t k \Big)\\
 &\doteq g \t \Big ( a,  (\d(a)^{-1} \t x)\, e\, x^{-1} , 
(x e^{-1}) \t' \Big \{a^{-1},  x^{-1}\Big\}^{-1}  \,\, \Big \{ x, e^{-1}a^{-1} \Big\}   \,\, \d(x)  \t k \Big).
\end{align*}
\end{Lemma}
Put $(1_G,x)=x$ and $(g,1_E)=g$. Particular cases of the first lifted action, which will be important later are:
\begin{equation}\label{gact}  g\bullet (a,e,k)=(g\t a, g \t e, g \t k). \end{equation}
\begin{equation}
x\bullet (a,e,k)=\Big( a,  (\d(a)^{-1} \t x)\, e\, x^{-1} ,
 \Big ((x e^{-1}) \t' \Big \{a^{-1},  x^{-1}\Big\}^{-1} \Big)    \Big \{ x, e^{-1}a^{-1} \Big\}    \d(x)  \t k \Big),
\end{equation}
\begin{equation}\label{m2}
x \bullet (\de(k),1,l)=(\de(k),1, k^{-1}\d(x) \t (kl)), \quad \wh  x \in E, \an k,l \in L,
\end{equation}
\begin{equation}\label{m3}
x\bullet (1,e,k)=\Big( a,  x  e x^{-1} ,
\{ x, e^{-1}\}  \,\,\d(x) \t k\Big),
\end{equation}
\begin{equation}\label{m5}
(g,x)\bullet (1,e,k)=\Big( 1, g \t (x e x^{-1} ), g \t \big (
  \{ x, e^{-1} \}    \d(x)  \t k \big) \Big),
\end{equation}
\begin{equation}\label{m4}
x\bullet (1,\de(k),k^{-1})=\Big( 1,  x  \de(k) x^{-1} ,
\{ x, \de(k)^{-1}\}  \,\,\d(x) \t k^{-1}\Big)=\Big( 1,  x  \de(k) x^{-1} ,x \t' k^{-1}\Big),
 \end{equation}
\begin{equation}\label{referright}
(g,x)\bullet (1,\de(k),k^{-1})=\Big(1,  g \t (x \de(k) x^{-1} ), g \t \big (
  \{ x, \de(k)^{-1} \}  \,\,  \d(x)  \t k^{-1} \big) \Big)=\Big(1,  g \t \de (x \t' k ), g \t \big (
 x \t' k^{-1} \big) \Big),
\end{equation}
\begin{equation}
x \bullet (a,1,1)=\Big(a,  \d(a)^{-1} \t x\,\,  x^{-1} ,  \big (x \t' \big \{a^{-1},  x^{-1}\big\}^{-1} \big)   \big \{ x, a^{-1} \big\}\Big),
\end{equation}
\begin{equation}
 \de(k) \bullet (a,1,1)=\Big(a,  \d(a)^{-1} \t  \de(k)  \,\,   \de(k)^{-1} , k \,\,\d(a) \t k^{-1}\Big),
\end{equation}
If  $x$ is such that $\{x,e\}=\{e,x\}=1$, for all $e \in E$, we have:
\begin{equation}\label{ker}
x\bullet (a,e,k)=\Big( a,  (\d(a)^{-1} \t x) e x^{-1} ,   \d(x)  \t k \Big).
\end{equation}

Let ${{\A}}=(L \ra{\de} E \ra{\d} G, \t,\{,\})$ be a 2-crossed module. If $\ad$ is the adjoint action, it follows easily that:
\begin{Lemma}[Second lifted action]
There is a left action $\bullet$ of $G \ltimes_\t E$ in $L \ltimes_\ad L$, by automorphisms (called the second lifted action of $\A$), which has the form ({where $x\in E$, $g \in G$ and $(k,l) \in L \ltimes_\ad L$}):
$$x \bullet (k,l)=(k, k^{-1}\d(x) \t (kl)) \an g \bullet (k,l)=(g \t k,g \t l).$$ 
\end{Lemma}

\subsubsection{Definition of the path space $\P({{\A}})$ of a 2-crossed module ${{\A}}$}\label{path}
By using  \eqref{m2}, \eqref{k4} and \eqref{gl}, we can see that  the maps 
\begin{equation}\label{defab}
\begin{split} 
(a,e,k)\in E \ltimes_* (E \ltimes_{\t'} L) &\stackrel{\beta}{\longmapsto} (\d(a),e) \in G \ltimes_\t E\\  
  (k,l) \in L \ltimes_{\ad} L & \stackrel{\a}{\longmapsto} (\de(k),1_E,l) \in E \ltimes_* (E \ltimes_{\t'} L)  
\end{split}
\end{equation}
are $(G \ltimes_{\t} E)$-equivariant group morphisms, with respect to the lifted actions $\bullet$ of ${{\A}}$ and the adjoint action of  $G \ltimes_\t E$ on itself.  This defines {a chain complex} of groups:
\begin{equation}\label{gcci}
L \ltimes_{\ad} L \ra{\alpha} E \ltimes_* (E \ltimes_{\t'} L) \ra{\beta} G \ltimes_\t E.
\end{equation}
{The Peiffer pairing in the pre-crossed module  $  E \ltimes_* (E \ltimes_{\t'} L) \ra{\beta} G \ltimes_\t E$ was calculated in \cite{FM3}:
\begin{equation}\label{ppfs}
 \langle ( a,e,k),(a',e',k')\rangle=\big (\langle a,a'\rangle,1,\{a,a'\}^{-1} \{ae\de(k),a'e'\de(k')\} \big).
\end{equation}
{By using  example \ref{gh}, and the structure of the group complex \eqref{gcci}, namely equations (\ref{k4}) and the form of the product in $L \ltimes_{\rm \ad} L$, there follows that there exists a}  2-crossed module structure in \eqref{gcci}, {whose Peiffer lifting $|,|\colon \big(E \ltimes_* (E \ltimes_{\t'} L)  \big) \times \big( E \ltimes_* (E \ltimes_{\t'} L) \big) \to L \ltimes_{\ad} L$ takes the following form:}
\begin{equation}\label{Pp}
 |(a,e,k),(a',e',k')|=\big (\{a,a'\},\{a,a'\}^{-1} \{ae\de(k),a'e'\de(k')\} \big).
\end{equation}
\begin{Definition} \label{ppsp} For a 2-crossed module $\A=(L \to E \to G,\t,\{\,\})$,
the 2-crossed module \begin{equation}\label{mnbv}
\P({{\A}})= \left ( L \ltimes_{\ad} L \ra{\alpha} E \ltimes_*( E \ltimes_{\t'} L)  \ra{\beta} G \ltimes_\t E,\bullet, |,| \right)
\end{equation} just defined will be called the (pointed) path-space of ${{\A}}$.  Clearly the path-space construction $\P$ is functorial with respect to 2-crossed module morphisms.
\end{Definition}

 By straightforward calculations we conclude that:
\begin{Theorem}\label{mor}Let ${{\A}}=(L \ra{\de} E \ra{\d} G, \t,\{,\})$ be a 2-crossed module.
 The maps $(k,l)\in L \ltimes_{\ad} L \stackrel{p'}{\longmapsto} kl \in   L$ and  $$(a,e,k) \in E \ltimes_*( E \ltimes_{\t'} L)\stackrel{q'}{\longmapsto} ae\de(k)\in E, \quad  (g,e) \in G \ltimes E \stackrel{r'}{\longmapsto} g \d(e) \in G;$$ and 
also:
$$(k,l)\in L \ltimes_{\ad} L \stackrel{p}{\longmapsto} k \in  L, \quad (a,e,k) \in E \ltimes_*( E \ltimes_{\t'} L)\stackrel{q}{\longmapsto} a \in E, \quad (g,e) \in G \ltimes E \stackrel{r}{\longmapsto} g  \in G$$ 
 are group morphisms. Moreover the triples ${\rm Pr}_0^{{\A}}\doteq (p,q,r)$ and ${\rm Pr}_1^{{\A}}\doteq (p',q',r')$ define surjective morphisms  $\P({{\A}}) \to {{\A}}$ of 2-crossed modules. We also have an inclusion map $i_{{\A}}\colon {{\A}} \to \P({{\A}})$ such that:
\begin{align*}
 &g \longmapsto (g,1), & e& \longmapsto (e,1,1), & k& \longmapsto (k,1).
\end{align*}
Therefore $ {\rm Pr}_0^{{\A}}\circ i_{{\A}}$ and  $ {\rm Pr}_1^{{\A}}\circ i_{{\A}}$ each are the identity of ${{\A}}$. Moreover the map $({\rm Pr}_0^{{\A}},{\rm Pr}_1^{{\A}})\colon \P({{\A}}) \to {{\A}} \times {{\A}}$ is a fibration of 2-crossed modules, considering the model category structure in the category of 2-crossed modules defined in \cite{CG}, see the introduction. This is because both maps $(p,p')\colon L \ltimes_\ad L \to L \times L$ and $(q,q')\colon   E \ltimes_*( E \ltimes_{\t'} L) \to E \times E$ are surjective. Therefore $\P({{\A}})$ is a good path-space for ${{\A}}$, \cite{DS}, {since clearly $i_{{\A}}\colon {{\A}} \to \P({{\A}})$ induces isomorphism at the level of homotopy groups.}
\end{Theorem}
It is convenient to ``visualise'' the path space of a 2-crossed module ${{\A}}$ in the following way, putting emphasis on the projection maps $\Pr_0^\A,\Pr_1^\A\colon \P({{\A}}) \to {{\A}}$:
\begin{equation}\label{pathv}
\xymatrix{\pathspace}.
\end{equation}
The elements $(g,e) \in \Gr_0(\P({{\A}}))$ and $(a,e,k)\in \Gr_1(\P({{\A}}))$ (example \ref{RGR2}) will several times be denoted as:
\begin{equation}\label{qas}
(g,e)=\big(\pu{g}{e}\big) \textrm{ and } (a,e,k)=\big(\pt{a}{e}{k} \big).
\end{equation}
The 2-crossed module boundary map $\beta \colon \Gr_1(\P({{\A}}))\to \Gr_0(\P({{\A}}))$ being (in this notation):
$$ \big(\pt{a}{e}{k} \big) \stackrel{\beta}{\longmapsto} \big(\pu{\d(a)}{e} \big) .$$
Note  that (in the notation of theorem \ref{mor}), if $(a,e,k)\in \Gr_1(\P(\A))$, $(g,e) \in \Gr_0(\P(\A))$,		 $g \in G$ and $a \in E$:
\begin{align}\label{qaa}
&\Pr_0^\A\big( \pt{a}{e}{k}   \big)=a&   &\Pr_0^\A\big(   \pt{a}{e}{k} \big)=a e \de(k)\\
&i_\A(g)=\big(g\ra{1_E} g\big)& & i_\A(a)=\big(a \ra{(1_E,1_L)} a\big)
\end{align}
\subsection{The double path-space $\P(\P({{\A}}))$ of a 2-crossed module ${{\A}}$}\label{dps}

Given a 2-crossed module  ${{\A}}=(L \ra{\de} E \ra{\d} G, \t,\{,\})$, we can iterate the path-space construction, obtaining a 2-crossed module $\P(\P({{\A}}))$, the double path-space of ${{\A}}$, with underlying  group complex:
\begin{multline}\label{aa}
(L \ltimes_{\ad} L) \ltimes_{\ad} (L \ltimes_{\ad} L) \ra{\alpha^\sharp} (E \ltimes_* (E \ltimes_{\t'} L))\ltimes_*((E \ltimes_* (E \ltimes_{\t'} L)\ltimes_{\bullet'}   (L \ltimes_{\ad} L))\\ \ra{\beta^\sharp} (G \ltimes_\t E)\ltimes_\bullet (E \ltimes_* (E \ltimes_{\t'} L) ) ,
\end{multline}
and lifted actions (now denoted by $\square$) of the first group on the remaining ones. Here $\bullet'$ denotes the secondary action of $\P({{\A}})$; lemma \ref{sa}, whereas $*$ denotes the derived actions either of ${{\A}}$ or of $\P(\A)$; lemma \ref{derived}.

There are four natural 2-crossed module maps $\P(\P({{\A}})) \to \P({{\A}})$. Namely the maps ${\rm Pr}^{\P({{\A}})}_1$  and ${\rm Pr}^{\P({{\A}})}_0$ of theorem \ref{mor}, which have the form (respectively):
\begin{multline*}
 (k,l,k',l')\longmapsto (kk',k'^{-1}  l k' l'), \quad \quad (a,e,k,a',e',k',l,l') \longmapsto (a,e,k)(a',e',k')(\de(l),1,l'), \\ (g,x,a,e,k) \longmapsto (g,x)(\partial(a),e)=(g\d(a), \d(a)^{-1} \t x \,\, e)
\end{multline*}
and
\begin{align*}
 (k,l,k',l')&\longmapsto (k,l), &(a,e,k,a',e',k',l,l') &\longmapsto (a,e,k),& (g,x,a,e,k) &\longmapsto (g,x).
\end{align*}
Also, by applying the path-space functor $\P$ to the 2-crossed module maps ${\rm Pr}^{{{\A}}}_1$ and  ${\rm Pr}^{{{\A}}}_0$, from $\P({{\A}})$ to ${{\A}}$, yields 2-crossed module maps  $\P({\rm Pr}^{{{\A}}}_1), \P({\rm Pr}^{{{\A}}}_0)\colon \P(\P({{\A}})) \to \P({{\A}})$, which have the form, respectively:
\begin{align*}
 (k,l,k',l')&\longmapsto (kl,k'l'), &(a,e,k,a',e',k',l,l') &\longmapsto (ae\de(k),a'e'\de(k'),ll'), &(g,x,a,e,k) &\longmapsto (g\d(x), a e \de(k)),\label{expl}\\
 (k,l,k',l')&\longmapsto (k,k'), &(a,e,k,a',e',k',l,l') &\longmapsto (a,a',l),& (g,x,a,e,k) &\longmapsto (g,a).
\end{align*}

Let us now change notation and put (all of these are maps $\P(\P({{\A}})) \to \P({{\A}})$):
\begin{align*}
&d_0=\P({\rm Pr}^{{{\A}}}_1); &d_1={\rm Pr}^{\P({{\A}})}_1; && &d_2 ={\rm Pr}^{\P({{\A}})}_0; &d_3=\P({\rm Pr}^{{{\A}}}_0).
\end{align*}
It is convenient to visualise the double path space in the form below,  emphasising the boundary maps to the path space (compare with \eqref{pathv}):
\begin{equation}\label{z1}
\xymatrix{\sqspace{{{\A}}}{{{\A}}}{{{\A}}}{{{\A}}}{\P{(\A)} }{\P{(\A)} }{\P({{\A}})} {\P{({{\A}})} }{\P(\P({{\A}}))}} 
\end{equation}
In conformity with this notation, we can denote $(a,e,k,a',e',k',l,l') \in  \Gr_1\big(\P(\P({{\A}}))\big)$ 
as:
\begin{equation}\label{z2}
\hskip-1cm\xymatrix{\sqt{a}{e}{k}{a'}{e'}{k'}{l}{l'}}
\end{equation}
and also $ (g,x,a,e,k)\in (G \ltimes_\t E)\ltimes_\bullet (E \ltimes_* (E \ltimes_{\t'} L) )=\Gr_0\big(\P(\P({{\A}}))\big)$ as:
\begin{equation}\label{z3}
\xymatrix{(g,x,a,e,k) =\\}\xymatrix{\squ{g}{x}{a}{e}{k}}  
\end{equation}
 Compare with the analogous notation for elements in the path-space in subsection \ref{tps}.

\subsubsection{Definition of the triangle space $\T({{\A}})$ of a 2-crossed module ${{\A}}$}\label{triangle}

Let $\G=(L \ra{\de} E \ra{\d} G, \t,\{,\})$ be a 2-crossed module. The triangle space $\T({{\A}})$ of  ${{\A}}$, will (not surprisingly) have a primary importance in the construction of the {concatenation} of 2-crossed module homotopies. This triangle space is included inside the double path space as the limit of the diagram: 
\begin{equation}
\xymatrix{&&\P(\P({{\A}}))  \ar[d]^{\P({\rm Pr}^{{{\A}}}_0)=d_3}\\ &{{\A}} \ar[r]_{i_{{{\A}}}} &\P({{\A}})}
\end{equation}
Looking at \eqref{z1}, \eqref{z2}, \eqref{z3}, we are thus making the left pointing upwards arrows to consist only of identities. 

Let us make this more explicit in dimensions one and two. Clearly:
\begin{align*} \Gr_1(\T({{\A}}))&=(E \ltimes_* (E \ltimes_{\t'} L))\ltimes_* \big((\{1\} \ltimes_* (E \ltimes_{\t'} L)\ltimes_{\bullet'}   (\{1\} \ltimes_{\ad} L)\big)\\
\Gr_0(\T({{\A}}))&= (G \ltimes_\t E)\ltimes_\bullet \big(\{1\} \ltimes_* (E \ltimes_{\t'} L) \big).
\end{align*}
The 2-crossed module boundary map  $\beta'\colon \Gr_1(\T({{\A}})) \to \Gr_0(\T({{\A}})) $ being:
\begin{equation}\label{bp}
\b'\big( (a,e,k), (1,f,l),(1,m)\big)=(\d(a),e,1,f,l).
\end{equation}
The restrictions of the projection maps $d_2={\rm Pr}_0^{\P({{\A}})}$ and $d_1={\rm Pr}_1^{\P({{\A}})}\colon\P(\P({{\A}})) \to \P({{\A}})$ to $\T({{\A}})$ are:
\begin{align}
d_2(a,e,k,1,f,l,1,m)&=(a,e,k), &  d_2(g,x,1,e,k)&=(g,x),\label{ppp}\\
d_1(a,e,k,1,f,l,1,m)&=\big(a,ef,(f^{-1} \t' k) l m  \big), &  d_1(g,x,1,e,k)&=(g,xe). \label{pp}
\end{align}
On the other hand the restrictions of $d_0=\P({\rm Pr}_1^{{\A}})$ and $d_3=\P({\rm Pr}_0^{{\A}})$ to $\T({{\A}})$ and   are given by, respectively:
\begin{align}
 d_0(a,e,k,1,f,l,1,m)& =(ae\de(k),f\de(l),m), & d_0(g,x,1,e,k)&=(g\d(x),e\de(k)), \label{def2v}\\
 d_3(a,e,k,1,f,l,1,m)& =(a,1,1), & d_3(g,x,1,e,k)&=(g,1).
\end{align}

As we did for the double path space, we will visualize the triangle space $\T({{\A}})$ in the form below \eqref{trianglev}, emphasising the three non-trivial boundary maps $\T({{\A}}) \to \P({{\A}}):$

\hskip-0.9cm\begin{minipage}{.4\linewidth}
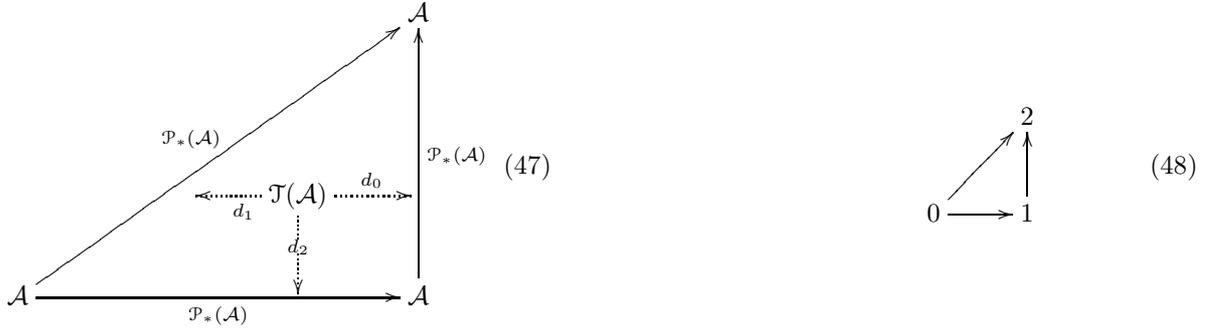
\begin{equation*}
\xymatrix{\trispace{{{\A}}}{{{\A}}}{{{\A}}}{\P{({{\A}})} }{\P{({{\A}})} }{\P{({{\A}})} }{\T({{\A}})}}
\end{equation*}
\end{minipage}
\begin{minipage}{.1\linewidth}
\begin{equation}\label{trianglev}
\end{equation}
\end{minipage}
\begin{minipage}{.2\linewidth} $$ $$
\end{minipage}
\begin{minipage}{.2\linewidth}
\begin{equation*}
 \xymatrix{&&2\\&0 \ar[ur]\ar[r]&1\ar[u] }
\end{equation*}
\end{minipage}
\begin{minipage}{.1\linewidth}
 \begin{equation}\quad\label{trinum}
 \end{equation}

\end{minipage}
In \eqref{trinum} we present the numbering of the vertices in \eqref{trianglev}.
The consistent notation for elements of $\Gr_0(\T({{\A}}))$ and of $\Gr_1(\T({{\A}}))$ is, respectively:
\begin{equation}\label{triel}
\xymatrix{\\(g,x,1,e,k)=\\}\xymatrix{\trixu{g}{x}{1}{e}{k}}
\end{equation}
and
\begin{equation}
\xymatrix{\\(a,e,k,1,e',k',1,l')=\\}\xymatrix{\trit{a}{e}{k}{1}{e'}{k'}{1}{l'}}
\end{equation}
The boundary map $\beta'\colon \Gr_1(\T({{\A}})) \to \Gr_0(\T({{\A}}))$ in the 2-crossed module $\T({{\A}})$ being:
\begin{equation}\label{hgf}
\xymatrix{\trit{a}{e}{k}{1}{e'}{k'}{1}{l'}}\xymatrix{\\ &\stackrel{\beta'}{\longmapsto}\\} \xymatrix{\trixu{\d(a)}{e}{1}{e'}{k'}}.
\end{equation}

 Consider the pull back 2-crossed modules $\P({{\A}}) {_{\Pr_1^{{\A}}}\times_{\Pr_0^{{\A}}}} \P({{\A}})$, $\P({{\A}}) {_{\Pr_1^{{\A}}}\times_{\Pr_1^{{\A}}}} \P({{\A}})$ and $\P({{\A}}) {_{\Pr_0^{{\A}}}\times_{\Pr_0^{{\A}}}} \P({{\A}})$. Looking at \eqref{trianglev} and \eqref{pathv}, we naturally denote these as (respectively):
$$\xymatrix{&&{{\A}}\\&{{\A}}\ar[r]_{\P({{\A}})} &{{\A}}\ar[u]_{\P({{\A}})}}  \quad \quad \quad \quad \quad \quad  \xymatrix{&&{{\A}}\\ &{{\A}}\ar[ru]^{\P({{\A}})} &{{\A}}\ar[u]_{\P({{\A}})}}  \quad \quad \quad \quad \quad \quad  \xymatrix{&&\A\\ &{{\A}}\ar[ru]^{\P({{\A}})}\ar[r]_{\P({{\A}})} &{{\A}}}$$
c.f. \eqref{pathv}.
We have obvious boundary maps from the triangle space $\T({{\A}})$ of ${{\A}}$ to each of these 2-crossed modules. Call these 2-crossed module maps $d_{+-},d_{++}$ and $d_{--}$. 
Now a simple, however essential, result:
\begin{Lemma}\label{techlemma}
 There exists a pull-back diagram of groups, where the boundary maps in the relevant crossed modules are denoted by $\beta'$:
$$\xymatrix{&&\Gr_1(\T(\A))\ar[dl]_{\beta'}\ar[rd]^{\Gr_1(d_{+-})} &\\ &\Gr_0(\T({{\A}}))\ar[rd]_{\Gr_0(d_{+-})} & &  \Gr_1\left(\P({{\A}}) {_{\Pr_1^{{\A}}}\times_{\Pr_0^{{\A}}}} \P({{\A}})\right)\ar[dl]^{\beta'}  \\ && \Gr_0\left(\P({{\A}}) {_{\Pr_1^{{\A}}}\times_{\Pr_0^{{\A}}}} \P({{\A}})\right) &}$$
The same statement is true for  $\P({{\A}}) {_{\Pr_1^{{\A}}}\times_{\Pr_1^{{\A}}}} \P({{\A}})$ and $\P({{\A}}) {_{\Pr_0^{{\A}}}\times_{\Pr_0^{{\A}}}} \P({{\A}})$, with the obvious adaptations. 
\end{Lemma}
\begin{Proof}
 This is a simple inspection. Clearly (by \eqref{qas} and \eqref{qaa}) the group  $\Gr_1\left(\P({{\A}}) {_{\Pr_1^{{\A}}}\times_{\Pr_0^{{\A}}}} \P({{\A}})\right)$ is given by all pairs of elements in 
$\Gr_1(\P({{\A}}))$ of the form: $$\big((a,e,k),(a e \de(k),e',k')\big), $$
being $\beta'\big((a,e,k),(a e \de(k),e',k')\big)=\big((\d(a),e),(\d(ae),e'\big),$ by \eqref{defab}.
 On the other hand the group pullback: $$\Gr_0(\T(\A)) {}_{\Gr_0(d_{+-})}\times_{\beta'} \Gr_1\left(\P({{\A}}) {_{\Pr_1^{{\A}}}\times_{\Pr_0^{{\A}}}} \P({{\A}})\right)$$ is given by all pairs of elements in $\Gr_0(\T({{\A}}))\times \Gr_1\left(\P({{\A}}) {_{\Pr_1^{{\A}}}\times_{\Pr_0^{{\A}}}} \P({{\A}})\right)$ of the form: 
$$\Big((\d(a),e,1,x,l), \big((a,e,k),(a e \de(k),x\de(l),k')\big)\Big). $$
Clearly the map:
$$\Gr_1(\T({{\A}})) \ra{\big(\beta',\Gr_1(d_{+-}) \big)} \Gr_0\big(\T({{\A}})\big) {}_{\Gr_0(d_{+-})}\times_{\beta'} \Gr_1\left(\P({{\A}}) {_{\Pr_1^{{\A}}}\times_{\Pr_0^{{\A}}}} \P({{\A}})\right), $$
made explicit below, is a group bijection:
$$(a,e,k,1,e',k',1,l') \longmapsto \Big ( (\d(a),e,1,e',k'), \big ( (a,e,k), (ae\de(k), e'\de(k'),l') \big) \Big).$$
\end{Proof}
\begin{remark}[Explicit operations in the triangle space]
By using \eqref{m5} and \eqref{gl} we can easily obtain an explicit form for the multiplication in $\Gr_0(\T(\A))=(G \ltimes_\t E)\ltimes_\bullet \big(\{1\} \ltimes_* (E \ltimes_{\t'} L) \big)$.
On the other hand the product in  $\Gr_1(\T(\A))=(E \ltimes_* (E \ltimes_{\t'} L))\ltimes_*((\{1\} \ltimes_* (E \ltimes_{\t'} L)\ltimes_{\bullet'}   (\{1\} \ltimes_{\ad} L))$ is:
\begin{multline}\label{pcomp}
\big( (a,e,k), (1,f,l),(1,m)\big)\big( (a',e',k'), (1,f',l'),(1,m')\big)\\=\Big((a,e,k)(a',e',k'), \big ( (\d(a'),e')^{-1} \bullet (1,f,l)\,\,  (1,f',l') \big),\\ (1,(f'\de(l'))^{-1} \t'
 \{(a' e' \de(k'))^{-1} ,(f \de(l))^{-1}\}^{-1}\big) \big( 1,(f'\de(l') )^{-1} \t'
 (a'e'\de(k'))^{-1} \t' m \big) (1,m') \Big).
\end{multline}
The restriction  $\triangle$ of the lifted action of $\P({{\A}})$ to an action of $\Gr_0(\T({{\A}}))$ on $\Gr_1(\T({{\A}}))$ is:
\begin{multline}\label{acomp}
 (g,x,1,z,w)\triangle \big ((a,e,k),(1,f,l),(1,m)\big)
=(g,x) \bullet \Big( (a,e,k), (\d(a),e)^{-1} \bullet (1,z,w) \,\, (1,f,l) \,\, (1,z,w)^{-1},\\  (1, (z \de(wl^{-1}) f^{-1}) \t' ( \{ (a e \de(k))^{-1},\de(w)^{-1}z^{-1}\}^{-1} ) \,\,\{z \de(w), (ae\de(k) f \de(l))^{-1}\} \,\,\d(z) \t m \big)\Big)
\end{multline}
None of this formulae will be explicitely used.
\end{remark}

\subsubsection{Definition of the  disk space $\D({{\A}})$ of a  2-crossed module ${{\A}}$}\label{disk}
Let ${{\A}}$ be a 2-crossed module. Recall the construction of the triangle space $\T({{\A}})$, the path-space $\P(\A)$ and also \eqref{trianglev}.
The disk space $\D({{\A}})$ of $\A$ is defined as being the limit of the diagram:  
$$\xymatrix{&\T({{\A}})\ar[d]^{d_0}&\\ &\P({{\A}}) &{{\A}}\ar[l]_<<<<<{i_{{\A}}}} $$
We are thus making the pointing upwards arrow of \eqref{trianglev} to consist only of identities. 
On the contrary of the triangle space, the operations on the disk space are quite simple. Let ${{\A}}=(L \ra{\de} E \ra{\d} G, \t,\{,\})$ be a 2-crossed module. The  disk space $\D({{\A}})$ is an embedded  2-crossed module, within the double path space $\P(\P({{\A}}))$ of ${{\A}}$. In the last bit of \eqref{aa}, instead of $E \ltimes_* (E \ltimes_{\t'} L) $, we put the subgroup isomorphic to $L$, of elements of the form $(1,\de(k),k^{-1})$, where $k\in L$. The group law is $$(1,\de(k),k^{-1}) (1,\de(l),l^{-1})=(1,\de(kl),(kl)^{-1});$$ see \eqref{a2}}.  Under the identification $k=(1,\de(k),k^{-1})$, the restriction of the first lifted action  $\bullet$ of ${{\A}}$ to this subgroup is, by \eqref{referright}:
\begin{equation}\label{lkjh}
(g,x) \bullet k=g \t (x \t' k).
\end{equation}

As far as the second group  of the disk $\D({{\A}})$   is concerned, we consider the subgroup of $(E \ltimes_* (E \ltimes_{\t'} L))\ltimes_{\bullet'}   (L \ltimes_{\ad} L)$, isomorphic to $L$, of elements of the form:
$$(1,\de(k),k^{-1},1,1). $$

By equations \eqref{Pp} and \eqref{v} , we have, where $|,|$ is the Peiffer lifting in $\P({{\A}})$: $$|(a,e,k),(1,\de(l),l^{-1})|=(1_L,1_L)=|(1,\de(l),l^{-1}),(a,e,k)|.$$
This will be used several times in the following calculations.
The restriction of the derived  action $*$ of $ E \ltimes_* (E \ltimes_{\t'} L))$ on $(E \ltimes_* (E \ltimes_{\t'} L))\ltimes_{\bullet'}   (L \ltimes_{\ad} L)$ to this isomorphic image of $L$ is,  given by $$(a,e,l) *\big (1,\de(k),k^{-1},1,1\big)=\big( \beta(a,e,l) \bullet (1,\de(k),k^{-1}),  1,1\big)=\big( 1,\d(a) \t (e  \de(k) e^{-1}) , \d(a) \t (e \t' k^{-1}),  1,1\big).$$
Thus under the identification $k=(1,\de(k),k^{-1},1,1)$ we have $$(a,e,l)*k=\d(a) \t (e \t' k).$$

We now describe the restriction of the lifted action $\square$ of $\P({{\A}})$. By \eqref{ker}:
\begin{multline}
(1,\de(k),k^{-1}) \square \big (a,e,l,1,\de(m),m^{-1},1,1\big)\\=\Big(a,e,l,1, e^{-1} \,\, ( \d(a)^{-1} \t \de(k) ) \,\,e\,\, \de(m)\,\, \de(k)^{-1}, k m^{-1} e^{-1} \t' ( \d(a)^{-1} \t k^{-1} ) ,1,1\Big).
\end{multline}

 Denoting the Peiffer lifting in $\D(\A)$ and $\P(\A)$ by $\{,\}$ and $|,|$, respectively, we have,  by \eqref{Pp}:
\[
\{(a,e,k,1,\de(l),l^{-1},1,1), (a',e',k',1,\de(l'),l'^{-1},1,1)\}
=\big(|aek,a'e'k'|,|aek,a'e'k'|^{-1} | aek\de(l) l^{-1},a'e'k'\de(l') l'^{-1}| \big),
\]
where we used the notation of \eqref{mnn} and \eqref{lkj}. Now note that, by \eqref{a1}, where $\{,\}$ is the Peiffer lifting in ${{\A}}$:
\begin{align*}
| aek\de(l) l^{-1},a'e'k'\de(l') l'^{-1}| &=| ae\de(l)\, (\de(l^{-1}) \t' k)\, l^{-1},a'e'\de(l') \,(\de(l')^{-1} \t' k' )\, l'^{-1}|\\
 &=(\{a,a'\},\{a,a'\}^{-1}\{ ae\de(l) \de(\de(l^{-1}) \t' k l^{-1}),a'e'\de(l') \de(\de(l')^{-1} \t' k' l'^{-1})|\\
 &=(\{a,a'\},\{a,a'\}^{-1}\{ ae\de(k),a'e'\de(k')\})=|aek,a'e'k'|.
\end{align*}
Therefore
$$\{(a,e,k,1,\de(l),l^{-1},1,1), (a',e',k',1,\de(l'),l'^{-1},1,1)\}
=\big(|aek,a'e'k'|,1,1\big) $$
We thus have the following  essential theorem:
\begin{Theorem}
 Given a 2-crossed module ${{\A}}=(L \ra{\de} E \ra{\d} G, \t,\{,\})$, there exists a 2-crossed module $\D({{\A}})$, called the  disk-space of ${{\A}}$, with underlying complex of groups:
$$L \ltimes_{\ad} L \ra{\alpha^2} (E \ltimes_* (E \ltimes_{\t'} L))\ltimes_* L \ra{\beta^2} (G \ltimes_\t E)\ltimes_\bullet L,$$
where (recall remark \ref{sd}):
\begin{align}
(g,e) \bullet k&= g \t (e \t' k) &(a,e,k) *l &=\d(a) \t (e \t' l).
\end{align}
The underlying action (denoted by $\square$)  of $(G \ltimes_\t E)\ltimes_\bullet L$ on  $(E \ltimes_* (E \ltimes_{\t'} L))\ltimes_* L$ is given by:
\begin{align}\label{l1}
(g,e) \square (a,f,l,l')&=\big((g,e) \bullet (a,f,l), g \t (e\t' l')\big), \textrm{ where } (g,e)=(g,e,1_L) \in (G \ltimes_\t E)\ltimes_\bullet L \\ \label{l2}
k \square (a,e,l,l') &=(a,e,l,  e^{-1} \t' ( \d(a)^{-1} \t k ) \,\, l'\,\, k ^{-1}), \textrm{ where } k=(1_G,1_E,k) \in (G \ltimes_\t E)\ltimes_\bullet L .
\end{align}
The action $\square$ of    $(G \ltimes_\t E)\ltimes_\bullet L$  on $L \ltimes_{\ad} L$ has the form:
\begin{equation}\label{mnb}
(g,e,l) \square (k,k')=(g,e) \bullet (k,k').
\end{equation}
The boundary maps are:
$$\alpha^2(k,l)=(\de(k),1,l,1) \an \beta^2(a,e,k,l)=(\d(a),e,l). $$ 
And finally the  Peiffer lifting is (recall \eqref{Pp}):
\begin{equation}\label{PLDS}
\begin{split}
\{(a,e,k,l),((a,e,k,l)\}&= |(a,e,k),(a',e',k')|=\big (\{a,a'\},\{a,a'\}^{-1} \{ae\de(k),a'e'\de(k')\} \big)  .
\end{split}
\end{equation}
The disk space $\D({{\A}})$ has  an inclusion map into $\P(\P({{\A}}))$, factoring through $\T({{\A}})$, with the form:
\begin{align*}
&(g,x,k) \longmapsto (g,x,1,\de(k),k^{-1}),
&(a,e,k,l)\longmapsto (a,e,k,1,\de(l),l^{-1},1,1), 
&&(k,l)\longmapsto (k,l,1,1) .
\end{align*}
Moreover the maps $d_2=(p,q,r)\colon \D({{\A}}) \to \P({{\A}})$ and $d_1=(p',q',r')\colon \D({{\A}}) \to \P({{\A}})$, where:
\begin{align}
p(k,l)&=(k,l), &q(a,e,k,l)&=(a,e,k), & r(g,e,k)&=(g,e),\label{poi1}\\
p'(k,l)&=(k,l), &q(a,e,k,l)&= (a,e\de(l),l^{-1}  k),  & r'(g,e,k)&=(g,e\de(k)),\label{poi3}
\end{align}
are morphisms of 2-crossed modules, and are obtained from the composition of the inclusion map $\D({{\A}}) \to \P(\P({{\A}}))$ and the projection maps $d_2={\rm Pr}_1^{\P({{\A}})},d_1={\rm Pr}_0^{\P({{\A}})}\colon \P(\P({{\A}})) \to \P({{\A}})$; see \eqref{pathv}.
Note: $$(a,e,k)(1,\de(l),l^{-1})=(a,e\de(l),(\de(l)^{-1} \t k) \,l^{-1})=(a,e\de(l),l^{-1} k ).$$
\end{Theorem}
{Note that the maps $d_1,d_2\colon \D(\A) \to \P({{\A}})$ can (compare with \eqref{pathv} and \eqref{trianglev}) be visualised as:}
\begin{equation}\xymatrix{\diskspace{{{\A}}}{{{\A}}}{\P({{\A}})}{\P({{\A}})}{\D({{\A}})}}
\end{equation}

\subsubsection{The tetrahedron space of a 2-crossed module}\label{tetrahedron}
Let ${{\A}} =(L \ra{\de} E \ra{\d} G, \t,\{,\})  $ be a 2-crossed module. We freely use the notation of \ref{triangle} and \ref{path}. The tetrahedron space $\Tet(\A)$ of ${{\A}}$ (or more precisely its underlying group $\Gr_0(\Tet(\A)$) will be used for proving that the concatenation of homotopies is associative, and several other places.
Consider the 2-crossed module $\P(\T({{\A}}))$, the path-space of the triangle space of $\A$. Its underlying group is: 
$$\Gr_0(\T({{\A}})) \ltimes_\Delta \Gr_1(\T({{\A}}))=\big( (G \ltimes_\t E)\ltimes_\bullet (\{1\} \ltimes_* (E \ltimes_{\t'} L) )\ltimes_\Delta\big( E \ltimes_* (E \ltimes_{\t'} L))\ltimes_*((\{1\} \ltimes_* (E \ltimes_{\t'} L)\ltimes_{\bullet'}   (\{1\} \ltimes_{\ad} L))\big) $$

By looking at \eqref{trianglev}, we have three maps $\P(\T({{\A}})) \to \P(\P({{\A}}))$ being $\P(d_i), i=0,1,2$, and two additional maps $\Pr_0^{\T({{\A}})},\Pr_1^{\T({{\A}}))}\colon \P(\T({{\A}})) \to \T({{\A}})$; see  theorem \ref{mor} and definition \ref{ppsp}. This corresponds to the five faces of a triangular prism. 
Let $\Tet({{\A}})$ be the 2-crossed module which is the limit of the diagram:
$$ \xymatrix{&&\P(\T({{\A}}))\ar[d]^{\P(d_2)}\\ 
              &\P({{\A}})\ar[r]_{i_{\P({{\A}})}} &\P(\P({{\A}})) }$$
We will only make use of the underlying group $\Gr_0(\Tet({{\A}}))$ of it, which is: 
$$\big ( (G \ltimes_\t E)\ltimes_\bullet (\{1\} \ltimes_* (E \ltimes_{\t'} L) )\big ) \ltimes_{\triangle} \big (\{1\} \ltimes_* (\{1\}\ltimes_{\t'} \{1\})\ltimes_*((\{1\} \ltimes_* (E \ltimes_{\t'} L)\ltimes_{\bullet'}   (\{1\} \ltimes_{\ad} L)\big).$$
The restrictions of the maps $\P(d_0),\P(d_1),\Pr_0^{\T(A)}$ and $\Pr_1^{\T(A)}$ to $\Gr_0(\Tet({{\A}}))$ are, respectively, and renaming them with a  simplicial notation:
\begin{align}\label{morB}
 \big((g,x,1,z,w),(1,1,1,1,f,l,1,m) \big)&\stackrel{d_1}{\longmapsto} \big(g,xz,1,f,lm\big)\\
\label{morC} \big((g,x,1,z,w),(1,1,1,1,f,l,1,m) \big) &\stackrel{d_0}{\longmapsto} \big(g\d(x),z\de(w),1,f\de(l),m\big)\\
\label{morA}
 \big((g,x,1,z,w),(1,1,1,1,f,l,1,m) \big)&\stackrel{d_3}{\longmapsto} \big(g,x,1,z,w\big)\\\label{morD}
 \big((g,x,1,z,w),1,1,1,1,f,l,1,m)\big) &\stackrel{d_2}{\longmapsto} \big(g,x,1,zf,f^{-1} \t' w \,\, l\big)
\end{align}
As we did for the triangle and disk space, we will represent the elements of $\Gr_0(\Tet({{\A}}))$ graphically, emphasising their boundary maps to $\Gr_0(\T(\A))$. We therefore put:
\begin{equation}\label{tetnot}\xymatrix{\\\big((g,x,1,z,w),(1,1,1,1,f,l,1,m) \big) =\\}\xymatrix{\te{g}{x}{1}{z}{w}{f}{l}{1}{m} }
\end{equation}
(we added squares to the graphical notation \eqref{triel} for elements of $\Gr_0(\T(\A))$, in order to emphasise the elements of $L$ assigned to the faces.)
We consider the following numbering of the vertices of the tetrahedron:
$$\xymatrixrowsep{6pt}
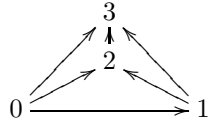
\xymatrix{&&3\\ &&2\ar[u]&\\&0\ar[rr]\ar[ru]\ar[ruu]&&1\ar[lu]\ar[luu]}  $$ Equation \eqref{tetnot} leaves out the $d_3$ boundary of an element of $\big((g,x,1,z,w),(1,1,1,1,f,l,1,m) \big) \in \Gr_0(\Tet({{\A}}))$,  which is:
\begin{equation}\label{mordd}
 \xymatrix{\trixumod{g}{x}{1}{zf}{f^{-1} \t' w \,\,l}}
\end{equation}

\section{Pointed homotopy of 2-crossed module maps}
\subsection{Quadratic  derivations and 1-fold homotopy of 2-crossed modules}\label{qder}
Suppose we have two pre-crossed modules $(E \to G,\t)$ and $(E' \to G',\t)$. 
\begin{Definition}\label{derdef}
Let $\f\colon G \to G'$ be a group morphism. A $\phi$-derivation $s\colon G \to E'$ is a set map such that, for each $g,h \in G$, we have:
$$s(gh)=\phi(h)^{-1} \t s(g)\,\,s(h). $$
\end{Definition}
Note that if $s\colon G \to E'$ is a derivation then:
\begin{equation}\label{deri}
 s(1_G)=1_{E'} \an s(g^{-1})=\phi(g) \t s(g)^{-1}.
\end{equation}

\begin{remark}\label{derprop}
 By looking at remark \ref{sd}, $\phi$-derivations are\footnote{This was pointed out to us by Ronnie Brown.}  in one-to-one correspondence with group maps $G \to G' \ltimes_{\t} E'$, where $g \longmapsto (\phi(g),s(g))$. In particular if $G$ is free, a $\phi$ derivation $s\colon G \to E'$ can be specified (and uniquely) by its {values on} a free basis $B$ of $G$. Therefore a set map $s_0\colon B \to E'$ uniquely extends to a $\phi$-derivation $s\colon B \to E'$.   We will very frequently use this line of thinking.
\end{remark}

Let $$\A=\left (L \ra{\de} E \ra{\d} G,  \t,\{,\}\right) \an \A'=\left (L' \ra{\de} E' \ra{\d} G', \t,\{,\}\right) $$
be 2-crossed modules. Let also $f=({\mu}, \psi,\phi)\colon \A \to \A'$ be a 2-crossed module morphism.
\begin{Definition}
A pair $(s,t)$ of maps 
$s\colon G \to E'$ and $t\colon E \to L'$ will be called a {\it quadratic $f$-derivation} if $s$ and $t$ satisfy, for each $g,h \in G$ and $a,b \in E$:
\begin{equation}\label{h1} s(gh)=\left(\f(h)^{-1} \t s(g)\right) s(h),
\end{equation}
\begin{equation}\label{h2} t(ab)=\big(\p(b)((s \circ \d))(b)\big)^{-1} \t' \left ( \left \{\psi(b), \f(\d (b)^{-1}) \t s (\d (a))^{-1}\right \}t(a) \right )t(b),
 \end{equation}
or alternatively, \eqref{iNvv}:
\begin{equation}\label{h3alt} t(ab)=\big((s \circ \d) (b) \big)^{-1} \t' \left ( \left \{\psi(b)^{-1},  s (\d (a))^{-1}\right \}^{-1}  \,\, \p(b)^{-1} \t' t(a) \right )\,\, t(b),
 \end{equation}
 and also (for each $g\in G$ and $a \in E$):
\begin{multline}\label{h3}
  t(g \t a)=\f(g) \t \Big (s(g) s(\d(a))^{-1} \t' \Big \{\p (a)^{-1},  s(g)^{-1})\Big\}^{-1} \Big) \\  \f(g) \t \Big \{ s(g), s(\d(a))^{-1}\p(a)^{-1} \Big\}   \big(\f(g)(\d \circ s)(g)\big) \t t(a).
\end{multline}
\end{Definition}
\begin{Lemma}[Pointed homotopy of 2-crossed module maps]\label{connect}
In the condition of the previous definition, if $(s,t)$ is a quadratic $f$-derivation, and if we define
$f'=({\mu}',\psi',\phi')\colon \A \to \A'$  as:
\begin{align} {\mu}'(l)&={\mu}(l) \,\,(t \circ \de)(l), \wh l \in L \\
              \psi'(a)&=\psi(a)\,\,((s \circ \d))(a) \,\,(\de \circ t)(a), \wh a \in E\\ 
              \phi'(g)&=\phi(g)\,\,(\d \circ s)(g), \wh g \in G
\end{align} 
then $f'$  a morphism of 2-crossed modules $\A \to \A'$. In this case we put: 
$$f\ra{(f,s,t)} f', $$
and say that $(f,s,t)$ is a   homotopy (or 1-fold homotopy) connecting $f$ and $f'$.
\end{Lemma}
This is  proved by using the following lemma,  noting that $f'= {\rm Pr}_1^{\A'} \circ H$ and  $f={\rm Pr}_0^{\A'}  \circ H$; see theorem \ref{mor}.
\begin{Lemma}\label{esc}
 Given a 2-crossed module morphism $f=({\mu},\psi,\phi)\colon \A \to \A'$, the pair of  maps $t\colon E\to L'$ and $s\colon G \to E'$ is a quadratic $f$-derivation if, and only if, $H=(i_3,i_2,i_1)\colon \A \to \P(\A')$ is a 2-crossed module morphism, where
\begin{align*} 
l \in L &\stackrel{i_3}{\longmapsto} \big({\mu}(l),t \circ \de(l)\big) \in L' \ltimes_{\ad} L'\\ 
 a \in E  & \stackrel{i_2}{\longmapsto}  \big(\p(a),(s \circ \d)(a),t(a)\big) \in E' \ltimes_* (E' \ltimes_{\t'} L') \\  g \in G &\stackrel{i_1}{\longmapsto}  \big(\f(g),s(g)\big) \in G' \ltimes_\t E'.
\end{align*}
Also, $H=(i_3,i_2,i_1)$ is a 2-crossed module map if, and only if, $(i_2,i_1)$ is a pre-crossed module map.
\end{Lemma}
\begin{Proof} Conditions (\ref{h1}), (\ref{h2}), (\ref{h3}) express exactly that $(i_2,i_1)$ is a pre-crossed module morphism.
 That $H$ defines a morphism of complexes $\A \to \P(\A')$, which is equivariant with respect to the actions of $G$ and $G'\ltimes_\t E'$, if and only if, $(t,s)$ is a quadratic $f$- derivation follows from the explicit construction of $\P(\A')$. We now need to prove that $H$ always preserves the Peiffer lifting, which follows immediately from the equation:
\begin{equation}\label{mnnbv}
t\big(\langle a,b\rangle\big)=\big\{\psi(a),\psi(b)\big\}^{-1}\big\{\psi(a) \, s(\d(a)) \, \de(t(a)), \psi(b) \, s(\d(b)) \, \de(t(b))\big\}, \textrm{ for each } a,b \in E.
\end{equation}
This equation was proven in \cite{FM3}.
\end{Proof}
\begin{remark}\label{wisk} Note that if $(s,t)$ is a quadratic $f$ derivation, where $f\colon  \A \to \A'$, connecting $f$ and $f'$, and if $g\colon \A' \to \A''$ is 2-crossed module map, then $(g\circ s, g\circ t)$ is a quadratic $(g \circ f)$-derivation connecting $g \circ f$ and $g\circ f'$. Analogously if $h\colon \A''' \to \A$ is a 2-crossed module map then $(s\circ h, t\circ h)$ is a quadratic $(f\circ h)$-quadratic derivation connecting $f\circ h$ and $f'\circ h$.
\end{remark}

\begin{remark}\label{visform}
Given $f,f'\colon \A \to \A'$, consider a homotopy $f=(\mu,\psi,\eta)\ra{(f,s,t)} (\mu',\psi',\eta')=f'$. The underlying maps $G\to \Gr_0(\P({{\A}}'))$ and $E \to \Gr_1(\P({{\A}}'))$ can be visualised (see \eqref{pathv}) in  the form:
$$g \in G \longmapsto  \Big(\pu{\phi(g)\,}{s(g)}\Big) \an e \in E  \longmapsto \Big({\psi(e)}\ra{\big ( s(\d(e)), t(e) \big )} \psi(e) \,s(\d(e))\, \de(t(e))\Big ).  $$
\end{remark}
\subsubsection{{Homotopy between 2-crossed module maps is not an equivalence relation}}\label{fail}

As mentioned in the Introduction, given 2-crossed modules $\A$ and $\B$, and considering 2-crossed module maps $f,g\colon \A \to \B$,  saying that $f \sim g$ if there exists a quadratic $f$-derivation $(s,t)$ with $f \ra{(f,s,t)} g$, does not define an equivalence relation. This is in contrast with the crossed complex case treated in \cite{BH,BHS}, but it is however not surprising for model category theoretical reasons \cite{DS}, given that this should only be the case (in general) when $\A$ is cofibrant and $\B$ is fibrant. Note that by Theorem \ref{mor} and  Lemma \ref{esc}, homotopy between 2-crossed module maps is defined via a good path space in the model category of 2-crossed modules. 

Nevertheless, let us give an explicit example to show that 2-crossed module homotopy does not define an equivalence relation between 2-crossed module maps. In the example below, symmetry fails to hold.

Let $\E=2\Z$ be the group of even integers. Let $\Z_2=\Z/\E=\{0,1\}.$ Let $\d\colon \Z \to \Z_2$ be the quotient map. Let $\de\colon \E \to \Z$ be the inclusion. We have an action of $\Z_2$ on $\Z$ and on $\E$ by automorphisms, being $0 \t  n=n$ and $1 \t n=-n$, for each $n \in \Z$.

 Clearly $(\d\colon \Z \to \Z_2, \t)$ is a pre-crossed module.  The Peiffer pairing is (given $m,n \in \Z$): 
\[\langle m,n \rangle   = \left\{ \begin{array}{ll} 0 & ,\textrm{ if } m \textrm{ is even}
       \\ 2n &,\textrm{ if } m \textrm{ is odd }
         \end{array} \right. \]

Consider now the 2-crossed modules: $$\A=	\big(\{0\} \to \{0\} \to \Z_2 \big) \textrm{ and } \B=\big ( \E \ra{\de}  \Z \ra{\d} \Z_2,\t,\{,\}\big),$$
where for $\A$ all actions are trivial, as is  the Peiffer lifting, and the Peiffer lifting in $\B$ has the form:  $\{m,n\}=\langle m, n\rangle$, for $m,n \in \Z$. See Example \ref{gh} for a similar construction. 

Let us consider $f=(\eta,\psi,\phi)\colon \A \to \B$, such that $\eta$ and $\psi$ are trivial, and $\phi\colon \Z_2 \to \Z_2$ is the identity map. This clearly is a 2-crossed module map. Let us also consider $f'=(\eta,\psi,\phi')\colon \A \to \B$, such that $\eta$, $\psi$  and $\phi'$ all are trivial, thus their images consist solely of the identity element. 

We now consider the $f$-quadratic derivation $(s,t)$ such that $s(1)=1,s(0)=0$ and $t(0)=0$.  It is clear that  \eqref{h2} and \eqref{h3} are satisfied. As for condition \eqref{h1}, we check it case by case:
\begin{align*}
 0&=s(1+1)=(-\phi(1)) \t s(1) +s(1)=-s(1)+s(1)=0\\
 1&=s(0+1)=(-\phi(1)) \t s(0) +s(1)=-s(0)+s(1)=1\\
 1&=s(1+0)=(-\phi(0)) \t s(1) +s(0)=s(1)+s(0)=1\\
 0&=s(0+0)=(-\phi(0)) \t s(0) +s(0)=-s(0)+s(0)=0\\
\end{align*}

 We clearly have that (by Lemma \ref{connect}):
$$f \ra{(f,s,t)} f'. $$
Let us see that there does not exist an $f'$-quadratic derivation $(s',t')$ such that $f' \ra{(f',s',t')} f. $ Indeed, for this to happen, we would need a $\phi'$-derivation $s'\colon \Z_2 \to \Z$, such that $\d(s(1))=1$. However, since $\phi'(\Z_2)=\{0\}$, $s'$ would need to be (by \eqref{h1}) a non-trivial group morphism $\Z_2 \to \Z$, which is impossible.

\subsection{Quadratic 2-derivations and 2-fold homotopy of 2-crossed modules}\label{2folddef}
Let $$\A=\left (L \ra{\de} E \ra{\d} G,  \t,\{,\}\right) \an \A'=\left (L' \ra{\de} E' \ra{\d} G', \t,\{,\}\right) $$
be 2-crossed modules. 
Let  $f=(\mu,\psi,\phi)\colon \A \to \A'$ be a 2-crossed module morphism. Consider a quadratic  $f$-derivation $(s,t)$. Recall the construction of the  disk space $\D(\A')$ of $\A'$;  \ref{disk}.

\begin{Definition}
 We say that a map $k\colon G \to L'$ is a quadratic $(f,s,t)$ 2-derivation if for each $g,h \in G$:
\begin{equation}\label{2der} k(gh)=\Big (  s(h)^{-1} \t'\big (\phi(h)^{-1} \t k(g)\big) \Big)\,\, k(h).
\end{equation}
\end{Definition}
\begin{Lemma}\label{dereprop}
 The map $k\colon G \to L'$ is a quadratic $(f,s,t)$ 2-derivation if and only if 
$g \longmapsto(\f(g),s(g),k(g))$ is a group morphism $G \to (G \ltimes_ \t E)\ltimes_\bullet L=\Gr_0(\D(\A'))$.
In particular $k(1_G)=1_K$ and for each $g \in G$:
$$k(g^{-1})=\phi(g) \t (s(g) \t' k(g)^{-1}) .$$
Moreover if $g,h,i \in G$:
$$ k(g^{-1}hi)=\Xi^{(\phi,s,k)}(g,h,i),$$ where we have defined:
\begin{equation}\label{defxi}
\Xi^{(\phi,s,k)}(g,h,i)=
s(i)^{-1} \t'\Big(\phi(i)^{-1} \t \Big (  s(h)^{-1} \t'\big((\phi(h^{-1})\phi(g)) \t (s(g) \t' k(g)^{-1})\big) \,\, k(h) \Big) \Big)\,\,k(i) .
 \end{equation}
\end{Lemma}
\begin{Proof}
 The first assertion is immediate from remark \ref{sd} and  \eqref{lkjh}. Also:
\begin{align*}
 k(g^{-1}hi)&=s(i)^{-1} \t'\Big(\phi(i)^{-1} \t \big (  s(h)^{-1} \t'(\phi(h)^{-1} \t k(g^{-1})) \,\, k(h) \big) \Big)\,\,k(i) \\
&=s(i)^{-1} \t'\Big(\phi(i)^{-1} \t \big (  s(h)^{-1} \t'(\phi(h^{-1}g) \t (s(g) \t' k(g)^{-1})) \,\, k(h) \big) \Big)\,\,k(i) 
\end{align*}

\end{Proof}
\begin{Corollary}\label{extt}
 In the condition of the previous lemma, if $G$ is a free group on the basis $B$ a quadratic $(f,s,t)$ 2-derivation $k\colon G \to L'$ can be specified, uniquely, by its valued in the basis $B$. More precisely given a set map $k_0\colon B \to L'$ there exists a unique quadratic $(f,s,t)$ 2-derivation $k\colon G \to L'$ extending $k_0$.
\end{Corollary}

\begin{Lemma}
The set map $k\colon G \to L'$ is a quadratic $(f,s,t)$  2-derivation if, and only if, the map $H^2=(j_3,j_2,j_1)$, below, is a 2-crossed module morphism from $\A$ into $\D(\A')$.
\begin{align*}
 j_1(g)&=\big(\f(g),s(g),k(g)\big)\\
 j_2(e)&=\big(\psi(e),(s \circ \d)(e),t(e), (k \circ \d)(e) \big)\\
 j_3(l)&=\big(l,(t \circ \d(l)\big)
\end{align*}
\end{Lemma}
\begin{Proof}
  We already know by lemma \ref{esc} that, forgetting the last component of $j_1$ and $j_2$,  we have a 2-crossed module map $\A \to \P(\A')$. Given the form of $\D(\A')$ we clearly get a morphism of group complexes $\A \to \D(\A')$, which is compatible with the action of the first groups of the complexes on the remaining. To prove this note the following calculation, and compare with \eqref{l1}, \eqref{l2}, \eqref{lkjh} and \eqref{mnb}:
\begin{align*}
 k(\d(g \t e))=k(g\d(e)g^{-1})&=(\phi(g) \t s(g))\t' \big(\phi(g) \t k(g \d(e))\big) \,\, \phi(g) \t (s(g) \t' k(g)^{-1}) \\
&=\phi(g) \t\Big (s(g) \t' \big( k(g \d(e) )\,\, k(g)^{-1}\big) \Big) \\
&=\phi(g) \t\big (s(g) \t' \left(  \Big (  s(\d(e))^{-1} \t'\big (\phi(\d(e))^{-1} \t k(g)\big) \Big)\,\, k(\d(e)) \,\, k(g)^{-1}\right) \big).
\end{align*}

 We now   need to prove that $H^2$ preserves the Peiffer lifting. This follows from lemma \ref{esc} and the form of the Peiffer lifting on $\D(\A')$ and $\P(\A')$, see \eqref{Pp} and \eqref{PLDS}.
\end{Proof}

Therefore, looking at the maps $d_2,d_1\colon \D(\A') \to \P(\A')$ of \eqref{poi1}  and \eqref{poi3} (and their composition with $H^2\colon \A \to \D(\A')$), lemma \ref{esc} and the previous one, given a quadratic $(f,s,t)$ 2-derivation $k$ then $s'\colon G \to E'$ and $t'\colon E \to L'$, defined as:
\begin{align}\label{tar}
&s'(g)=s(g)\,( \de \circ k) (g), \quad \quad \quad  \textrm{ and }
&t'(e)=(k\circ \d)(e)^{-1} \,\, t(e)
\end{align}
yield a   quadratic  $f$-derivation $(s',t')$, and in this case we put:
$$(f,s,t) \ra{(f,s,t,k)} (f,s',t') .$$

\begin{remark}\label{coher}
Note that if
$f \ra{(f,s,t)} f' $
then we also have $f \ra{(f,s',t')} f'$.
\end{remark}

\begin{Definition}
 We say that the quadratic $f$-derivations $(s,t)$ and  $(s',t')$ are 2-fold homotopic if there exists a quadratic $(f,s,t)$ 2-derivation such that $(f,s,t) \ra{(f,s,t,k)} (f,s',t') .$ The quadruple $(f,s,t,k)$ will be called a 2-fold homotopy connecting  $(f,s,t)$ and  $(f,s',t')$.
\end{Definition}

\begin{remark}\label{tte}
Looking at the definition of the disk space 2-crossed module from the triangle space 2-crossed modules, a quadratic $(f,s,t)$ 2-derivation is a map $k\to L'$, such that  
$g\longmapsto \big(\phi(g),s(g), 1,\de(k(g)),k(g)^{-1}\big)$
is a group morphism $ G \to (G' \ltimes_\t E')\ltimes_\bullet (\{1\} \ltimes_* (E' \ltimes_{\t'} L') )=\Gr_0(\T(\A'))$. This also follows directly from \eqref{referright}. (We have put $f=(\mu,\psi,\phi)$.)  We can thus visualise a quadratic $(f,s,t)$ 2-derivation $k\colon G \to L'$ as a map:
$$\xymatrix{\\ g \in G \mapsto\\} \xymatrix{\trixumodb{\phi(g)}{s(g)}{ 1}{\de(k(g))}{k(g)^{-1}}} \xymatrix{\\ \in \Gr_0(\T(\A'))\\} $$
\end{remark}

\subsection{A groupoid of 2-crossed module maps and their homotopies (in the free up to order one case)}
In this subsection, let us fix two 2-crossed modules:
$$\A'=\left (L' \ra{\de} E'\ra{\d} F,  \t,\{,\}\right) \an \A=\left (L \ra{\de} E \ra{\d} G,  \t,\{,\}\right)  .$$
 Suppose also that $F$ is a free group, with a chosen (free) basis $B \subset F$.  Let us define a groupoid $[\A',\A]_1^B$, with objects the 2-crossed module maps $\A' \to \A$ and morphisms their homotopies. This will explicitly depend on the chosen basis $B$ of $F$. We will freely use the notation of subsections \ref{tps}, \ref{triangle} and \ref{qder}.

\subsubsection{Concatenating homotopies in the free up to order one case}\label{comphom}
Consider homotopies of 2-crossed module maps {$\A' \to \A$}:
$$f=({\mu},\psi,\phi) \ra{(f,s,t)} f'=({\mu}',\psi',\phi') \an ({\mu}',\psi',\phi')   \ra{(f',s',t')}  f''=({\mu}'',\psi'',\phi'')$$
(where $f,f',f''\colon \A' \to \A)$.
We therefore have pre-crossed module maps $H,H'\colon \A' \to \P(\A)$; lemma \ref{esc}. We consequentely have an induced map to the group pull-back: 
$$(H,H') \colon  \A' \to \P({{\A}}) {_{\Pr_1^{{\A}}}\times_{\Pr_0^{{\A}}}} \P({{\A}}). $$
Let us define the {concatenation} of homotopies:
$$f=({\mu},\psi,\phi) \ra{(f,s{\otimes} s',t{\otimes} t')} f''=({\mu}'',\psi'',\phi'').$$
It suffices to define the associated pre-crossed module map (lemma \ref{esc}) $H\tn H' =(f, s \tn s',t \tn t') \colon  \A' \to \P(\A)$.

The derivation {$(s {\otimes} s')\colon F \to E$} is the unique $\phi$-derivation (see definition \ref{derdef} and remark \ref{derprop})  which on the chosen basis $B$ of $F$  has the form \begin{equation} b \in B \longmapsto (s {\otimes} s')(b)=s(b) s'(b) \in E.\end{equation}

There is another piece of information that we will use,  namely a set map $\w^{(s,s')}\colon F \to L$, measuring the difference {(for each $g \in F$)} between $s(g)s'(g)$ and $(s\otimes s')(g)$; this difference is null in a crossed module. Recall the construction of the triangle space 2-crossed module  \ref{triangle}. Given that $F$ is free, there exists a unique group  map 
$X^{(s,s')}\colon F \to  (G \ltimes_\t E)\ltimes_\bullet (\{1\} \ltimes_* \big(E \ltimes_{\t'} L) \big)=\Gr_0(\T(\A))$, say: 
$$\xymatrix{\\g \in F \stackrel{X^{(s,s')}}{\longmapsto} \big(\phi(g),s(g), 1, \zeta(g),\w^{(s,s')}(g)\big)=\\}\xymatrix{\trixu{\phi(g)}{s(g)}{1}{ \zeta(g)}{\w^{(s,s')(g)}}}$$ which on the chosen basis $B$ of $F$ takes the form:
\begin{equation}\label{ghb}
\xymatrix{&\\ b \stackrel{X^{(s,s')}}{\longmapsto} \big(\phi(b),s(b),1,s'(b),1\big) = \\} \xymatrix{\trixuvb{\phi(b)}{s(b)}{1}{s'(b)}{1}} 
\end{equation}

Looking, see \eqref{trianglev}, at the three $\P({{\A}})$ sides of the triangle space $\T({{\A}})$, we have maps from $F$ to the underlying group $G' \ltimes_\t E'=\Gr_0(\P({{\A}}'))$   of the path space $\P({{\A}}')$. By remark \ref{derprop}, since there are unique group maps $F \to G' \ltimes E'$ extending 
$b\longmapsto (\phi(b), s(b)s'(b))$,  $b\longmapsto (\phi(b)\d(s(b)), s'(b))$ and $b\longmapsto (\phi(b), s(b))$, it follows at once that $X^{(s,s')}$ has the form:
$$\xymatrix{\\g \stackrel{X^{(s,s')}}{\longmapsto} \big(\phi(g),s(g), 1, \zeta(g),\w^{(s,s')}(g)\big)=\\}\xymatrix{\trixuu{\phi(g)}{s(g)}{1}{ \tn s'}{\w^{(s,s')}(g)}}$$
Therefore $(s {\otimes} s') (g)=s(g) \, \zeta(g),$ and 
 $\zeta(g)\, \de(\w^{(s,s')}(g))=s'(g)$, or
\begin{equation}\label{circ2}
 s(g)\,\zeta(g) \,\de(\w^{(s,s')}(g))= s(g)\, s'(g), \fo g \in F.
\end{equation}

 We have proven:
\begin{Lemma}\label{good}
 There exists a unique group morphism $X^{(s,s')}\colon F \to \Gr_0(\T(\A))$, with: 
$$\xymatrix{\\g \longmapsto \big(\phi(g),s(g),1, s'(g) \de(\w^{(s,s')}(g))^{-1}, \w^{(s,s')}(g)\big)=\\} \xymatrix{\trixul{\phi(g)}{s(g)}{1} {s'(g) \de(\w^{(s,s')}(g))^{-1}}{\w^{(s,s')}(g)}}  $$ 
for each $g$ in $F$, such that $\w(b)=1$ (see \eqref{ghb}) on the free generators $b \in B$ of $F$.
In particular by \eqref{circ2}:
\begin{equation}\label{q}(s {\otimes} s') (g)=s(g) \,s'(g) \, \de(\w^{(s,s')}(g))^{-1},\textrm{ for each } g \in F. \end{equation}
\end{Lemma}

{We now define $(t {\otimes} t')\colon E' \to L$.  For $e \in E'$ put:}
\begin{equation}\label{ttn}
(t {\otimes} t') (e)=\big(\w^{(s,s')}(\d(e)) \big )  \, \big( s'(\d(e))^{-1} \t' t(e)\big) \,\, t'(e).
\end{equation}
We prove that $(s{\otimes} s', t {\otimes} t')$ is a quadratic $({\mu},\psi,\phi)$-derivation. 
Consider the commutative triangle: $$\xymatrix{ & F\ar[r]^{X^{(s,s')}} \ar[rd]_{\Gr_0(H,H')} &\Gr_0(\T(\A))\ar[d]^{\Gr_0(d_{+-})}  \\ &&  \Gr_0\left(\P({{\A}}) {_{\Pr_1^{{\A}}}\times_{\Pr_0^{{\A}}}} \P({{\A}})\right)}$$ Define $Y^{(s,s')}\colon E' \to \Gr_1(\T(\A))$ as being the unique group map  making the diagram below commutative:
 \begin{equation}\label{diag}
\xymatrixcolsep{4pc}
\xymatrix{ & & E'\ar@{-->}[r]^{Y^{(s,s')}} \ar[dll]_\d \ar[dr]|<<<<<<<<<<<<\hole_{\Gr_1(H,H')} &\Gr_1(\T(\A))\ar[d]^{\Gr_1(d_{+-})}  \ar[dll]_{\beta'}   \\ F\ar[r]_{X^{(s,s')}} \ar[rd]_{\Gr_0(H,H')} &\Gr_0(\T(\A))\ar[d]^{\Gr_0(d_{+-})}  & &\Gr_1\left(\P({{\A}}) {_{\Pr_1^{{\A}}}\times_{\Pr_0^{{\A}}}} \P({{\A}})\right) \ar[dll]_{\beta'}\\ &  \Gr_0\left(\P({{\A}}) {_{\Pr_1^{{\A}}}\times_{\Pr_0^{{\A}}}} \P({{\A}})\right)}
\end{equation} (We are using lemmas \ref{cart} and \ref{techlemma}.)
In particular $(X^{(s,s')},Y^{(s,s')})$ will be a pre-crossed module map. By composing with $d_1$ in \eqref{trianglev} will give $H\tn H' =(f, s \tn s',t \tn t') \colon  \A' \to \P(\A)$.

To determine $Y^{(s,s')}$ we only need (remark \ref{car2}) to find a set map that makes the diagram \eqref{diag} commutative.  Looking  at the construction of the triangle pre-crossed module (subsection \ref{triangle}),
 consider the set map:
 $$Y^{(s,s')} \colon E' \to (E \ltimes_* (E \ltimes_{\t'} L))\ltimes_*((\{1\} \ltimes_* (E \ltimes_{\t'} L)\ltimes_{\bullet'}   (\{1\} \ltimes_{\ad} L))=\Gr_1(\T(\A)),$$ definition \ref{princ}, of the form:
$$e \stackrel{Y^{(s,s')}}{\longmapsto} \Big(\psi(e),s (\d(e)), t(e), 1,s'(\d(e))\,\, \de(\w^{(s,s')} (\d(e)))^{-1},\w^{(s,s')} (\d(e)), 1,t'(e)\Big). $$
That is:
$$ Y^{(s,s')}(e)=\xymatrix{\tritl{\psi(e)}{s (\d(e))}{ t(e)}{ 1}{s'(\d(e)) {\de(\w^{(s,s')} (\d(e)))^{-1}}}{\w(\d(e))}{ 1}{t'(e)}}$$
Note that $\b'\circ Y^{(s,s')}=X^{(s,s')}\circ \d,$ where $\beta\colon \Gr_1(\T(\A)) \to \Gr_0(\T(\A))$ is the boundary map; \eqref{bp} or \eqref{hgf}. Clearly this makes the diagram \eqref{diag}  commute thus $(Y^{(s,s')},X^{(s,s')})$ is a pre-crossed module map.

In particular, by looking at the oblique side of $Y^{(s,s')}$ and $X^{(s,s')}$ it follows by remark \ref{esc} and \eqref{trianglev} that $(s{\otimes} s', t {\otimes} t')$ is  a quadratic $({\mu},\psi,\phi)$-derivation.

\begin{Lemma} 
 Consider homotopies of 2-crossed module maps $$f=({\mu},\psi,\phi) \ra{(f,s,t)} f'=({\mu}',\psi',\phi') \an ({\mu}',\psi',\phi')   \ra{(f',s',t')}  f''=({\mu}'',\psi'',\phi'').$$
Then the $({\mu},\psi,\phi)$-quadratic derivation $\big( (s{\otimes} s'), (t{\otimes} t')\big) $ connects 
$f$ and $f''$.
\end{Lemma}
\begin{Proof}At the pre-crossed module level follows by construction. Let us give full details however.
 Let $\bar{f}=(\bar{{\mu}},\bar{\psi},\bar{\phi})$ be the 2-crossed module morphism defined from $f$ and the quadratic derivation $(s{\otimes} s',t{\otimes} t')$; lemma \ref{connect}. We must prove that $\bar{f}=f''$.
Let $b$ be a free generator of $F$. Then 
$$\bar{\phi}{(b)}=\phi(b) \, \d(s {\otimes} s'(b))= \phi(b) \, \d(s(b) \, s'(b))=\phi(b) \, \d(s(b) ) \, \d(s'(b)))=\phi'(b)\, \d(s'(b))=\phi''(b).$$
In particular it follows that  $\bar{\phi}{(g)}=\phi''(g)$, for each $g \in F$. 

{Given $e \in E'$ we have (we use \eqref{q} and the  rule $\de(e \t' k)=e \de(k) e^{-1}$, for each $e \in E $ and $k \in L$):}
\begin{align*}
 \bar{\psi}(e)&=\psi(e)\, (s {\otimes} s')(\d(e)) \,\, \de((t {\otimes} t')(e))\\
              &=\psi(e)  \, s(\d(e)) \,s'(\d(e)) \, \de(\w(\d(e)))^{-1} \,\, \de \Big(\big((\w \circ \d)(e) \big )  \, \big( s'(\d(e))^{-1} \t' t(e)\big) \, t'(e)\Big)\\
              &=\psi(e)  \, s(\d(e)) \,\de(t(e)) \, \, s'(\d(e)) \,\de(t'(e))=\psi'(e)\, \, s'(\d(e)) \,\de(t'(e))
              =\psi''(e).
\end{align*}
Finally, given $k\in L'$ we have (since $\d \circ \de(k)=1$ for each $k \in L$):
$$\bar{{\mu}}(k)={\mu}(k) \, (t {\otimes} t')(\de(k))= {\mu}(k) \,t(\de(k)) \, t'(\de(k))={\mu}''(k).$$
Note $s'(1_F)=1_{E}$ and $\w(1_F)=1_{L}$.
\end{Proof}

\begin{remark}[Some properties of $\w^{(s,s')}\colon F \to L$]\label{propomega}
Clearly 
\begin{align*}
&\w^{(s,s')}(1_F)=1_L  & \textrm{and }&   &\w^{(s,s')}(b)=1_L, \textrm{ for each } b \in B.
\end{align*}
Also, by remark \ref{sd} and equation \eqref{m5}, for each $g,h\in F$:
\begin{equation*}
 \w^{(s,s')}(gh)=\w^{(s,s')}(h)\,\,\, s'(h)^{-1} \t' \Big( \f(h)^{-1} \t \big\{\phi(h) \t s(h)^{-1}, \de(\w^{(s,s')}(g)) s'(g)^{-1}\big\} \,\, \phi'(h)^{-1} \t \w^{(s,s')}(g) \Big)
\end{equation*}
and
\begin{equation*}
 \w^{(s,s')}(g^{-1})=\phi(g)\t\big\{s(g),s'(g)(\w^{(s,s')}(g))^{-1}\big\}\,\,\phi'(g)\t\big(s'(g) \t'(\w^{(s,s')}(g))^{-1}\big).
\end{equation*}
In particular, if $b,b',b''$ are free generators of $F$:
\begin{align*}
 &\w^{(s,s')}(bb')= s'(b')^{-1} \t' \big( \{ s(b')^{-1}, \phi(b')^{-1} \t s'(b)^{-1}\} \big),
&\textrm{and }\quad \quad \w^{(s,s')}(b^{-1})&=\phi(b)\t\big\{s(b),s'(b)\big\}.
\end{align*}
Thus: $$\w^{(s,s')}(b^{-1}b'b'')=\Theta^{(s,s')}(b,b',b''),$$
where by definition:
\begin{multline}\label{deftheta}
 \Theta^{(s,s')}(b,b',b'')=  s'(b'')^{-1} \t' \big( \{ s(b'')^{-1}, \phi(b'')^{-1} \t s'(b')^{-1}\} \big)\\ \big(\phi'(b'')^{-1} \t s'(b')\,\,s'(b'')\big)^{-1} \t' \Big(  \left\{     \big(\phi(b'')^{-1} \t s(b')\,\,s(b'')\big)^{-1},( \f(b'b'')^{-1} \phi(b)) \t \big (s(b)\, s'(b) \, s(b)^{-1}\big)\right\} \\ \big (\phi'(b'b'')^{-1} \phi(b)\big) \t \left\{s(b),s'(b)\right\} \Big).
\end{multline}

\end{remark}

\subsubsection{The {concatenation} of homotopies is associative (in the free up to order 1 case)}
We freely use the notation of subsection \ref{triangle}, and we resume the notation and context of \ref{comphom}.

\begin{Proposition}
 The {concatenation} of homotopies is associative.
\end{Proposition}
\begin{Proof}
Choose a chain of homotopies of 2-crossed module maps $\A' \to \A$:
 $$f=({\mu},\psi,\phi) \ra{(f,s,t)} f'=({\mu}',\psi',\phi')   \ra{(f',s',t')}  f''=({\mu}'',\psi'',\phi'')  \ra{(f'',s'',t'')}  f'''=({\mu}''',\psi''',\phi''').$$
It is immediate that $(s{\otimes} s') {\otimes} s''=s{\otimes}( s' {\otimes} s'')$, since this is true in a free basis of $F$. 

Let us now see that $(t{\otimes} t') {\otimes} t''=t{\otimes}( t' {\otimes} t'')$. We now have, {for each $e \in E'$:}
$$(t \otimes t')(e)= \big((\w^{(s,s')} \circ \d)(e) \big )  \, \big( s'(\d(e))^{-1} \t' t(e)\big) \, t'(e),$$ 
$$(t' {\otimes} t'')(e)= \big((\w^{(s',s'')} \circ \d)(e) \big )  \, \big( s''(\d(e))^{-1} \t' t'(e)\big) \, t''(e),$$ 
\begin{align*}
 \big(t {\otimes} (t' {\otimes} t'')\big)(e)&=(\w^{(s,s' \tn s'')} \circ \d)(e) \,\, (s'{\otimes} s'')(\d(e))^{-1} \t' t(e) \,\, \big((\w^{(s',s'')} \circ \d)(e) \big )  \, \big( s''(\d(e))^{-1} \t' t'(e)\big) \, t''(e) \\
&=(\w^{(s,s' \tn s'')} \circ \d)(e) \,\, \w^{(s',s'')}(\d(e)) \,\, (s'(\d(e) s''(\d(e))^{-1} \t' t(e) \,\,  \big( s''(\d(e))^{-1} \t' t'(e)\big) \, t''(e),
\end{align*}
$$\big((t {\otimes} t'){\otimes} t''\big)(e)=(\w^{(s\tn s', s'')} \circ \d)(e)\,\, s''(\d(e))^{-1} \t' \Big(\big((\w^{(s,s')} \circ \d)(e) \big )  \, \big( s'(\d(e))^{-1} \t' t(e)\big) \, t'(e)\Big)  \, t''(e).$$ 
To prove associativity, we therefore need to prove that for each $e \in E'$:
$$\w^{(s,s' \tn s'')}( \d(e)) \,\, \w^{(s',s'')}(\d(e))  =\w^{(s\tn s', s'')} (\d(e))\,\, s''(\d(e))^{-1} \t' \w^{(s,s')}(\d(e))  .$$
We will prove that for each $g\in F$ we have:
\begin{equation}\label{provenow}
\w^{(s,s' \tn s'')}( g)  =\w^{(s\tn s', s'')} (g)\,\,\,\big( s''(g)^{-1} \t' \w^{(s,s')}(g)\big)\,\,\, \w^{(s',s'')}(g) ^{-1}  .
\end{equation}
(It is a nice exercise to prove that this is coherent with \eqref{q}.)
To prove \eqref{provenow}, consider the unique group map: 
\begin{equation*}
\begin{split} W \colon F &\to \big ( (G \ltimes_\t E)\ltimes_\bullet (\{1\} \ltimes_* (E \ltimes_{\t'} L) )\big ) \ltimes_{\triangle} \big (\{1\} \ltimes_* (\{1\}\ltimes_{\t'} \{1\})\ltimes_*((\{1\} \ltimes_* (E \ltimes_{\t'} L)\ltimes_{\bullet'}   (\{1\} \ltimes_{\ad} L)\big)\\
& \quad =\Gr_0(\Tet(\A)),
\end{split}
\end{equation*}
which on the (chosen) free basis $B$ of $F$ is $b\mapsto (\phi(b),s(b),1,s'(b),1,1,1,1,1,s''(b),1,1,1)$, thus:
\begin{equation}\label{defw} \xymatrix{\\W(b)=\\}\xymatrix{\teu{\phi(b)}{s(b)}{1}{s'(b)}{1}{s''(b)}{1}{1}{1}} \end{equation}
By using the morphisms in \eqref{morA}, \eqref{morB} and \eqref{morC} (corresponding to the shown faces of the tetrahedron), and also lemma \ref{good}, we can conclude   that, {for each $g \in F$} (by looking at the value of the compositions of these {morphisms} with $W$, in the chosen free basis of $F$):
\begin{multline*}
\xymatrix{\\W(g)=\\}\xymatrix{\ted{\phi(g)}{s(g)}{1}{s'(g)\de(\w^{(s,s')}(g))^{-1}}{\w^{(s,s')}(g)}{s''(g)\de(\w^{(s\tn s', s'')}(g)^{-1}}{\w^{(s\tn s', s'')}(g)\, \big(\w^{(s',s'')}(g)\big)^{-1}}{1}{\w^{(s',s'')}(g)}}
\end{multline*}
By composing $W\colon F \to \Gr_0(\Tet(\A))$ with $d_3\colon \Gr_0(\Tet(\A)) \to \Gr_0(\T(\A))$, see \eqref{morD}  and \eqref{mordd},  yields a group map $Z\colon F \to \Gr_0(\T(\A))$,
  whose  value (for each {$g \in F$)} is: 
$$\xymatrix{\\Z(g)=\\}\xymatrix{\trixumodd{\phi(g)}
 {s(g)}  {1}
{ s'(g)\,\de(\w^{(s,s')}(g))^{-1} \,s''(b)\,\de( \w^{(s\tn s', s'')} (g)^{-1})}   
{ \w^{(s\tn s', s'')}(g) \,\, (s''(g)^{-1} \t' \w^{(s,s')}(g)) \,\, \w^{(s',s'')}(g)^{-1}}   }  $$
On  free generators $b \in B \subset F$, the map $Z$ has the form:
$$Z(b)=(\phi(b),s(b),1,s'(b)s''(b),1)= \hskip-1.3cm\xymatrix{\trixumodddd{\phi(b)}{s(b)}{1}{s'(b)s''(b)}{1}} $$ 
thus (by using lemma \ref{good} again) {for each $g \in F$}:
$$\xymatrix{\\Z(g)=} \xymatrix{\trixumoddd{\big(\phi(g)}{s(g)}{1}{(s'{\otimes} s'')(g) \de(\w^{(s,s'\tn s''))}(g))^{-1}}{\w^{(s,s'\tn s'')}(g)}} $$
from which we have 
$\w^{(s,s' \tn s'')}( g)  =\w^{(s\tn s', s'')} (g)\,\, s''(g)^{-1} \t' \w^{(s,s')}(g)\,\, \w^{(s',s'')}(g) ^{-1}  ,$  {for each $g \in F$.}
\end{Proof}

\subsubsection{Existence of units}\label{eunits}
Consider a homotopy between 2-crossed module maps:
$$f=({\mu},\psi,\phi) \ra{(f,s,t)} f'=({\mu}',\psi',\phi').$$
Let $(s_0^{f},t_0^{f})$ be the trivial quadratic $f$-derivation, such that $s_0^f(g)=1$ for each $g \in F$ and $t_0^f(e)=1$, if  $e \in E'$.  
\begin{Lemma}\label{uni}
 We have that $$ \w^{(s_0^f,s)}(g)=1 \an \w^{(s,s_0^{f'})}(g)={1, \fo g \in F}.$$
\end{Lemma}
\begin{Proof}
Consider the unique group map 
$$X^{(s_0^f,s)} \colon  F\to (G \ltimes_\t E)\ltimes_\bullet (\{1\} \ltimes_* (E \ltimes_{\t'} L) )=\Gr_0(\T(\A)),$$ of lemma \ref{good}, with, for each free generator $b \in B\subset F$:
$$X^{(s_0^f,s)}(b)=(\f(b),s_0^f(b),1,s(b),1). $$
Thus for each {$g\in F$}:
$$X^{(s_0^f,s)}(g)={\big (\f(g),s_0^f(g),1,s(g)\de(\w^{(s_0^f,s)}(g))^{-1}, \w^{(s_0^f,s)}(g) \big).}$$
{Since $g \in F\longmapsto \big (\f(g),1,1,s(g), 1 \big)\in \Gr_0(\T(\A))$ is also a group morphism, by \eqref{gact} and remark \ref{derprop}, which  extends the value of $X'$ in $B$, it follows in particular that  $\w^{(s_0^f,s)}(g)=1$ for each {$g \in F$.}} 

 {Similarly, consider the map $X\colon  F \to (G \ltimes_\t E)\ltimes_\bullet (\{1\} \ltimes_* (E \ltimes_{\t'} L) )$ of lemma \ref{good} with $s'=s_0^{f'}$. On free generators $X(b)=(\phi(b),s(b),1,1,1)$. Since we have a group morphism with $g \in F \longmapsto (\phi(g),s(g),1,1,1)\in \Gr_0(\T(\A))$, there follows that $X(g)=(\phi(g),s(g),1,1,1)$. In particular $\w^{(s,s_0^{f'})}(g)=1$, for each $g \in F$.}
\end{Proof}
By \eqref{ttn} it  follows that $(s,t) {\otimes} (s_0^{f'},t_0^{f'})=(s,t)$ and that $(s_0^{f},t_0^{f}) {\otimes} (s,t)=(s,t). $
\begin{remark}\label{unip}
 By the proof of lemma \ref{uni}, we in  particular, for a $\phi$-derivation $s\colon F \to E$, have maps:
$$g \in F \mapsto \xymatrix{& & & \phi(g) \, \d (s(g)) \\
&\phi(g)  \ar @/^1pc/ @{{}{ }{}} [rrrr]|{\boxed{1}} \ar[rru]^{s(g)}\ar[rrrr]_{1} & &&& \phi(g) \ar[ull]_{s(g)} 
} \in \Gr_0(\T(\A)) $$
and also:
$$g \in F \mapsto \xymatrix{& & & \phi(g) \, \d (s(g)) \\
&\phi(g)  \ar @/^1pc/ @{{}{ }{}} [rrrr]|{\boxed{1}} \ar[rru]^{s(g)}\ar[rrrr]_{s(g)} & &&& \phi(g) s(g) \ar[ull]_{1} 
} \in \Gr_0(\T(\A))   $$
\end{remark}

\subsubsection{Inverting homotopies in the free up to order one case}\label{ih}
We freely use the notation of subsections \ref{triangle} and \ref{tetrahedron}. 
Let $f,f'\colon \A'=\left (L' \ra{\de} E'\ra{\d} F,  \t,\{,\}\right) \to \A=\left (L \ra{\de} E \ra{\d} G,  \t,\{,\}\right) $ be 2-crossed module maps.
Consider a 2-crossed module homotopy:
$$f=({\mu},\psi,\phi) \ra{(f,s,t)} f'=({\mu}',\psi',\phi').$$
Let us define its inverse:
$$f'=({\mu}',\psi',\phi') \ra{(f',\bar{s},\bar{t})} f=({\mu},\psi,\phi).$$

The derivation $\bar{s}\colon F \to E$ is the unique $\phi'$-derivation which on the chosen basis $B$ of $F$ takes the form: 
$$\bar{s}(b)=s(b)^{-1}. $$
We clearly have $(s {\otimes} \bar{s})(g)=1=(\bar{s} \otimes s)(g)$ for each $g$ in $F$, for this is true in a free basis of $F$, remark \ref{derprop}.
Looking at equation \eqref{q}, for each {$g \in F$} we thus have:
\begin{equation}\label{invq}
 \bar{s}(g)=\big(s(g)\big)^{-1} \de (\w^{(s,\overline{s})})=\de (\w^{(\overline{s},s)})\,\,\big(s(g)\big)^{-1}.
\end{equation}

\begin{Lemma}\label{winv} For each {$g \in F$} we have
 $\w^{(\overline{s},s)}(g)=s(g)^{-1} \t' \w^{(s,\overline{s})}(g)$.
\end{Lemma}
\begin{Proof} 
Consider the map $W \colon F  \to \Gr_0(\Tet(\A))$, 
 in equation \eqref{defw}, for $s'=\bar{s}$ and $s''=s$. Thus, if $b \in B$:
$$\xymatrix{\\W(b)=\\}\xymatrix{\ten{\phi(b)}{s(b)}{1}{s(b)^{-1}}{1}{s(b)}{1}{ 1}{ 1}} $$
 Then, since $(s {\otimes} \overline{s})(g)=1$, for each $g \in F$, and lemmas \ref{uni} and \ref{good}, together with remark \ref{unip} we have:
$$\xymatrix{\\W(g)=\\}\xymatrix{\tep{\phi(g)}{s(g)}{1} {s(g)^{-1}}{\w^{(s,\overline{s})}(g)}{ s(g)}{(\w^{(\overline{s},s)}(g))^{-1}}{ 1}{ \w^{(\overline{s},s)}(g) }}	. $$
Composing $W\colon F \to \Gr_0(\Tet(\A))$ with $d_3\colon \Gr_0(\Tet(\A)) \to \Gr_0(\T(\A))$, see \eqref{morD}  and \eqref{mordd},  gives a  map:
$$\xymatrix{\\ g \in F\mapsto\\} \xymatrix{\trixumodn {\phi(g)}{ {s}(g)} {1}{1} {(s(g)^{-1} \t'\w^{(s,\overline{s})}(g))\,\, \w^{(\overline{s},s)}(g)^{-1}} }\xymatrix{\\ \in \Gr_0(\Tet(\A))} $$
which on the chosen basis $B$ of $F$ has the form: 
$$\xymatrix{\\ b {\longmapsto} \\}\xymatrix{\trixumodp{\phi(b)}{ {s}(b)}{ 1}{1}{ 1} }$$
By lemmas \ref{uni} and \ref{good},  and remark \ref{unip}, it follows that $s(g) \t' \w^{(s,\overline{s})}(g)\,\, \big(\w^{(\overline{s},s)}(g)\big)^{-1}=1$, for each $g \in F$.
\end{Proof}

{We now define $\overline{t}\colon E' \to L.$ For an $e \in E'$, put:}
\begin{equation}\label{deftb}
\bar{t}(e)=\big(\w^{(s,\overline{s})}(\d(e))\big)^{-1} \,\, (s\d(e)) \t'  t(e)^{-1}. 
\end{equation}
\begin{Lemma} The pair
 $(\bar{s},\bar{t})$ is an $f'$-quadratic derivation. 
\end{Lemma} We will use lemmas \ref{cart} and \ref{techlemma}, similarly to the construction of the concatenation of homotopies.\\
\noindent\begin{Proof}
Consider the group map (lemma \ref{good}) $M\colon F \to (G \ltimes_\t E)\ltimes_\bullet (\{1\} \ltimes_* (E \ltimes_{\t'} L) )=\Gr_0(\T(\A))$ with  $$g \stackrel{M}{\longmapsto} \big(\phi(g),s(g),1,\overline{s}(g) \de(\w^{(s,\overline{s})}(g))^{-1}, \w^{(s,\overline{s})}(g)\big)=\big(\phi(g),s(g),1,({s}(g))^{-1}, \w^{(s,\overline{s})}(g)\big	),$$
thus:
$$\xymatrix{\\ g \in F \stackrel{M}{\longmapsto}\\} \hskip-.4cm \xymatrix{\trixu {\phi(g)} {s(g)}{1}{ ({s}(g))^{-1}} { \w^{(s,\overline{s})}(g)}} \xymatrix{\\ =\\} \hskip-.9cm \xymatrix{\trixun {\phi(g)} {s(g)}{1}{ ({s}(g))^{-1}} { \w^{(s,\overline{s})}(g))}} $$
Consider the set map: 
$N \colon E' \to (E \ltimes_* (E \ltimes_{\t'} L))\ltimes_*((\{1\} \ltimes_* (E \ltimes_{\t'} L)\ltimes_{\bullet'}   (\{1\} \ltimes_{\ad} L))=\Gr_1(\T(\A))$ with:
$$e \in E' \stackrel{N}{\longmapsto} \Big (\psi(e),s (\d(e)), t(e), 1,(s(\d(e)))^{-1},\w^{(s,\overline{s})}(\d(e)), 1,\big(\w^{(s,\overline{s})}(\d(e))\big)^{-1} \,\, (s(\d(e)) \t'  t(e)^{-1}\Big), $$
thus:
$$\xymatrix{\\ e \in E'\stackrel{N}{\longmapsto} \\} \xymatrix{ & & & \psi(e)\\ \\
&\psi(e)  \ar[rruu]^{(1,1)}\ar[rrrr]_{\big ( s(\d(e),t(e)\big)} & &&& \psi(e)\,s(\d(e))\,\de(t(e)) \ar[uull]_{\quad \quad \quad \big( \overline{s}(e), (\w^{(s,\overline{s})}(\d(e)))^{-1} \,\, (s(\d(e)) \t'  t(e)^{-1}  \big)}}
$$

Let us see that $N$ is a group morphism and also that $(N,M)$ is a pre-crossed module map from  $(\d\colon E' \to F,\t)$ to the underlying pre-crossed module of  $\T(\A)$.
By  lemmas \ref{cart} and \ref{techlemma} it suffices to check that the diagram below in the category of sets commutes, which is straightforward (we note \eqref{invq}):
 $$\xymatrixcolsep{4pc}
\xymatrix{ & & E'\ar@{-->}[r]^{N} \ar[dll]_\d \ar[dr]|<<<<<<<<<<<<\hole_{\Gr_1(H',H)} &\Gr_1(\T(\A))\ar[d]^{\Gr_1(d_{--})}  \ar[dll]_{\beta'}   \\ F\ar[r]_{M} \ar[rd]_{\Gr_0(H',H)} &\Gr_0(\T(\A))\ar[d]^{\Gr_0(d_{--})}  & &\Gr_1\left(\P({{\A}}) {_{\Pr_0^{{\A}}}\times_{\Pr_0^{{\A}}}} \P({{\A}})\right) \ar[dll]_{\beta'}\\ &  \Gr_0\left(\P({{\A}}) {_{\Pr_0^{{\A}}}\times_{\Pr_0^{{\A}}}} \P({{\A}})\right)}$$ 
Here $H\colon {{\A}}' \to \P({{\A}})$ is associated to the quadratic $f$-derivation $(s,t)$ and $H'$ is associated with the trivial quadratic $f$-derivation $(s_0^f,t_0^f)$; see lemma \ref{esc} and \ref{eunits}. We therefore have an associated map to the 2-crossed module pull-back: $(H',H) \colon \A \to \P({{\A}}) {_{\Pr_0^{{\A}}}\times_{\Pr_0^{{\A}}}} \P({{\A}})$.

Thus $(N,M)$ is a pre-crossed module map. Composing $(N,M)$ with the pre-crossed module map $d_0\colon \T(\A) \to \P(\A)$ (see \eqref{def2v}), yields a pre-crossed module map from $(\d\colon E' \to F,\t) $ to the underlying pre-crossed module of $\P(\A)$, which has the form:
\begin{align*}
{e \in E'}&\longmapsto  \Big(\psi(e) \, s (\d(e)) \,\de (t(e)), (s(\d(e)))^{-1} \de(\w^{(s,\overline{s})}(\d(e))), \w^{(s,\overline{s})}(\d(e))\big)^{-1} \,\, (s\d(e)) \t'  t(e)^{-1}\Big)\\ &={\big(\psi'(e),\bar{s}(\d(e)), \overline{t}(e)\big) \in \Gr_1(\P(\A))} 
\end{align*}
and
$${g \in F \longmapsto \big(\phi(g)\,\d(s(g)),\overline{s}(g)\big)=\big ( \phi'(g),\overline{s}(g)\big) \in \Gr_0(\P(A)),} $$
hence {(by lemma \ref{esc})} it follows that  $(\bar{s},\bar{t})$ an $f'$-quadratic derivation. 
\end{Proof}

Now note that obviously $s{\otimes} \overline{s}=s_0^f$ and $\overline{s} {\otimes} s=s_0^{f'}$, since the same is true in a free basis of $F$. On the other hand (by  \eqref{invq}), {if $e \in E'$}:
\begin{align*}
 (t{\otimes} \overline{t})(e)&=\w^{(s,\overline{s})}(\d(e))\,\, \big(\overline{s}(\d(e))^{-1} \t' t(e) \big) \,\, \big(\w^{(s,\overline{s})}(\d(e))\big)^{-1} \,\, s(\d(e)) \t'  t(e)^{-1}\big)\\\ &=\w^{(s,\overline{s})}(\d(e))\,\,\big( (\de(\w^{(s,\overline{s})}(\d(e)))^{-1} \, s(\d(e))\big) \t' t(e) \,\, \big(\w^{(s,\overline{s})}(\d(e))\big)^{-1} \,\, s(\d(e)) \t'  t(e)^{-1}=1,
\end{align*}
where we used the crossed module rule $\de(k) \t' l=klk^{-1}$. Also {(we use \eqref{invq} again)}:
\begin{align*}
(\overline{t} {\otimes} t)(e)&= \w^{(\overline{s},s)}(\d(e)) \,\, s^{-1}(\d(e)) \t' \overline{t}(\d(e))\,\, t(e) \\
&=  s(\d(e))^{-1} \t' \w^{(s,\overline{s})}(\d(e))\,\, (s(\d(e)))^{-1} \t' \Big(\big(\w^{(s,\overline{s})}(\d(e))\big)^{-1} \,\, (s(\d(e)) \t'  t(e)^{-1}\Big) \,\, t(e) =1.
\end{align*}
Thus we proved that $(\overline{s},\overline{t})$ is an inverse of $(s,t)$.
We have  finished proving the  main result of this subsection:
\begin{Theorem}\label{jh}
Let  $\A'=\left (L' \ra{\de} E'\ra{\d} F,  \t,\{,\}\right)$ and $\A=\left (L \ra{\de} E\ra{\d} G ,  \t,\{,\}\right)$ be 2-crossed modules, where $F$ is a free group, with a chosen basis $B$. We can define a groupoid $[\A',\A]_1^B$ of 2-crossed module maps $\A' \to \A$, and their homotopies.
\end{Theorem}
In the next subsection we will see that this construction can be expanded to be a 2-groupoid $\HOM_B(\A',\A)$, by considering 2-fold homotopies between 2-crossed module homotopies.
\begin{Corollary}
 Let $\A$ and $\A'$ be 2-crossed modules. If $\A'$ is free up to order one then homotopy between 2-crossed module maps $\A' \to \A$ yields an equivalence relation.
\end{Corollary}

\subsection{A 2-groupoid of 2-crossed module maps, 1- and 2-fold homotopies (in the free up to order one case)}
For the definition of a 2-groupoid see \cite{HKK}.
\subsubsection{A groupoid of 2-crossed module  homotopies and their 2-fold homotopies }

Let $\A'=\left (L' \ra{\de} E'\ra{\d} G',  \t,\{,\}\right)$ and $\A=\left (L \ra{\de} E\ra{\d} G ,  \t,\{,\}\right)$ be 2-crossed modules. (In this subsubsection, only, we will not need to suppose that $G'$ is a free group).  Consider two 2-crossed module maps $f,f'\colon \A' \to \A$. Let us define a groupoid  $[f,f']$, with objects the homotopies $(f,s,t)$ connecting $f$ and $f'$ (so the pair $(s,t)$ is a quadratic $f$-derivation), the 2-morphisms being constructed from  2-fold homotopies (quadratic 2-derivations) $k\colon G' \to L$.  For nomenclature and notation we refer to subsection \ref{2folddef}.

The set of object of $[f,f']$ is the set of triples $(f,s,t)$ where $(s,t)$ is  a quadratic $f$-derivation  connecting $f$ and $f'$. The set of 1-morphisms $(f,s,t) \to (f,s',t')$ is made out of quadruples $(f,s,t,k)$, where $k\colon G' \to L$  is a quadratic  $(f,s,t)$ 2-derivation,  such that    $(f,s,t) \ra{(f,s,t,k)} (f,s',t') .$ If we have a chain of arrows:  $$(f,s,t) \ra{(f,s,t,k)} (f,s',t') \ra{(f,s',t',k')} (f,s'',t''),$$
then their {concatenation} is given by the map $k\diamond k' \colon G' \to L$, such that,for each $g \in G'$:
 $$(k\diamond k')(g)=k(g)k'(g).$$
\begin{Lemma}\label{qaz}
The map $(k\diamond k')\colon G' \to L$ is a quadratic  $(f,s,t)$ 2-derivation.
\end{Lemma}
 \begin{Proof}
By equations  \eqref{2der} and \ref{tar}, and since $(\delta\colon L \to E,\t')$ is a crossed module,  we have, for each $g,h \in G'$:
  \begin{align*}
  (k\diamond k') (gh)&\doteq k(gh) k'(gh)\\
&=\Big (  s(h)^{-1} \t'\big (\phi(h)^{-1} \t k(g)\big) \Big)\,\, k(h)\,\, \Big (  s'(h)^{-1} \t'\big (\phi(h)^{-1} \t k'(g)\big) \Big)\,\, k'(h)\\
&=\Big (  s(h)^{-1} \t'\big (\phi(h)^{-1} \t k(g)\big) \Big)\,\, k(h)\,\, \Big ( ( \de(k(h))^{-1}  s(h)^{-1}) \t'\big (\phi(h)^{-1} \t k'(g)\big) \Big)\,\, k'(h)\\
&=\Big (  s(h)^{-1} \t'\big (\phi(h)^{-1} \t k(g)\big) \Big)\,\, \Big (  s(h)^{-1} \t'\big (\phi(h)^{-1} \t k'(g)\big) \Big)\,\,k(h) \,\, k'(h).\\
&=\Big (  s(h)^{-1} \t'\big (\phi(h)^{-1} \t (k\diamond k')(g)\big) \Big)\,\, (k\diamond k')(h).
  \end{align*}
 \end{Proof}

\noindent Note that by equation \eqref{tar}, it follows that $$(f,s,t) \ra{(f,s,t,k\diamond k')} (f,s'',t'').$$
This {concatenation} of  quadratic  $(f,s,t)$ 2-derivations is clearly associative and it has units; the   quadratic $(f,s,t)$-derivation such that $k(g)=1_L,\forall g \in G'$. The fact we have a groupoid  $[f,f']$ follows from:
\begin{Lemma}
 If $(f,s,t) \ra{(f,s,t,k)} (f,s',t'),$ then the map $\overline{k}\colon G' \to L$ such that 
$\overline{k}(g)=k(g)^{-1}$ for each $g \in G'$ is a quadratic  $(f,s',t')$ 2-derivation.
\end{Lemma}
\begin{Proof}By equation  \eqref{2der}, and since $(\delta\colon L \to E,\t')$ is a crossed module,  we have, for each $g,h \in G'$:
\begin{align*}\overline{k}(gh)\doteq k(gh)^{-1} &= \Big (\big (  s(h)^{-1} \t'\big (\phi(h)^{-1} \t k(g)\big) \big) \,\,k(h)\Big)^{-1}\\
 &= \Big (k(h)\,\, \big (  (\de(k(h))^{-1} s(h)^{-1}) \t'\big (\phi(h)^{-1} \t k(g)\big) \big)\Big)^{-1}\\
&=\big (  s'(h)^{-1} \t'\big (\phi(h)^{-1} \t \overline{k}(g)\big) \big)\,\, \overline{k}(h).
\end{align*}
\end{Proof}

Given two 2-crossed modules $\A'$ and  $\A$ we have a groupoid $[\A',\A]_2$, whose objects are arbitrary 2-crossed module homotopies $f\ra{(f,s,t)} f'$, where $f,f'\colon \A \to \A'$, the morphisms being the 2-crossed module 2-fold homotopies. Suppose that $\A'$ is free up to order one, with a chosen basis. To define a 2-groupoid $\HOM_B(\A',\A)$, we now need compatible left and right  actions of the groupoid $[\A',\A]_1^B$ on $[\A',\A]_2$ {(see theorem \ref{jh});} in other words we need whiskering operators. 

\subsubsection{Right whiskering 2-fold homotopies by 1-fold homotopies (in the free up to order one case).}\label{rwisk}

Let $\A'=\left (L' \ra{\de} E'\ra{\d} F ,  \t,\{,\}\right)$ and $\A=\left (L \ra{\de} E\ra{\d} G ,  \t,\{,\}\right)$ be 2-crossed modules. We now go back to assuming  $F$ to be a free group over the chosen basis $B$.
Let $f,f'\colon \A' \to \A$ be 2-crossed module maps. Suppose that  we have two homotopies $(f,s,t)$ and $(f,s',t')$ connecting $f$ and $f'$. 
Suppose  that we have a quadratic  $(f,s,t)$ 2-derivation  $k\colon F \to L$ connecting $(s,t)$ and $(s',t')$, thus $(f,s,t) \ra{(f,s,t,k)} (f,s',t').$ Since the latter will be a 2-morphism in $\HOM_B(\A',\A)$,  we now represent it    as:
 $$\xymatrix{ &f\ar@/^1pc/[rr]^{(f,s',t')}\ar@/_1pc/[rr]_{(f,s,t)} & \Uparrow(f,s,t,k)  &f'}.$$
Let $f''\colon \A' \to \A$ be another 2-crossed module map. Suppose we also have a homotopy: $$f'=({\mu}',\psi',\phi') \ra{(f',s'',t'')} f''=({\mu}'',\psi'',\phi''),$$
so what we have diagrammatically is:
 $$\xymatrix{ &f\ar@/^1pc/[rr]^{(f,s',t')}\ar@/_1pc/[rr]_{(f,s,t)} & \Uparrow(f,s,t,k)  &f' \ar[rr]^	{(f',s'',t'')} & &f''}.$$
Let us define the whiskering:
$$(f,s,t,k){\otimes} (f',s'',t'') =(f,s{\otimes} s'',t{\otimes} t'',k {\otimes} s''),$$
such that {$k {\otimes} s''$ connects $(s{\otimes} s'' ,t{\otimes} t'')$ and $(s'{\otimes} s'' ,t'{\otimes} t'')$; diagrammatically:}
 $$\xymatrix{ &f\ar@/^1pc/[rr]^{(f,s'{\otimes} s'',t'{\otimes} t'')}\ar@/_1pc/[rr]_{(f,s{\otimes} s'' ,t{\otimes} t'')} & \Uparrow (f,s,t,k){\otimes} (f',s'',t'')  &f''}.$$
By definition (see corolary \ref{extt}), $k {\otimes} s''$ is the unique quadratic $(f,s{\otimes} s'',t'{\otimes} t'')$ 2-derivation $F \to L$, which on the chosen basis $B$ of $F$ has the form:
$$ (k {\otimes} s'')(b)=s''(b)^{-1} \t' k(b). $$
Then we have for each $g \in F$:
\begin{equation}\label{clear}
(s {\otimes} s'')(g)\, \de((k {\otimes} s'')(g))=(s'{\otimes} s'')(g);
\end{equation}
c.f equation \eqref{tar}. This is because (remark \ref{derprop}), on the free generators $b \in B \subset  F$, we have:
$$ (s {\otimes} s'')(b)\, \de((k {\otimes} s'')(b))=s(b)\, s''(b)\,\de(s''(b)^{-1} \t' k(b))=s(b)\,\de(k(b)) \, s''(b)=s'(b)\, s''(b)=(s' {\otimes} s'')(b).  $$
\begin{Lemma}
The following holds for each $e \in E'$ (c.f. equation \eqref{tar}):
\begin{equation}\label{teq}
(k {\otimes} s'')(\d(e))^{-1}\,\, (t {\otimes} t'')(e)=(t'{\otimes} t'')(e).
\end{equation}
\end{Lemma}
\begin{Proof} We freely use the notation introduced in \ref{triangle}, \ref{disk} and \ref{tetrahedron}.
We have, {for each $e \in E'$:}
\begin{align*}
(k {\otimes} s'')(\d(e))^{-1}\,\, (t {\otimes} t'')(e)&= (k {\otimes} s'')(\d(e))^{-1}\,\, \w^{(s,s'')}(\d(e))\,\,s''(\d(e))^{-1} \t' t(e)\,\, t''(e), 
\end{align*}
whereas (putting $s\, \de(k)\colon F \to L$ as being the derivation $g \in F \mapsto s(g)\, \de(k(g) \in L$):
\begin{align*}
(t' {\otimes} t'')(e)&=\w^{(s',s'')}(\d(e))\,\,  s''(\d(e))^{-1} \t' t'(e) \,\, t''(e)\\ 
&=\w^{(s\, \de(k),s'')}(\d(e))\,\, s''(\d(e))^{-1} \t' ((k\circ \d(e))^{-1} t(e)) \,\, t''(e).
\end{align*}
Therefore \eqref{teq} is equivalent to:
$$(k {\otimes} s'')(\d(e))^{-1}\,\, \w^{(s,s'')}(\d(e))=\w^{(s\de(k),s'')}(\d(e))\,\, s''(\d(e))^{-1} \t' ((k\circ \d(e))^{-1}, {\fo e \in E'}.$$

Let us then prove that for any $g\in F$ we have:
\begin{equation}\label{geq}
 (k {\otimes} s'')(g)^{-1}\,\, \w^{(s,s'')}(g)=\w^{(s',s'')}(g)\,\, s''(g)^{-1} \t' k(g)^{-1}.
\end{equation}
The technique of proof is entirely analogous to the proof that the concatenation of homotopies is associative. Consider the unique map $K \colon F \to \Gr_0(\Tet(\A))$,
which on the free basis $B$ of $F$ is:
$$ K(b)= (\phi(b),s(b),1,s''(b),1,1,1,1,1,\de(s''(b)^{-1} \t' k(b)),s''(b)^{-1} \t' k(b)^{-1} ,1,1),$$
that is, for each $b \in B$:
\begin{equation}\label{defk}
 \xymatrix{\\K(b)= \\} \xymatrixcolsep{4pc} \hskip-1cm \xymatrix{
\teo
{\phi(b)}
{s(b)}
{1}
{{s''(b)}}
{1}
{  {\de(s''(b)^{-1} \t' k(b))}  }
  {({s''}(b)^{-1} \t' k(b)^{-1}) }   
{1}
{1}
}
\end{equation}
By using the morphisms in \eqref{morA}, \eqref{morB} and \eqref{morC},  \ref{comphom}, lemma \ref{good} and remarks \ref{tte} and \ref{unip}, we can conclude:
$$\xymatrix{\\ K(g)=\\}\xymatrixcolsep{4pc} \hskip-1cm  \xymatrix{\tei{\phi(g)}{s(g)}{1}{s''(g) \de(\w^{(s,s'')})(g)^{-1}}{\w^{(s,s'')}(g)}{\de(k{\otimes} s'')(g)}{(k{\otimes} s'')(g)^{-1}}{1}{1}}$$
for each $g\in G$.
Consider the unique map $K' \colon F \to \Gr_0(\Tet(\A))$,
which on the free basis $B$ of $F$ is:
$$ K'(b)= (\phi(b),s(b),1,\de(k(b)),k(b)^{-1},1,1,1,1,s''(b),1,1,1),$$
that is, if $b \in B$:
\begin{equation}\label{defkp} \xymatrix{\\
K'(b)= \\}\xymatrix{\ter{\phi(b)}{s(b)}{1}{\de(k(b))}{k(b)^{-1}}{s''(b)}{1}{1}{1}}
\end{equation}
Clearly, by lemma \ref{good} and remarks \ref{unip}, \ref{tte} and \ref{derprop}, we have if $g \in F$:
$$\xymatrix{\\K'(g)=\\}\xymatrix{\tee{\phi(g)}{s(g)}{1}{\de(k(g))}{k(g)^{-1}}{ s''(g)\,\,\de(\w^{(s',s'')}(g))^{-1}}{ \w^{(s',s'')}(g)}{1}{1}}$$

 By composing $K$ and $K'$ with the group morphism $d_3\colon \Tet(\A) \to \T(\A) $, of \eqref{morD} and \eqref{mordd},  yields two group morphisms $F \to  \Gr_0(\T(\A))$, namely:
$$\xymatrix{ g  \in F \longmapsto \\}\xymatrix{\trixumods{\phi(g)}{s(g)}{1}{s''(g) \de(\w^{(s,s'')})(g)^{-1}\,\,\de(k{\otimes} s'')(g)}{(k{\otimes} s'')(g)^{-1}\, \w^{(s,s'')}(g))}}$$
and 
$$ g \in F \longmapsto \xymatrix{\trixumoda{\phi(g)}{s(g)}{1}{ \de(k(g))\, s''(g)\,\,\de(\w^{(s',s'')}(g))^{-1}}{ \w^{(s',s'')}(g) \,\,  s''(g)^{-1}\t' k(g)^{-1}}} $$
(It is an instructive exercise to check that the last two equations are coherent with \eqref{clear} and \eqref{q}.)
Since these morphisms agree on the chosen basis $B$ of $F$, they coincide, thus equation \eqref{geq} follows.  \end{Proof}

Therefore, by \eqref{clear} and \eqref{teq}, we have that (by equation \eqref{tar}):
$$(f,s{\otimes} s'', t {\otimes} t'') \ra{(f, s{\otimes} s'', t {\otimes} t'',k {\otimes}s'')}(f,s'{\otimes} s'', t ' {\otimes} t'') .$$

\begin{Lemma}[Functoriality of the right whiskering] {Let $f,f',f''\colon \A' \to \A$ be 2-crossed module maps. Suppose that  we are given homotopies  $(f,s,t)$, $(f,s',t')$ and $(f,s'',t'')$, connecting $f$ and $f'$, as well as 2-fold homotopies $(f,s,t) \ra{(f,s,t,k)} (f,s',t')$ and $(f,s',t') \ra{(f,s',t',k')} (f,s'',t'')$. Suppose we are  also given a homotopy $f' \ra{(f',u,v)} f''$. Diagrammatically we have:}
 $$\xymatrix{ &f\ar[rr]|{(f,s',t')}\ar@/_3pc/[rr]_{(f,s,t)} \ar@/^3pc/[rr]^{(f,s'',t'')} \ar@/_1pc/ @{{}{ }{}}[rr]|{\Uparrow(f,s,t,k)}  \ar@/^1pc/ @{{}{ }{}}[rr]|{\Uparrow(f,s',t',k')}  & &f' \ar[rr]^	{(f',u,v)} & &f''}.$$
Then:
$$\big ( (f,s,t,k)\diamond (f,s',t',k') \big) \otimes (f',u,v) = \big ( (f,s,t,k) \otimes (f',u,v)  \big)  \diamond \big (  (f,s',t',k') \otimes (f',u,v)\big). $$
\end{Lemma}
\begin{Proof}
 In the left-hand-side, since $(f,s,t,k)\diamond (f,s',t',k') =(f,s,t,kk')$, the underlying quadratic $(f,s\otimes u, t \otimes v)$ 2-derivation $(kk')\otimes u$  is the unique quadratic $(f,s\otimes u, t \otimes v)$ 2-derivation $F \to L$  which on the chosen basis $B$ of $F$ is $b \longmapsto u(b)^{-1} \t' \big (k(b) k'(b)\big). $ On the right hand side we have the  quadratic $(f,s\otimes u, t \otimes v)$ 2-derivation $F \to L$, which is the product of $k \otimes u$ and  $k' \otimes u$. On the chosen basis $B$ of $F$ it takes the form:
$$b \longmapsto \big( u(b)^{-1} \t' k(b) \big) \,\, \big(  u(b)^{-1} \t' k'(b)\big)= u(b)^{-1} \t' \big (k(b) k'(b)\big).$$
Therefore: $$(kk') \otimes u=(k\otimes u) \,(k' \otimes u),$$
since the same is true in a free basis of $F$.
\end{Proof}

Analogously:
\begin{Lemma}
 Suppose  we have a 2-fold  $(f,s,t)$ homotopy $(f,s,t) \ra{(f,s,t,k)} (f,s',t')$. Consider also a chain of homotopies:   $f' \ra{(f',u,v)} f'' \ra{(f'',u',v')} f'''$; diagrammatically:
 $$\xymatrix{ &f\ar@/^1pc/[rr]^{(f,s',t')}\ar@/_1pc/[rr]_{(f,s,t)} & \Uparrow(f,s,t,k)  &f'  \ar[rr]^{(f',u,v)} & &f'' \ar[rr]^{(f',u',v')} & & f''' }.$$
Then:
$$(f,s,t,k) \otimes (u \otimes u', v \otimes v')=\big( (f,s,t,k) \otimes (u,v)\big) \otimes (u',v'). $$
\end{Lemma}
\begin{Proof}
 In $B$, the underlying quadratic $(f,s \tn u \tn u',t \tn v \tn v')$ 2-derivation in the left-hand-side is:
$$b\longmapsto (u\otimes u') (b)^{-1} \t'k(b)=(u(b)u'(b))^{-1} \t' k(b),$$
whereas the underlying quadratic $(f,s \tn u \tn u',t \tn v \tn v')$   2-derivation on the right-hand-side restricts to:
$$b\longmapsto  u'(b)^{-1} \t' (u (b)^{-1} \t'k(b)).$$
Now apply corollary \ref{extt}.
\end{Proof}

Therefore 
\begin{Proposition}
Let $\A'=\left (L' \ra{\de} E'\ra{\d} F ,  \t,\{,\}\right)$ and $\A=\left (L \ra{\de} E\ra{\d} G ,  \t,\{,\}\right)$ be 2-crossed modules, where $F$ is a free group over the chosen basis $B$.  Whiskering on the right gives a right action (by groupoid morphisms) of the groupoid $[\A',\A]_1^B$  of  maps $\A \to \A'$ and their homotopies, on the groupoid $[\A',\A]_2$ of 2-crossed module homotopies  and their  2-fold homotopies. 
\end{Proposition}

\subsubsection{Left whiskering 2-fold homotopies by 1-fold homotopies (in the free up to order one case)}\label{next}
 We freely use the notation introduced in \ref{triangle}, \ref{disk} and \ref{tetrahedron}. The discussion is very similar to the one in \ref{rwisk}.
Let $\A'=\left (L' \ra{\de} E'\ra{\d} F ,  \t,\{,\}\right)$ and $\A=\left (L \ra{\de} E\ra{\d} G ,  \t,\{,\}\right)$ be 2-crossed modules, where  $F$ is a free group over the chosen basis $B$.
If we have 2-crossed module maps $f,f'\colon \A' \to \A$, homotopies $(f,s,t)$ and $(f,s',t')$ and a 2-fold homotopy $(f,s,t,k)$, all fitting into the diagram:
 $$\xymatrix{ &f\ar@/^1pc/[rr]^{(f,s',t')}\ar@/_1pc/[rr]_{(f,s,t)} & \Uparrow(f,s,t,k)  &f'},$$
and we also have a homotopy $$f''=({\mu}'',\psi'',\phi'') \ra{(f'',s'',t'')} f=({\mu},\psi,\phi),$$
so what we have is:
 $$\xymatrix{ &f'' \ar[rr]^{(f'',s'',t'')} & &f\ar@/^1pc/[rr]^{(f,s',t')}\ar@/_1pc/[rr]_{(f,s,t)} & \Uparrow(f,s,t,k) & f' },$$
let us define the whiskering:
$$(f',s'',t'') {\otimes} (f,s,t,k)  =(f,s''{\otimes} s,t''{\otimes} t, s''{\otimes} k),$$
such that we have:
 $$\xymatrix{ &f''\ar@/^1pc/[rr]^{(f'',s''{\otimes} s',t''{\otimes} t'')}\ar@/_1pc/[rr]_{(f'',s''{\otimes} s ,t''{\otimes} t)} & \Uparrow  (f'',s'',t'') {\otimes} (f,s,t,k)  &f'}.$$
We put $s'' {\otimes} k$ as being unique quadratic $(f'',s''{\otimes} s ,t''{\otimes} t)$ 2-derivation which on the basis $B$ of $F$ is:
$$(s''{\otimes} k)(b)=k(b) .$$
Then we have for each $g \in F$:
\begin{equation}\label{clear2}
 (s'' {\otimes} s)(g)\, \de((s'' {\otimes} k)(g))=(s'{\otimes} s'')(g);
\end{equation}
c.f. equation \eqref{tar}. This is because on the free generators $b \in F$ we have:
$$ (s'' {\otimes} s)(b)\, \de((s'' {\otimes} k)(b))=s''(b)\, s(b)\,\de( k(b))=s''(b)\,s'(b)=(s'' {\otimes} s')(b).  $$
\begin{Lemma}
We  have, for each $e \in E'$ (c.f. equation \eqref{tar}):
\begin{equation}\label{teq2}
( s'' {\otimes} k)(\d(e))^{-1}\,\, (t'' {\otimes} t)(e)=(t''{\otimes} t')(e).
\end{equation}
\end{Lemma}
\begin{Proof}
Equation \ref{teq2} follows if we prove that for each $g\in F$ we have:
\begin{equation}\label{pretar}
(s'' {\otimes} k)(g)^{-1}\,\,\w^{(s'',s)}(g)=\w^{(s'',s')}(g) \,\,k(g)^{-1}.
\end{equation}
Consider the unique map: $K \colon F \to \Gr_0(\Tet(\A)),$
which on the free basis $B$ of $F$ is:
$$ K(b)= (\phi''(b),s''(b),1,s(b),1,1,1,1,1, \de(k(b)), 1,1, k(b)^{-1} ),$$
or 
$$ K(b)= \xymatrix{\tew{\phi''(b)}{s''(b)}{1}{s(b)}{1}{ \de(k(b))}{1}{1}{ k(b)^{-1} }}$$
This map has the following form for each $g \in F$ (by lemma \ref{good} and remarks \ref{unip}, \ref{tte} and \ref{derprop}):
\begin{equation*}
g \in F \stackrel{K}{\longmapsto}\hskip-3cm\xymatrixcolsep{5pc}\xymatrix{ \teq{\phi''(g)}{s''(g)}{1}{s(g)\de(\w^{(s'',s)})(g)^{-1} }{\w^{(s'',s)}(g)}{\de(s''{\otimes} k(g))}{ (s''{\otimes} k(g))^{-1} k(g)}{1}{k(g)^{-1} }}
\end{equation*}
By composing $K\colon F \to \Gr_0(\Tet(\A)$ with the 2-crossed module morphism $d_3\colon \Tet(\A) \to \T(\A)$ in \eqref{morD} yields a group morphism $F \to  \Gr_0(\T(\A))$  with the form:
$$ g \in F\stackrel{d_3\circ K}{\longmapsto} \xymatrix{\trixumodg{\phi''(g)}{s''(g)}{1}{s(g)\,\,\de(\w^{(s'',s)})(g)^{-1}\,\de(s''{\otimes} k(g))} { (s''{\otimes} k(g))^{-1}\,\, \w^{(s'',s)}(g)\,  \, k(g)}}$$
{By lemma \ref{good}},  we have another group morphism $F \to  \Gr_0(\T(\A))$ with the form:
$$g\longmapsto \xymatrix{\trixumodf{\phi''(g)}{s''(g)}{1}{s'(g)\de(\w^{(s'',s')})(g)^{-1}}{\w^{(s'',s')}(g)}} $$
Since these morphisms agree on a basis of $F$, they coincide, from which \eqref{pretar}, thus \eqref{teq2}, follows. 
\end{Proof}

From  equation \eqref{tar}, by \eqref{clear2} and \eqref{teq2} we therefore have:
$$(f'',s''{\otimes} s, t'' {\otimes} t) \ra{(f'', s''{\otimes} s, t'' {\otimes} t,s''{\otimes} k)}(f,s''{\otimes} s', t '' {\otimes} t') .$$
(We have used the crossed module rule $\de(k) \t' l= klk^{-1}$ for all $k,l\in L$.)
As before we can easily prove that 
\begin{Proposition}
 Whiskering on the left gives a left action (by groupoid morphisms) of the groupoid $[\A',\A]_1^B$ of maps $\A \to \A'$ and their homotopies on the groupoid $[\A',\A]_2$ of 2-crossed module homotopies  and their  2-fold homotopies.
\end{Proposition}
\begin{Proposition}
 Whiskering on the right commutes with wiskering on the left. 
\end{Proposition}

If $\A'$ is free up to order one, with a chosen basis, we have therefore constructed a sesquigroupoid \cite{St} $\HOM_B(\A',\A)$ with objects being the 2-crossed module maps $f\colon \A' \to \A$, and the morphisms and 2-morphisms being homotopies and 2-fold homotopies. To prove that $\HOM_B(\A',\A)$ is a 2-
groupoid we now need to prove that it satisfies the interchange law.   
\subsubsection{The interchange law}

Let $\A'=\left (L' \ra{\de} E'\ra{\d} F ,  \t,\{,\}\right)$ and $\A=\left (L \ra{\de} E\ra{\d} G ,  \t,\{,\}\right)$ be 2-crossed modules, where  $F$ is a free group over the chosen basis $B$. Suppose we have the following diagram of 2-crossed module maps $f,f',f''\colon \A' \to \A'$, homotopies and 2-fold homotopies:
 $$\xymatrix{ &f\ar@/^1pc/[rr]^{(f,s',t')}\ar@/_1pc/[rr]_{(f,s,t)} & \Uparrow(f,s,t,k) & f' \ar@/^1pc/[rr]^{(f',u',v')}\ar@/_1pc/[rr]_{(f',u,v)} & \Uparrow(f,u,v,k') & f'' }.$$
Let us prove the interchange law:
$$\big((f,s,t,k) \otimes (f',u,v)\big) \diamond \big(  (f,s',t')\otimes (f,u,v,k') \big) =
\big ((f,s,t) \otimes  (f,u,v,k')\big) \diamond \big ( (f,s,t,k) \otimes  (f',u',v')\big). $$
In the left hand side we have the quadratic $(f,s\otimes u, t \otimes v)$ 2-derivation, which on the basis $B$ of $F$ is:
$$b \longmapsto (u^{-1} (b) \t' k(b) )\,\, k'(b).$$
In the right hand side we have the $(f,s\otimes u, t \otimes v)$ 2-derivation, which on the chosen basis $B$ of $F$ is: 
$$b \longmapsto k'(b) \,\, u'(b)^{-1} \t' k(b)= k'(b)\,\, (\de(k'(b)^{-1}) u(b)^{-1}) \t' k(b)= (u^{-1} (b) \t' k(b) )\,\, k'(b),$$
where we have used \eqref{tar} and the crossed module rules, recalling that $(\de\colon L \to E, \t')$ is a crossed module.

We therefore proved that:
\begin{Theorem}[Mapping space 2-groupoid]\label{ms2g}
Given two 2-crossed modules $\A=\left (L \ra{\de} E \ra{\d} G,  \t,\{,\}\right)$ and 
$\A'=\left (L' \ra{\de} E' \ra{\d} F ,  \t,\{,\}\right)$,
with $F$ free up to order 1, with a chosen basis $B$ of $F$, there exists a 2-groupoid:
$$\HOM_B(\A',\A) $$
of 2-crossed module maps, 1-fold homotopies between 2-crossed module maps, and 2-fold homotopies between 1-fold homotopies. 
\end{Theorem}
We note that $\HOM_B(\A',\A) $  explicitly depends on the chosen basis $B$ of $F$.
\section{Lax homotopy of 2-crossed modules}

\subsection{Definition of $Q^1(\A)$ and lax homotopy}
Consider a 2-crossed module of groups $\A=\left (L \ra{\de} E \ra{\d} G ,  \t,\{,\}\right).$
 We will consider a very natural partial resolution $Q^1(\A)$ of it, which is free up  to order one, with a chosen basis, together with a surjective projection ${\rm proj}\colon Q^1(\A) \to \A$,  defining isomorphisms at the level of 2-crossed module homotopy groups.
It is proven in \cite{Go} that $Q^1$  (clearly functorial by its construction) is a part of a  comonad. We will then use $Q^1(\A)$ to define lax homotopy of (strict) 2-crossed module maps.

\subsubsection{Construction of $Q^1(\A)$ and abstract definition of $Q^1$-lax homotopy}\label{coq}

Let $G$ be a group. The free group on the underlying set of $G$ is denoted by $\FG(G)$.  The inclusion (set) map $G \to \FG(G)$ is denoted by $g\in G \longmapsto [g]\in \FG(G)$. The projection (group) map sending $[g]\in \FG(G)$ to $g \in G$ is denoted by $p\colon \FG(G) \to G$. Note that we do not take $[1]$, where $1$ is the identity of $G$, to be the identity of $\FG(G)$, the latter being the empty word, denoted by $\e$.

For a 2-crossed module  $\A=\left (L \ra{\de} E \ra{\d} G ,  \t,\{,\}\right)$, we put:
$$Q^1(\A)=\left(L \ra{\de'} E {{}_\d \times p} \FG(G) \ra{\d'} \FG(G),\t,\{,\}\right),$$
where of course:
$$ E {{}_\d \times p} \FG(G)=\{(e,u)\in E \times \FG(G): \d(e)=p(u)\}.$$
Moreover $\d'(e,u)=u$, and, for all $u \in \FG(G)$, we put $u \t (e,u')=(p(u) \t e, uu'u^{-1})$.
Clearly the Peiffer pairing is  $\langle (e,u),(e',u')\rangle=(\langle e,e'\rangle,\e)$.
We also put $\de'(k)=(\de(k),\e)$ and also $u \t k=p(u) \t k$, for all $k \in L$ and $u \in \FG(G)$. It is immediate that with the Peiffer lifting:  $$\{(e,u),(e',u')\}=\{e,e'\},$$
this defines a 2-crossed module of groups (a very similar construction appears in \cite{AAO}). Also, $Q^1(\A)$ is free up to order one, and we choose the free basis $[G]=\{[g], g \in G\}$ of $\FG(G)$.

There is a projection ${\rm proj}=(r,q,p) \colon Q^1(\A) \to \A$, which, rather clearly,  yields isomorphisms at the level of 2-crossed module homotopy groups (the homology groups of the underlying complexes). It has the form:
\begin{equation}\label{defp1}
{\rm proj} \quad  =  \xymatrix{ &L \ar[r]^-{\de'}\ar[d]_r& E {{}_\d \times p} \FG(G) \ar[r]^-{\d'}\ar[d]_q& \FG(G)\ar[d]^p\\
&L \ar[r]_\de & E \ar[r]_\d  & G}
\end{equation}
where $r=\id$ and $q(e,u)=e$, for each $(e,u) \in  E {{}_\d \times p} \FG(G) $.

 Consider 2-crossed modules $\A=\left (L \ra{\de} E \ra{\d} G ,  \t,\{,\}\right)$ and 
$\A'=\left (L' \ra{\de} E' \ra{\d} G' ,  \t,\{,\}\right)$. If we have a 2-crossed module morphism $f\colon \A \to \A'$ then $f\circ {\rm proj}$ is a 2-crossed module morphism $Q^1(\A) \to \A'$. This yields an injective map ${\rm proj}^*\colon \Hom(\A,\A') \to \Hom(Q^1(\A),\A')$, since ${\rm proj}\colon Q^1(\A) \to \A$ is surjective. Here $\Hom(\A,\A')$ denotes the set of 2-crossed module maps $\A\to \A'$.

\begin{Definition}[Lax mapping space]\label{laxmspace}
 A morphism $Q^1(\A)\to \A'$ is said to be  a strict map $\A \to \A'$ if it factors (uniquely) through  ${\rm proj}\colon Q^1(\A) \to \A$.  The lax mapping space 2-groupoid: $$\EuScript{HOM}_{\rm LAX}(\A,\A')$$ is the full sub-2-groupoid of $\HOM_{[G]}(Q^1(\A),\A')$, theorem \ref{ms2g}, with objects the strict maps $f\colon \A \to \A'$, each  uniquely identified with $f\circ {\rm proj}\colon Q^1(\A) \to \A'$ . The objects of $\EuScript{HOM}_{\rm LAX}(\A,\A')$ are therefore in one-to-one correspondence with  2-crossed module maps $\A \to \A'$, and we call the 1- and 2-morphisms of $\EuScript{HOM}_{\rm LAX}(\A,\A')$ ``lax homotopies'' and ``lax 2-fold homotopies''. 
\end{Definition}

\subsubsection{The structure of $Q^1(\A)$}

We now want to completely unpack the definition of $\EuScript{HOM}_{\rm LAX}(\A,\A')$, definition \ref{laxmspace}. To unravel the structure of $Q^1(\A)$, we prove an  auxiliary lemma,  describing the kernel of the projection map ${\rm proj}\colon Q^1(\A) \to \A$ in \eqref{defp1}. To characterize this kernel it  suffices to elucidate the  kernel $\ker(p)$ of the obvious projection $p\colon \FG(G) \to G$.  Then the kernel of the projection map ${\rm proj}\colon Q^1(\A) \to \A$ is the 2-crossed module:
$$\ker({\rm proj})=\{{1} \to \{1_E\} \times \ker(p) \to \ker(p)\}, $$
with action by conjugation and  trivial Peiffer lifting.

For details on the definition of free crossed modules and free pre-crossed modules (possibly with ulterior relations), we refer to \cite{FM1,BHS}. Let $G$ be a group. Given $g,h \in G$ put: 
\begin{equation}
[g,h]= [gh]^{-1}[g][h] \in \ker(p) \subset \FG(G).
\end{equation}
 Note that we always have:
\begin{equation}
[gh,i]\,\,[i]^{-1}\,\,[g,h]\,\,[i]=[g,hi]\,\,[h,i], \textrm{ where } g,h,i \in G.
\end{equation}
Also if $g,h \in G$:
\begin{align}
&[g]\,\,[h]\,\,[g,h]^{-1}=[gh], 
&&[1]=[1,1],
&&[g^{-1}]=[g]^{-1}\,\,[1]\,\,[g,g^{-1}].
 \end{align}
Moreover:
\begin{equation}
 \begin{split}
[ghg^{-1}]&=[g]\,\,[hg^{-1}]\,\,[g,hg^{-1}]^{-1}= [g]\,\,[h]\,\,[g^{-1}]\,\,[h,g^{-1}]^{-1}\,\,[g,hg^{-1}]^{-1}\\
&=[g]\,\,[h]\,\,[g]^{-1}\,\, [1]\,\,[g,g^{-1}]\,\, [h,g^{-1}]^{-1}\,\,[g,hg^{-1}]^{-1}.
\end{split}
\end{equation}

\begin{Lemma}\label{key}
 The inclusion map $\iota\colon \ker(p) \to \FG(G)$, together with the action $\t$ of $\FG(G)$ on $\ker(p)\subset \FG(G)$ by conjugation,  is isomorphic to the crossed module, over $\FG(G)$, formally generated by the elements $(g,h)$, where $g,h \in G$, with:
$$\iota(g,h)= [g,h], \textrm{ for each } g,h \in G ,$$
modulo the relations:
\begin{equation}\label{relations}
(gh,i)\,\, [i]^{-1} \t (g,h)=(g,hi)\,\, (h,i), \textrm{ where } g,h,i \in G .
\end{equation}
In particular,  $\ker(p) \to \FG(G)$ is isomorphic to the  pre-crossed module, formally generated by  the symbols $(g,h)$, for all $g,h \in G$, with $\iota(g,h)=[g,h]$, modulo the relations \eqref{relations}, as well as the following relations (where $g,g',h,h' \in G$ and $k,k' \in \FG(G)$), enforcing the second Peiffer condition in   definition \ref{cm}:
$$\big(\iota (k \t (g,h))\big) \t  \big (k' \t (g',h')\big)=\big(k \t (g,h) \big)\,\, \big (k' \t (g',h')\big)\,\, \big(  k \t (g,h)\big)^{-1}. $$
(It suffices to consider the case $k=\e$.)

Also, we have that $\ker(p)\subset \FG(G)$, as a group,  is generated by all conjugates (under $\FG(G)$) of elements $[g,h]\in \FG(G)$, where $g,h \in G$, and the relations \eqref{relations} are the only relations between these. 
\end{Lemma}
We will give a topological proof of this lemma. Recall \cite{BHS} that, if $(X,Y)$ is a pair of path-connected spaces, then the boundary map $\pi_2(X,Y) \to \pi_1(Y)$, together with the standard action of $\pi_1$ on $\pi_2$, defines a crossed module $\Pi_2(X,Y)$, a result due to Whitehead. If $X$ is obtained from $Y$ by attaching 2-cells,  then $\Pi_2(X,Y)$ is the free crossed module on the attaching maps of the 2-cells of $X$ in $\pi_1(Y)$, a fact usually known as Whitehead theorem \cite{W1,W2,W3}. If $X$ is a CW-complex, then $X^i$ denotes the $i$-skeleton of $X$.

\noindent \begin{Proof}
Let $K$ be the simplicial set which is the nerve of $G$, thus the geometric realisation of $K$ is the usual classifying space of $G$; see for example \cite{BHS}. We thus have a unique 0-simplex, the 1-simplices of $K$ are in one-to-one correspondence with elements of $G$, and we denote these by $[g]$, where $g \in G $.  The 2-simplices of $K$ are in one-to-one correspondence with pairs $(g,h),$ where $g,h \in G$, being:
\begin{align*}
\d_0(g,h)&=[h] &\d_1(g,h)&=[gh] &\d_2(g,h)&=[g].
\end{align*}
The $3$-simplices of $K$ are triples $(g,h,i)$ of elements of $G$, being:
\begin{align*}
\d_0(g,h,i)&=(h,i) &\d_1(g,h,i)&=(gh,i) &\d_2(g,h,i)&=(g,hi) &\d_3(g,h,i)&=(g,h).
\end{align*}
Let us consider the fat geometric realization $X$ of $K$ (forgetting about the degeneracy maps of $K$, thus looking at $K$ merely as being a ${\triangle}$-complex, \cite{H}.) It is well known that the fat and standard geometric realisations of a simplicial set are homotopy equivalent (see for example \cite{BS}), thus $X$ is an aspherical CW-complex with $\pi_1(X)\cong G$.

Consider the exact sequence $\{1\} \to \pi_2(X,X^1) \to \pi_1(X^1) \to \pi_1(X)\cong G$. These are the final groups of the long homotopy exact sequence of the pair $(X,X^1)$, where $X^1$ is the 1-skeleton of $X$. Then $\pi_1(X^1) \to \pi_1(X)$ is exactly the map $p\colon \FG(G) \to G$, so we only need to determine $\pi_2(X,X^1)$. The crossed module  $(\pi_2(X,X^1) \to \pi_1(X^1))=\Pi_2(X,X^1)$ is a quotient of the crossed module $(\pi_2(X^2,X^1) \to \pi_1(X^1))=\Pi_2(X^2,X^1)$, which, by Whitehead theorem, is  the free crossed module on the boundary maps of the 2-cells of $X$ in $\pi_1(X^1)$. Thus $\pi_2(X^2,X^1)$ is the principal  group of the pre-crossed module,  over $\pi_1(X^1)$, generated by all pairs $(g,h)$ where  $g,h \in G$, with $\iota(g,h)=[g,h]$,  modulo the relations (for $g,g',h,h' \in G$ and $k,k' \in \FG(G)$), enforcing the second Peiffer condition in definition \ref{cm}:
$$\big(\iota (k \t (g,h))\big) \t \big (k' \t (g',h')\big)=\big(k \t (g,h) \big)\,\, \big (k' \t (g',h')\big)\,\, \big(  k \t (g,h)\big)^{-1}. $$
To obtain $\pi_2(X,X^1)$ from $\pi_2(X^2,X^1)$, we now need to add one extra relation for each 3-cell of $X$ (see \cite{FM1}), yielding that we should have:
$$(gh,i)\,\, [i]^{-1} \t (g,h)=(g,hi)\,\, (h,i), \textrm{ where } g,h,i \in G ,$$
which arise from all 3-simplices of $K$.
\end{Proof}

\subsubsection{Explicit form of a lax homotopy between two strict 2-crossed module maps}\label{j1}
 We use the notation of subsection \ref{qder} and \ref{coq}. Consider 2-crossed modules $\A=\left (L \ra{\de} E \ra{\d} G ,  \t,\{,\}\right)$ and 
$\A'=\left (L' \ra{\de} E' \ra{\d} G' ,  \t,\{,\}\right)$. 
Recall the construction of $\EuScript{HOM}_{\rm LAX}(\A,\A')$, definition \ref{laxmspace}, from $Q^1(\A)=\left(L \ra{\de'} E {{}_\d \times p} \FG(G) \ra{\d'} \FG(G),\t,\{,\}\right)$.
Let $f_1=({\mu}_1,\psi_1,\phi_1)$ and $f_2=({\mu}_2,\psi_2,\phi_2)$ be 2-crossed module maps $\A \to \A'$. Consider a lax homotopy connecting $f_1$ and $f_2$, namely: $$({\mu}_1',\psi_1',\phi_1')=f_1'\doteq (f_1\circ{\rm proj}) \ra{(f_1\circ{\rm proj},s,t)} (f_2\circ{\rm proj})\doteq f_2'=({\mu}_2',\psi_2',\phi_2'),$$
thus $(s,t)$ is a quadratic $f_1'$-derivation connecting $f_1'=(f_1 \circ {\rm proj})$ and $f_2'=(f_2\circ {\rm proj})$.
Therefore $\phi_1'([g])=\phi_1(g)$ and $\phi_2'([g])=\phi_2(g)$, for each $g \in G$. 
Also $\phi_1'([g,h])=\phi_2'([g,h])=1_G$ and $\psi_1'(1,[g,h])=\psi_2'(1,[g,h])=1_E$, for each $g,h \in G$.  Following the notation of lemma \ref{key}, let us put $(g,h)=(1,[g,h]) \in  E {{}_\d \times p} \FG(G)$, thus $\partial'(g,h)=[g,h]=[gh]^{-1}[g][h]$.

Let us look at the associated group map  (lemma \ref{esc}): $$(\phi_1',s)\colon \FG(G) \to G'\ltimes_{\t} E'.$$
Note the following equation, which will be used several times (we use \eqref{deri} and $\phi_1'([gh])=\phi_1(gh)$):
\begin{equation}
\begin{split}
(s \circ  \partial')(g,h)&= s([gh]^{-1}[g][h]))= \phi_1(gh)^{-1} \t s([gh]^{-1})\,\, \phi_1(h)^{-1} \t s([g])\,\, s([h]) \\&=s([gh])^{-1}\,\, \phi_1(h)^{-1} \t s([g])\,\, s([h]).
\end{split}
\end{equation}
The fact that $\f_2'(\ker(p))=1$ and $\f_1'(\ker(p))=1$  tells us  that for each $g,h \in G$:
\begin{align*}
 1&=\f_2'([g,h])=\f_2'([gh]^{-1}[g][h])=\f_1'([gh]^{-1}[g][h])\,\, \d(s([gh]^{-1}[g][h]))=\d\big(s([gh]^{-1}[g][h])\big)\\&=\d\big( s([gh])^{-1}\,\, \phi_1(h)^{-1} \t s([g])\,\, s([h])\big).
\end{align*}
Thus we must have that,  for each $g,h \in G$:
\begin{equation}\label{lax1}
\d \big(s([gh])\big)=\d\big(\phi_1(h)^{-1} \t s([g])\,\, s([h])\big) .
\end{equation}

Let us now look at the group map  (lemma \ref{esc}): $$(\psi_1',s\circ \d', t)\colon  E {{}_\d \times p} \FG(G)  \to E'\ltimes_{*} (E' \ltimes_{\t'} L').$$
Let $$\Pi(g,h)= t(g,h), \textrm{ where } g,h \in G.$$
Noting that   $\psi_1'(g,h)=1$, we have: 
\begin{align*}
 (\psi_1',s\circ \d', t) (g,h)&= \big(1,(s\circ \d')(g,h), t(g,h)\big)=\big(1, s([gh])^{-1}\,\, \phi_1(h)^{-1} \t s([g])\,\, s([h]),  t(g,h)\big)\\
&=\big(1, s([gh])^{-1}\,\, \phi_1(h)^{-1} \t s([g])\,\, s([h]), \Pi(g,h)\big).
\end{align*}
The fact that $\psi'_2(g,h)=1$ tells us that:
\begin{align*}
1=(s\circ \d')(g,h)\,\, \de( \Pi(g,h)) = s([gh])^{-1}\,\, \phi_1(h)^{-1} \t s([g])\,\, s([h])\,\, \de( \Pi(g,h)).
\end{align*}
Therefore, for each $g,h \in G$, we have:
\begin{equation}\label{lax2}
 s([gh])=\phi_1(h)^{-1}\t  s([g])\,\, s([h])\,\,\de( \Pi(g,h)),
\end{equation}
thus also for $g,h \in G$:
\begin{equation}(\psi_1',s\circ \d', t)(g,h)=\big(1,\de(\Pi(g,h))^{-1},\Pi(g,h)\big),\end{equation}
and
\begin{equation}\label{oiu}
 s([g,h])=\de(\Pi(g,h))^{-1},
\end{equation}
thus in particular:
\begin{equation}\label{aaaa}
s([1])=s([1,1])=\de(\Pi(1,1))^{-1} . 
\end{equation}

Since $\big((\psi_1',s\circ \d', t),(\phi_1',s)\big)$ is, by lemma \ref{esc}, a pre-crossed module map into  $(\beta \colon E'\ltimes_{*} (E' \ltimes_{\t'} L') \to G'\ltimes_\t E', \bullet)$, we have:
\begin{align}
 (\psi_1',s\circ \d', t) (1,k \t [g,h])&=(\psi_1',s\circ \d', t) \big(  k \t (g,h)\big)=\big(\phi_1'(k),s(k)\big)\bullet (\psi_1',s\circ \d', t) (g,h) \nonumber\\&=\big(\phi_1(p(k)),s(k)) \bullet  \big (1, s([gh])^{-1}\,\, \phi_1(h)^{-1} \t s([g])\,\, s([h]), t(g,h))\big),\label{that}
\end{align}
where $k \in \FG(G)$ and $g,h \in G$. 
The fact that 
$\psi_2'(k \t (g,h))=1$ is equivalent {(by \eqref{that} and lemmas \ref{connect} and \ref{esc}) to:}
$$\phi_1(p(k))\d(s(k)) \t \Big( s([gh])^{-1}\,\, \phi_1(h)^{-1}\t s([g])\,\, s([h])\,\, \de( \Pi(g,h)) \Big)=1, $$
 therefore this fact is implied by \eqref{lax2}. 

Let us now find necessary and sufficient conditions for 
$$k\in \FG(G) \longmapsto \big(\phi_1(p(k),s(k)\big)\in G' \ltimes_\t E' $$
and
$$(e,k) \in  E {{}_\d \times p} \FG(G) \longmapsto  (\psi_1(q(e,k)),s(\d'(e,k)), t(e,k))  \in E'\ltimes_{*} (E' \ltimes_{\t'} L')  $$ 
to be a pre-crossed module morphism, yielding a lax homotopy between the strict 2-crossed module maps $f_1$ and $f_2$. As far as $\big(\phi_1(p(k),s(k)\big)$ is concerned, since $\FG(G)$ is free on the underlying set of $G$, there are no conditions on $s$ to add to \eqref{lax1}.

Given $(e,k) \in E {{}_\d \times p} \FG(G)$ we have 
$$(e,k)=\big(e,[\d(e)]\big)\,\,\big (1,[\d(e)]^{-1} k\big),$$
where clearly $[\d(e)]^{-1} k\in \ker(p)$.
By using this and lemma \ref{key}, we can easily see that the group $ E {{}_\d \times p} \FG(G)$ is the principal group of the  $\FG(G)$ pre-crossed module, formally generated by the symbols $[e]\doteq (e,[\d(e])$, with $e \in E$, and $(g,h) \doteq (1,(g,h))$, with $g,h \in G$, modulo the relations \eqref{conj}, below
(for $g,h \in G$, $e,f \in E$ and $l \in \FG(G)$), which (we leave the reader to verify this) do hold in  $E {{}_\d \times p} \FG(G)$:
\begin{equation}\label{conj}
\begin{split}
[e] [f]&= [ef]\,\,(\d(e),\d(f))\\
g \t [e]&=[g \t e ]\,\, (g,\d(e)g^{-1})\,\, (\d(e),g^{-1})\,\,(g,g^{-1})^{-1}\,\,(1,1)^{-1} \\
(gh,i)\,\,[i]^{-1}\t (g,h) &=(g,hi)\,\, (h,i)\\
\d(g,h) \t (l \t (g',h'))&=(g,h) \,\, (l \t (g',h') )\,\, (g,h)^{-1} \textrm{ or }  ([g,h]l) \t (g',h')=(g,h) \,\, (l \t (g',h') )\,\, (g,h)^{-1}.
 \end{split}
\end{equation}
Let $t(a)=t([a])$ and $s(g)=s([g])$, where $g \in G$ and $a \in E$. Note that $\psi_1'([e])=\psi_1(e)$ and $\psi_2'([e])=\psi_1(e)$. Then $t\colon E {{}_\d \times p} \FG(G) \to L'$ can be specified by $t(a)$ and $t(g,h)=\Pi(g,h)$, which must satisfy relations \eqref{conj}.  These translate into (in order of appearance):
$$\big(\psi_1(a),s(\d(a)),t(a)\big)\,\, \big(\psi_1(b),s(\d(b)),t(b)\big)=\big(\psi_1(ab),s(\d(ab)),t(ab)\big)\, \big(1,\de(\Pi(\d(a),\d(b))^{-1}, \Pi(\d(a),\d(b))\big) ,$$
which is the same as, for any $a,b \in G$, (by \eqref{gl}):
\begin{multline*}
\Big(\psi_1(a)\psi_1(b),  \phi_1(\d(b))^{-1} \t s(a) \,\, s(\d(b)),  \big((s \circ \d) (b) \big)^{-1} \t'  \left \{\psi_1(b)^{-1},\ s (\d (a))^{-1}\right \}^{-1}     \, \\\left(\big(\p_1(b)((s \circ \d))(b)\big)^{-1} \t' t(a) \right )\,\, t(b) \Big)=
\Big(\psi_1(ab), s(\d(ab))\,\, \de(\Pi(\d(a),\d(b))^{-1}, \Pi(\d(a),\d(b)) \,\, t(ab)\Big).
\end{multline*}
At the level of the first two components, equality always hold, given the calculations above. Therefore 
\begin{equation}\label{lax3}
 \Pi(\d(a),\d(b)) \,\, t(ab)= \big((s \circ \d) (b) \big)^{-1} \t' \left (   \left\{\psi_1(b)^{-1},  s (\d (a))^{-1}\right \}^{-1} \,\, \p_1(b)^{-1} \t'  t(a) \right )\,\, t(b).
\end{equation}
The second relation translates into (for any $g \in G$ and $a \in E$):
\begin{multline*}
(\phi(g),s(g))\bullet (\psi(a),s(\d(a)),t(a))=\big(\psi_1(g \t a),s(\d( g\t a)),t(g \t a)\big)\,
 \big(1,\de(\Pi(g,\d(a)g^{-1}))^{-1} ,\Pi(g,\d(a)g^{-1}))\big) \,\\
\big(1, \de(\Pi(\d(a),g^{-1}))^{-1},  \Pi(\d(a),g^{-1})\big)\,
 \big(1,\de(\Pi(g,g^{-1})),\Pi(g,g^{-1})^{-1}\big)\,
\big(1,\de(\Pi(1,1)),\Pi(1,1)^{-1}\big) 
\end{multline*}
which gives:
\begin{multline}\label{lax4}\f_1(g) \t \Big (s(g) s(\d(a))^{-1} \t' \Big \{\p_1 (a)^{-1},  s(g)^{-1})\Big\}^{-1} \Big) 	\,\, \f_1(g) \t \Big \{ s(g), s(\d(a))^{-1}\p_1(a)^{-1} \Big\} \,\,  \big(\f_1(g)(\d \circ s)(g)\big) \t t(a)\\=
\Pi(1,1)^{-1}\, \Pi(g,g^{-1})^{-1}\, \Pi(\d(a),g^{-1})  \,\Pi(g,\d(a)g^{-1})) \,t(g \t a).
\end{multline}
By using \eqref{m4} and \eqref{gact}, the third relation translates into the cocycle type identity (for each $g,h,i \in G$):
\begin{equation}\label{lax5}
 s(i) \t'  \big( \phi_1(i)^{-1} \t  \Pi(g,h)\big) \,\,\Pi(gh,i) =  \Pi(h,i) \,\, \Pi(g,hi).
\end{equation}
Let us see that the fourth relation is void in this case. Note that,  for $g,h \in G$ and $k \in \FG(G)$, we have:
$$\big(\psi_1',(s \circ \d'),t\big)( k \t (g,h))=\big(\phi_1(p(k)),s(k)\big) \bullet \big(\psi_1,(s \circ \d'),t\big)(g,h)=\big(\phi_1(p(k)),s(k)\big) \bullet  \big(1,\de(\Pi(g,h))^{-1},\Pi(g,h)\big)  .$$
Applying $\big(\psi_1',(s \circ \d),t\big)$ to the left hand side of the last relation of \eqref{conj} yields (since $\phi_1(p([g,h]))=1$ and by using \eqref{m4}):
\begin{align*}
\big(\phi_1(p(l)),s([g,h]l)\big) &\bullet  \big(1,\de(\Pi(g',h'))^{-1},\Pi(g',h')\big)\\
&=\phi_1(p_1(l)) \t  \big(1,s([g,h]l) \,\, \de(\Pi(g',h'))^{-1}\,\, s([g,h]l)^{-1}, s([g,h]l) \t' \Pi(g',h')\big).
\end{align*}
Applying  $\big(\psi_1',(s \circ \d),t\big)$  to the right-hand-side of the last relation of \eqref{conj} gives:
\begin{align*}
&\big(1,\de(\Pi(g,h))^{-1},\Pi(g,h)\big) \,\,\big(\phi_1(p(l)),s(l)\big) \bullet  \big(1,\de(\Pi(g',h'))^{-1},\Pi(g',h')\big) \,\, \big(1,\de(\Pi(g,h))^{-1},\Pi(g,h)\big)^{-1}\\
&=\big(1,\de(\Pi(g,h))^{-1},\Pi(g,h)\big) \,\,\phi_1(p(l))\t   \big(1,s(l) \,\, \de(\Pi(g',h'))^{-1}\,\,s(l)^{-1}, s(l) \t' \Pi(g',h')\big) \,\, \big(1,\de(\Pi(g,h))^{-1},\Pi(g,h)\big)^{-1}.
\end{align*}
Both sides applied by $\big(\psi_1',(s \circ \d),t\big)$ coincide since (note the fact that $(\de\colon L \to E,\t')$ is a crossed module):
\begin{align*}
\Pi(g,h)^{-1} \,\, \phi_1(p(l) \t (s(l) \t' \Pi(g',h'))\,\,  \Pi(g,h)&=\phi_1(p(l)) \t\Big ( \big(\de(\phi_1(p(l))^{-1} \t \Pi(g,h)^{-1}) \,\, s(l) \big) \t'\Pi(g',h')\Big)\\
&=\phi_1(p(l)) \t \Big(\big(\phi_1(p(l))^{-1} \t s([g,h]) \,\, s(l) \big) \t'\Pi(g',h')\Big)\\
&=\phi_1(p(l)) \t\Big( \big( s([g,h] l)  \big) \t'\Pi(g',h')\Big).
\end{align*}

We have (almost) proven:
\begin{Theorem}
Consider 2-crossed modules $\A=\left (L \ra{\de} E \ra{\d} G ,  \t,\{,\}\right)$ and 
$\A'=\left (L' \ra{\de} E' \ra{\d} G' ,  \t,\{,\}\right)$. 
Let $f_1=({\mu}_1,\psi_1,\phi_1)$ and $f_2=({\mu}_2,\psi_2,\phi_2)$ be 2-crossed module maps $\A \to \A'$. 
A lax homotopy connecting $f_1$ and $f_2$, definition \ref{laxmspace}, which we write as:
 $$f_1 \ra{(f_1,\hat s, \hat t,\Pi)} f_2 ,$$ is given by:
\begin{enumerate}
 \item A (set) map $\hat s\colon G \to E',$
 \item A (set) map $\hat  t\colon E \to L',$
 \item A (set) map $\Pi\colon G \times G \to L'.$
\end{enumerate}
These are to satisfy, for each $g,h \in G$ and $a,b \in E$, that:
\begin{equation}\label{1lax1}
\d(\hat s(gh))=\d\big(\phi_1(h)^{-1} \t \hat s(g)\,\, \hat s(h)\big), 
\end{equation}
\begin{equation}\label{1lax2}
 \hat s(gh)=\phi_1(h)^{-1}\t  \hat s(g)\,\,\hat  s(h)\,\,\de( \Pi(g,h)),
\end{equation}
\begin{equation}\label{1lax3}
 \Pi(\d(a),\d(b)) \,\, \hat t(ab)= \big(( \hat s \circ \d) (b) \big)^{-1} \t' \left ( \left \{\psi_1(b)^{-1}, \hat  s (\d (a))^{-1}\right \}^{-1} \,\, \p_1(b)^{-1} \t'\hat  t(a) \right )\,\,\hat t(b),
\end{equation}
\begin{multline}\label{1lax4}
\f_1(g) \t \Big (\hat s(g) \hat s(\d(a))^{-1} \t' \Big \{\p_1 (a)^{-1},  \hat s(g)^{-1}\Big\}^{-1} \Big)\,\,  \f_1(g) \t \Big \{\hat  s(g),\hat  s(\d(a))^{-1}\p_1(a)^{-1} \Big\}\,\,   \big(\f_1(g)(\d \circ\hat  s)(g)\big) \t\hat  t(a)\\=
\Pi(1,1)^{-1}\, \Pi(g,g^{-1})^{-1}\, \Pi(\d(a),g^{-1})  \,\Pi(g,\d(a)g^{-1})) \,\hat t(g \t a),
\end{multline}
\begin{equation}\label{1lax5}
  \hat s(i) \t'  \big( \phi_1(i)^{-1} \t  \Pi(g,h)\big) \,\,\Pi(gh,i) =  \Pi(h,i) \,\, \Pi(g,hi).
\end{equation}
And, moreover, if $g \in G$, $e \in E$ and $l \in L$:
\begin{equation}\label{lp}
\begin{split}
 \phi_2(g)&=\phi_1(g) \,\d(\hat s(g)),\\
 \psi_2(e)&=\psi_1(e) \,\hat  s(\d(e))\, \de(\hat t(e)),\\
 {\mu}_2(l)&=\mu_1(l) \,\,  \Pi(1,1)^{-1}  \,\, \hat{t}(\de(l)).\end{split}
\end{equation}
Moreover, the corresponding (strict) homotopy between strict 2-crossed module maps $Q^1(\A) \to \A'$, namely:
$$f_1'=f_1\circ {\rm proj} \ra{(f_1 \circ {\rm proj},{s},t)} f_2\circ {\rm proj}=f_2'$$
 is given by the $(f_1\circ {\rm 
proj})$-quadratic derivation $(s,t)$, where $s\colon \FG(G) \to E'$ is the unique $(\phi_1\circ p)$-derivation (definition \ref{derdef}) extending $\hat{s}\colon G \to E'$, and on the  group generators of  $E {{}_\d \times p} \FG(G)$ 
we have that $t(e,[\d(e)])=\hat{t}(e)$ and $t(g,h)=\Pi(g,h)$, where $e \in E$ and $g,h \in G$. Recall that $p\colon \FG(G) \to G$ is the obvious projection.
\end{Theorem}
{We thus have a lax analogue of the strict homotopy relation treated in subsection \ref{qder}.}\\
\begin{Proof}
We just need to check the last equation of \eqref{lp}. Note that if $l \in L$:
$$\de'(l)=(\de(l),\emptyset)=(\de(l),[1])\,(1,[1]^{-1})=[\de(l)]\,(1,1)^{-1}.$$
Therefore, by using \eqref{h3alt} and \eqref{aaaa}, and noting $\psi'([1])=1_{E'}$:
$$\hat{t}(\de'(l))=\hat{s}([1]^{-1})^{-1} \t' \hat{t}(\de(l))\,\,\hat{t}\big((1,1)^{-1}\big)=\hat{s}([1]^{-1})^{-1} \t' \big ( \hat{t}(\de(l))\,\, (\hat{t}(1,1))^{-1}\big).$$
Thus
$$\mu_2(l)=\mu_2'(l)=\mu_1(l) \,\, \hat{t}(\de'(l))=\mu_1(l) \,\,  \hat{s}([1]^{-1})^{-1} \t' \big ( \hat{t}(\de(l))\,\, \Pi(1,1)^{-1}\big).$$
And now note that $\hat{s}([1]^{-1})=\phi'_1([1])\t \hat s([1])^{-1}=\hat s([1])^{-1}$, since $\phi'_1([1])=1$, together with \eqref{aaaa}, and the second Peiffer condition in definition \ref{cm}.
\end{Proof}

Note that it follows by construction that $f_2=(\mu_2,\psi_2,\phi_2)$, defined in \eqref{lp}, is a 2-crossed module morphism $\A \to \A'$, if equations \eqref{1lax1} to \eqref{1lax5} are satisfied.  {This can easily be  proven directly.}

\subsubsection{Composition and inverses of lax homotopies}\label{j2}
We now freely use the notation of  \ref{comphom} and \ref{ih}, as well as definition \ref{laxmspace}.
\begin{Theorem}
Consider 2-crossed modules $\A=\left (L \ra{\de} E \ra{\d} G ,  \t,\{,\}\right)$ and 
$\A'=\left (L' \ra{\de} E' \ra{\d} G' ,  \t,\{,\}\right)$. 
Given lax homotopies of 2-crossed module maps $\A\to \A'$, say:  $$f=({\mu},\psi,\phi) \ra{(f,\hat s, \hat t,\Pi)} f'=({\mu}',\psi',\phi')   \ra{(f',  \hat {s'} , \hat{ t'},\Pi')}  f''=({\mu}'',\psi'',\phi'') , $$
the explicit form of their {concatenation}, denoted by $(f,\hat s, \hat t,\Pi) \hat{\otimes}(f',  \hat {s'} , \hat{ t'},\Pi') $, is:
  $$f \ra{(f,\hat s \hat \tn \hat {s'},\hat t \hat{\tn} \hat{ t'},\Pi \hat \tn \Pi')} f'' ,$$
where
\begin{equation}\label{Comp1}
 (\hat s\hat \tn \hat {s'})(g)=\hat s(g) \hat {s'}(g), \fo g \in G,
\end{equation}
\begin{equation}\label{Comp2}
 (\hat t \hat\tn \hat{ t'})(g)=\hat{s'}(\d(e))^{-1}\t' \hat t(e)\,\, \hat{ t'}(e), \fo e \in E,
\end{equation}and where
\begin{equation}\label{Comp3}
 ( \Pi \hat \tn  \Pi')(g,h)=\Theta^{(\hat s,\hat {s'})}\big([gh],[g],[h]\big) \,\,\Pi'(g,h) \,\, \Pi(g,h),
\end{equation}
where $\Theta^{(\hat s,\hat {s'})}$ was defined in \eqref{deftheta}.

Moreover, the inverse of $f=({\mu},\psi,\phi) \ra{(f,\hat s, \hat t,\Pi)} f'=({\mu}',\psi',\phi')$ is  
$$ f'=({\mu}',\psi',\phi')\ra{(f',\overline {\hat s}, \overline {\hat t},\overline \Pi)}f=({\mu},\psi,\phi) , $$
where {if $g,h \in G$ and $e \in E$:}
\begin{equation}
\begin{split}
 \overline{\hat{s}}(g)&=\hat{s}(g)^{-1},\\
\overline{\hat{t}}(e)&=\hat{s}(\d(e)) \t' \hat{t}(e)^{-1},\\
\overline{\Pi}(g,h)&=\Theta^{(\hat s, \overline{\hat{s}})}\big([gh],[g],[h] \big)^{-1} \,\, \Pi(g,h)^{-1}.
\end{split}
\end{equation}

\end{Theorem}
\begin{Proof}
As far as the {concatenation} of lax homotopies is concerned, we just need to {consider} the corresponding chain of strict homotopies, given by the previous theorem:
$$f\circ {\rm proj}\ra{(f\circ {\rm proj}, s, t)} f'\circ {\rm proj}\ra{(f'\circ {\rm proj_1}, s', t')}  f''\circ {\rm proj},$$
and look at the construction of their concatenation in \ref{comphom}, noting that the underlying set of $G$ is a free (chosen) basis of $\FG(G)$. 

We know that $s([g])=\hat{s}(g)$,  $s'([g])=\hat{s'}(g)$ and $(s\tn s')([g])=s([g])s'([g])$ for each $g \in G$. Thus $(\hat{s} \hat{\tn} \hat{s'})(g)=(s\tn s')([g])=\hat{s}(g)\hat{s'}(g)$ for each $g \in G$. 

Analogously, if $e \in E$, then $\hat{t}(e)=t(e,[\d(e)])$,  $\hat{t'}(e)=t'(e,[\d(e)])$  and  $(\hat{t} \hat{\tn }\hat{t'}) (e)= (t \tn {t'}) (e,[\d(e)])$, and
\begin{align*}
(t\tn t')(e,[\d(e)])&=\w^{(s,s')}(\d'(e,[\d(e)]))\,\,s'(\d'(e,[\d(e)]))^{-1} \t' t((e,[\d(e)])\,\, t'((e,[\d(e)]))\\
&=\w^{(s,s')}([\d(e)]))\,\,s'([\d(e)]))^{-1} \t' t((e,[\d(e)])\,\, t'((e,[\d(e)]))\\
&= s'([\d(e)]))^{-1} \t' t((e,[\d(e)])\,\, t'((e,[\d(e)]))=\hat{s'}(\d(e))^{-1}\t' \hat t(e)
\,\, \hat{ t'}(e),
\end{align*}
where we used remark \ref{propomega}. Thus \eqref{Comp2} follows. 

Recall that, given $g,h \in G$, then $\Pi(g,h)=t(1,[g,h])$ and $\Pi'(g,h)=t'(1,[g,h])$. We have:
\begin{align*}
&(\Pi\hat{\tn} \Pi')(g,h)=(t\tn t')(1,[g,h])=\w^{(s,s')}(\d'(1,[g,h] ))\,\,s'(\d'(1,[g,h]))^{-1} \t' t(1,[g,h])\,\, t'(1,[g,h])\\&=\w^{(s,s')}([gh]^{-1}[g][h]))\,\,s'([gh]^{-1}[g][h])^{-1} \t' t([1,[g,h])\,\, t'(1,[g,h])
\\&=\w^{(s,s')}([gh]^{-1}[g][h]))\,\,\de(\Pi'(g,h))\t' t([1,[g,h])\,\, t'(1,[g,h])
\\&=\w^{(s,s')}([gh]^{-1}[g][h]))\,\, \Pi'(g,h) \,\,\Pi(g,h)
\\&=\Theta^{(s,s')}\big([gh]^{-1},[g],[h]\big)\,\, \Pi'(g,h) \,\,\Pi(g,h)=\Theta^{(\hat s, \hat {s'})}\big([gh]^{-1},[g],[h]\big)\,\, \Pi'(g,h) \,\,\Pi(g,h),
\end{align*}
where we  used remark \ref{propomega} again, and noted $\hat{s}(g)=s([g])$ and $\hat{s'}(g)=s'([g])$, for each $g \in G$. 

 Inverses are handled in the same way. For instance  if $g,h \in G$ (we use \eqref{oiu}):
\begin{align*}
 \overline{\Pi}(g,h)&=\overline{t}(1,[g,h])=(\w^{( s, \overline{{s}})}(\d'(1,[g,h]))^{-1}\,\, s(\d'(1,[g,h]))\t' t(1,[g,h])^{-1}\\
&=\w^{( s, \overline{{s}})}([gh]^{-1}[g][h])^{-1}\,\,s([gh]^{-1}[g][h]) \t' \Pi(g,h)^{-1}\\
&=\Theta^{( s, \overline{{s}})}\big([gh],[g],[h] \big)^{-1} \de(\Pi(g,h))^{-1}\t' \Pi(g,h)^{-1}\\
&=\Theta^{( s, \overline{{s}})}\big([gh],[g],[h] \big)^{-1} \,\, \Pi(g,h)^{-1}=\Theta^{( \hat s, \overline{\hat{s}})}\big([gh],[g],[h] \big)^{-1} \,\, \Pi(g,h)^{-1}.
\end{align*}
\end{Proof}

\subsubsection{Lax 2-fold homotopy}\label{j3}
We now discuss lax 2-fold homotopy.
Consider 2-crossed modules $\A=\left (L \ra{\de} E \ra{\d} G ,  \t,\{,\}\right)$ and 
$\A'=\left (L' \ra{\de} E' \ra{\d} G' ,  \t,\{,\}\right)$. 
Given two lax homotopies between the 2-crossed module maps $f,f'\colon \A \to \A'$, say: $$\xymatrix{ &f\ar@/^1pc/[rr]^{(f,\hat {s'},\hat{ t'},\Pi')}\ar@/_1pc/[rr]_{(f,\hat s,\hat t,\Pi)} & &f'},$$
a lax 2-fold homotopy connecting them, say: 
$$\xymatrix{ &f\ar@/^1pc/[rr]^{(f,\hat {s'},\hat{ t'},\Pi')}\ar@/_1pc/[rr]_{(f,\hat s,\hat t,\Pi)} & \Uparrow(f,\hat s,\hat t,\Pi, \hat k)  &f'},$$
is given by a map $\hat{k}\colon G \to L'$, without any restrictions, apart from that it should relate the two lax homotopies as in  \eqref{2d1}, \eqref{2d2} and \eqref{2d3}, below. This is because to  $\hat{k}$ we associate the unique (strict) quadratic $(f \circ {\rm proj}, s,t)$ 2-derivation  $k\colon \FG(G) \to L'$ such that $k([g])=\hat{k}(g)$ for each $g \in G$ (corollary \ref{extt}). Therefore:
\begin{equation}\label{2d1}
 \hat{s'}(g)=\hat{s}(g)\,\, \de(\hat k(g)), \textrm{ for each } g \in G,
\end{equation}
and for each $e \in E$:
\begin{equation}\label{2d2} \begin{split}
 \hat{t'}(e)&=t'((e, [\partial(e)])=k(\d'(e, [\partial(e)]))^{-1}\,\,  t((e, [\partial(e)]))=
k([\d(e)])^{-1} \,\, t((e, [\partial(e)])\\
            &=\hat{k}(\d(e))^{-1}\,\, \hat{t}(e).
\end{split}
\end{equation}
Also, by using equation \eqref{defxi} in lemma \ref{dereprop}, for each $g,h \in G$:
\begin{equation}\label{2d3}
\begin{split}
 \Pi'(g,h)&=t'((1,[g,h]))=k(\d'((1,[g,h])))^{-1}\,\, t((1,[g,h]))=k([gh]^{-1}\,\,[g][h])^{-1}\,\, \Pi(g,h)\\&=\left(\Xi^{(\phi,s,k)}([gh],[g],[h])\right)^{-1} \,\,\Pi(g,h).
\end{split}
\end{equation}
Thus, by using equation \eqref{defxi}:
$$\Pi'(1,1)= \left(\Xi^{(\phi,s,k)}([1],[1],[1])\right)^{-1} \,\,\Pi(1,1)=k([1])^{-1} \,\,\Pi(1,1)=\hat{k}(1)^{-1} \,\, \Pi(1,1).  $$
In particular, we can prove that, if \eqref{2d1}, \eqref{2d2} and \eqref{2d3} are satisfied, then if $ f\ra{(f,\hat s,\hat{t},\Pi)} f'$, we must also have that  $ f\ra{(f,\hat {s'},\hat{t'},\Pi')} f'$, by using  \eqref{lp}.

If we have a chain of 2-fold homotopies  $(f,\hat s,\hat t,\Pi) \ra{(f,\hat s,\hat t,\Pi,  \hat k)} (f,\hat {s'},\hat{ t'},\Pi') \ra{(f,\hat {s'},\hat{ t'},\Pi', \hat{ k'})} (f,\hat{ s''},\hat{ t''},\Pi''),$
diagrammatically:
 $$\xymatrix{ &f\ar[rrr]|{(f,\hat {s'}, \hat{ t'},\Pi')}\ar@/_3pc/[rrr]_{(f,\hat s, \hat t,\Pi)} \ar@/^3pc/[rrr]^{(f,\hat{ s''},\hat{ t''},\Pi'')} \ar@/_1pc/ @{{}{ }{}}[rrr]|{\Uparrow(f,\hat s, \hat t,\Pi, \hat k)}  \ar@/^1pc/ @{{}{ }{}}[rrr]|{\Uparrow(f,\hat {s'},\hat{ t'},\Pi',\hat {k'})}  &&&f'}$$
Then their vertical 
{concatenation} is given by the map $\hat k\hat \diamond \hat {k'} \colon G \to L$ such that: 
 $$(\hat k\diamond \hat {k'})(g)=\hat k(g) \hat {k'}(g),$$
for each $g \in G$.
By construction,  $\hat k \hat \diamond \hat {k'}$ does connect $(f,\hat s,\hat t,\Pi)$ and $(f,\hat{ s''},\hat{ t''},\Pi'')$.

Suppose that  we have a lax 2-fold   homotopy  say $(f,\hat s,\hat t,\Pi) \ra{(f,\hat s,\hat t, \Pi, \hat k)} (f,\hat {s'},\hat{ t'},\Pi'),$ which we  write as:
 $$\xymatrix{ &f\ar@/^1pc/[rr]^{(f,\hat {s'},\hat{ t'},\Pi)}\ar@/_1pc/[rr]_{(f,\hat s,\hat t,\Pi)} & \Uparrow(f,\hat s,\hat t, \Pi, \hat k)  &f'},$$
and that we also have a lax homotopy: $$f'=({\mu}',\psi',\phi') \ra{(f',\hat{ s''},\hat{ t''},\Pi'')} f''=({\mu}'',\psi'',\phi''),$$
so what we have diagrammatically is:
 $$\xymatrix{ &f\ar@/^1pc/[rr]^{(f,\hat {s'},\hat{ t'},\Pi')}\ar@/_1pc/[rr]_{(f,\hat s,\hat t,\Pi)} & \Uparrow(f,\hat s,\hat t,\Pi, \hat k)  &f' \ar[rr]^	{(f',\hat{ s''},\hat{ t''},\Pi'')} & &f''}.$$
The whiskering:
$$(f,\hat  s,\hat t, \Pi, \hat k){\hat \otimes} (f',\hat{ s''},\hat{ t''}, \Pi'') =(f,\hat s{\hat \otimes} \hat{ s''},\hat t{\hat\otimes} \hat t '',{\Pi}  \hat{\tn} \Pi'', k {\hat \otimes} s'')$$
is given by the map  $\hat k {\hat \otimes} \hat{ s''}\colon G \to L'$, which has the form, for each $g \in G$:
$$ (\hat k {\hat \otimes} \hat{ s''})(g)=\hat{ s''}(g)^{-1} \t' \hat k(g). $$
By construction we have:
 $$\xymatrix{ &f\ar@/^2pc/[rr]^{(f,\hat {s'} \hat{\tn} \hat{ s''},\hat{ t'}  \hat{\tn} \hat{ t''},\Pi'  \hat{\tn} \Pi'')}\ar@/_2pc/[rr]_{(f,\hat s \hat{\tn} \hat{ s''},\hat t  \hat{\tn} \hat{ t''},\Pi  \hat{\tn} \Pi'')   } & \Uparrow (f,\hat s{\hat \otimes} \hat{ s''},\hat t{\hat\otimes} \hat{ t''},{\Pi}  \hat{\tn} \Pi'', \hat k {\hat \otimes} \hat{s''})  & f'}.$$

Similarly, suppose that  we have 2-crossed module maps $f,f'\colon \A' \to \A$, lax homotopies $(f,\hat s,\hat t,\Pi )$ and $(f,\hat {s'},\hat{ t'},\Pi')$, connecting $f$ and $f'$,  a lax  2-fold homotopy $(f,\hat s,\hat t,\Pi, \hat {k'})$, connecting them, and also  a  lax  homotopy $f''\ra{(f'',\hat {s''},\hat  t'',\Pi'')}f,$
diagrammatically:
 $$\xymatrix{ &f'' \ar[rr]^{(f'',\hat{ s''},\hat{ t''},\Pi'')} & &f\ar@/^1pc/[rr]^{(f,\hat {s'},\hat  t',\Pi')}\ar@/_1pc/[rr]_{(f,\hat s,\hat  t,\Pi)} & \Uparrow(f,\hat s,\hat t,\Pi,\hat k) & f' }.$$
The  whiskering: $$(f'',\hat {s''},\hat  t'',\Pi'')\hat{\tn} (f,\hat s,\hat t,\Pi, \hat {k'})=
(f'',\hat{ s''}{\hat \otimes}\hat  s,\hat{ t''}{\hat \otimes} \hat t, \Pi'' \hat{\otimes} \Pi, \hat{ s''}{\otimes} \hat k),$$
is such that, for each $g \in G$ we have  $(\hat{ s''} {\hat \otimes}  \hat k)(g)=\hat k(g). $ 
By construction:
 $$\xymatrix{ &f''\ar@/^2pc/[rr]^{(f'',s''{\otimes} s',t''{\otimes} t'')}\ar@/_2pc/[rr]_{(f'',\hat{ s''}{\hat \otimes}\hat  s,\hat{ t''}{\hat \otimes} \hat t, \Pi'' \hat{\otimes} \Pi) } &  \Uparrow (f'', s''{\hat \otimes}\hat  s,\hat{ t''}{\hat \otimes} \hat t, \Pi'' \hat{\otimes} \Pi, \hat{ s''}{\otimes} \hat k)  &f'}.$$

By definition (since this is simply an unpacked version of definition \ref{laxmspace}), it follows:
\begin{Theorem}
 Let $\A$ and $\A'$ be 2-crossed modules. There exists a 2-groupoid $\EuScript{HOM}_{\rm LAX}(\A,\A')$, whose objects are the (strict) 2-crossed module maps $\A \to \A'$, the morphisms are the (pointed) lax homotopies between  2-crossed module maps, and the 2-morphisms are the (pointed) lax 2-fold homotopies between lax homotopies, whose explicit descriptions, and various {concatenation}s, are described in \ref{j1}, \ref{j2} and \ref{j3}.
\end{Theorem}

\subsection{Composition of lax homotopies with strict 2-crossed module maps}\label{co}

\begin{Theorem}\label{LaxWisk1}
 Let $\A$, $\A'$ and $\A''$ be 2-crossed modules. Let $f,f' \colon \A \to \A'$ be  2-crossed module maps. Let also $h=(\mu,\psi,\phi)\colon \A' \to \A'' $ be another 2-crossed module map. 
If we have a lax 2-crossed module homotopy $(f,\hat{s},\hat{t},\Pi)$ connecting $f$ and $f'$, then $$(h \circ f, \psi \circ \hat{s}, \mu \circ \hat{t},\mu \circ \Pi)\doteq h \circ (f,\hat{s},\hat{t},\Pi) $$ is a lax homotopy connecting $g \circ f$ and $g \circ f'$.
\end{Theorem}
\begin{Proof}
 Equations \eqref{1lax1} to \eqref{lp} are satisfied since $h$ preserves all 2-crossed module operations, strictly.
\end{Proof}

\begin{Theorem}\label{LaxWisk2}
 Let $\A$, $\A'$ and $\A''$ be 2-crossed modules. Let $f,f' \colon \A \to \A'$ be  2-crossed module maps. Let also $h'=(\mu',\psi',\phi')\colon \A'' \to \A $ be a 2-crossed module morphism. 
If we have a lax 2-crossed module homotopy $(f,\hat{s},\hat{t},\Pi)$ connecting $f$ and $f'$. Then $$(f \circ h', \hat{s} \circ \phi',  \hat{t} \circ \psi', \Pi \circ (\phi' \times \phi'))\doteq (f,\hat{s},\hat{t},\Pi)\circ h' $$ is a lax homotopy connecting $ f \circ h'$ and $f'\circ h''$.
\end{Theorem}

The operators defined in theorems \ref{LaxWisk1} and \ref{LaxWisk2} will be called composition operators.
\begin{Theorem}
 The composition operators preserve {concatenation}s and inverses of lax homotopies. 
\end{Theorem}
 \begin{Proof}
  Immediate from the explicit form of the {concatenation}s and inverses of lax homotopies, and the fact that we only compose homotopies with strict 2-crossed module morphisms.
 \end{Proof}
  
We thus have a sesquicategory \cite{St}, whose objects are the 2-crossed modules, the morphisms are the 2-crossed module maps, and the 2-morphisms are the lax homotopies between then. It is important to note that this is not a 2-category, since the interchange law does not hold in general.

The composition operators are also defined, in the obvious way, and with the obvious properties, for lax 2-fold homotopies and strict 2-crossed module maps. Given that for any two 2-crossed modules $\A$ and $\A'$ we have a 2-groupoid $\EuScript{HOM}_{\rm LAX}(\A,\A')$,  we expect that this will give a Gray category \cite{Cr,GPS,Gu} of 2-crossed modules, (strict) 2-crossed module maps, lax homotopies and lax 2-fold homotopies. 
\subsection{Lax homotopy equivalence of 2-crossed modules}
In this subsection we make use of subsection \ref{co}.
Let $\A$ and $\A'$ be 2-crossed modules.
 \begin{Definition}
 We say that $f\colon \A \to \A'$ is a lax homotopy equivalence if there exists a  2-crossed module map  $g\colon \A' \to \A$, and lax  homotopies: $$\id_\A \ra{(\id_\A,\hat s, \hat t,\Pi )} g\circ f \an \id_{\A'} \ra{(\id_{\A'},\hat u, \hat v,J)} f\circ g. $$
In such a case $g$ is said to be a homotopy inverse of $f$.
\end{Definition}

\begin{Lemma}
 The composition of lax homotopy equivalences is a lax homotopy equivalence.
\end{Lemma}
\begin{Proof}
This result is almost immediate from the fact that we can concatenate lax homotopies between strict 2-crossed module maps. Let us give details. 
Let $\A,\A'$ and $\A''$ be 2-crossed modules. Suppose that $f\colon \A \to \A'$ and $f'\colon \A' \to \A''$ are lax homotopy equivalences. Let us see that $f'\circ f$ is a lax homotopy equivalence. Choose inverses up to homotopy $g$ and $g'$ of $f$ and $g'$, respectively. We thus have lax homotopies
$$\id_\A \ra{(\id_\A,\hat s, \hat t,\Pi )} g\circ f, \quad \id_{\A'} \ra{(\id_{\A'},\hat u, \hat v,J)} f\circ g ,    \quad  \id_{\A'} \ra{(\id_{\A'},\hat {s'}, \hat{ t'},\Pi' )} g'\circ f', \quad \id_{\A''} \ra{(\id_{\A''},\hat u', \hat v',J')} f'\circ g'  .$$

{We prove that $g \circ g'$ is a lax homotopy inverse of $f'\circ f$. This follows by considering the  {concatenations} below of  lax homotopies}: $$\id_\A \ra{(\id_\A,\hat s, \hat t,\Pi )} g\circ f=g \circ \id_{\A'} \circ f \ra{ g \circ (\id_{\A'},\hat {s'}, \hat{ t'},\Pi' )\circ f} g\circ g'\circ f' \circ f, $$
$$ \id_{\A''} \ra{(\id_{\A''},\hat u', \hat v',J')} f'\circ g'=f'\circ\id_{\A'} \circ g' \ra{f' \circ (\id_{\A'},\hat u, \hat v,J)\circ g' } f' \circ f\circ g \circ g'.$$

\end{Proof}

\begin{Proposition}
 Lax homotopy equivalences of 2-crossed modules have the two-of-three property \cite{DS}.
\end{Proposition}
\begin{Proof}
The complete proof is analogous to the proof of the particular case above.
\end{Proof}\\
By using the composition operators, we can easily see that a retract of a lax homotopy equivalence is again a lax homotopy equivalence.

\subsubsection*{Acknowledgements}

\noindent J. Faria Martins was partially supported by CMA/FCT/UNL, under the grant PEst-OE/MAT/UI0297/2011. B. Gohla was supported by FCT (Portugal) through the doctoral grant SFRH/BD/33368/ 2008. 
This work was  supported by FCT, by means of  the projects
PTDC/MAT/098770/2008 
and
PTDC/MAT/101503/2008. 
{We are grateful to the referee for the comments and suggestions made.}

\end{document}